\documentclass{article}

\usepackage{microtype}
\usepackage{graphicx}
\usepackage{subfigure}
\usepackage{booktabs} 
\usepackage{amsmath, amssymb, amsfonts, listings, mathrsfs}

\usepackage{hyperref}

\usepackage[accepted]{icml2020_modified}

\newcommand{\nr}{\par \noindent}
\newcommand{\Def}{\stackrel{\mathrm{def}}{=}}

\newcommand{\dom}{{\rm dom \,}}
\newcommand{\beq}{\begin{equation}}
\newcommand{\eeq}{\end{equation}}

\newcommand{\R}{\mathbb{R}}

\newcommand{\E}{\mathbb{E}}

\newcommand{\refLE}[1]{\ensuremath{\stackrel{(\ref{#1})}{\leq}}}
\newcommand{\refEQ}[1]{\ensuremath{\stackrel{(\ref{#1})}{=}}}
\newcommand{\refGE}[1]{\ensuremath{\stackrel{(\ref{#1})}{\geq}}}

\newtheorem{assumption}{Assumption}
\newtheorem{corollary}{Corollary}

\newtheorem{definition}{Definition}
\newtheorem{example}{Example}
\newtheorem{lemma}{Lemma}

\newtheorem{remark}{Remark}
\newtheorem{theorem}{Theorem}
\newcommand{\proof}{\bf Proof: \rm \nr}
\newcommand{\qed}{\hfill $\Box$ \nr \medskip}

\newcommand{\ba}{\begin{array}}
\newcommand{\ea}{\end{array}}
\newcommand{\beann}{\begin{eqnarray*}}
\newcommand{\eeann}{\end{eqnarray*}}
\newcommand{\bea}{\begin{eqnarray}}
\newcommand{\eea}{\end{eqnarray}}

\newcommand{\BT}{\begin{theorem}}
\newcommand{\ET}{\end{theorem}}
\newcommand{\BL}{\begin{lemma}}
\newcommand{\EL}{\end{lemma}}
\newcommand{\BC}{\begin{corollary}}
\newcommand{\EC}{\end{corollary}}
\newcommand{\BE}{\begin{example}}
\newcommand{\EE}{\end{example}}
\newcommand{\BD}{\begin{definition}}
\newcommand{\ED}{\end{definition}}
\newcommand{\BR}{\begin{remark}}
\newcommand{\ER}{\end{remark}}
\newcommand{\BAS}{\begin{assumption}}
\newcommand{\EAS}{\end{assumption}}
\newcommand{\BI}{\begin{itemize}}
\newcommand{\EI}{\end{itemize}}

\newcommand{\BMP}{\begin{minipage}{9.5cm}}
\newcommand{\EMP}{\end{minipage}}
\newcommand{\MPT}{\begin{minipage}{11.5cm}}
\newcommand{\EPT}{\end{minipage}}

\newcommand{\la}{\langle}
\newcommand{\ra}{\rangle}

\def\argmin{\mathop{\rm argmin}}
\def\Argmin{\mathop{\rm Argmin}}
\def\R{\mathbb{R}}

\def\E{\mathbb{E}}

\icmltitlerunning{Inexact Tensor Methods with Dynamic Accuracies}

\begin{document}
	
\twocolumn[
\icmltitle{Inexact Tensor Methods with Dynamic Accuracies}


\icmlsetsymbol{equal}{*}

\begin{icmlauthorlist}
	\icmlauthor{Nikita Doikov}{icteam}
	\icmlauthor{Yurii Nesterov}{core}
\end{icmlauthorlist}

\icmlaffiliation{icteam}{Institute of Information and Communication Technologies,
	Electronics and Applied Mathematics
	(ICTEAM), Catholic University of Louvain (UCL), Belgium}
\icmlaffiliation{core}{Center for Operations Research and Econometrics (CORE),
	Catholic University of Louvain (UCL), Belgium}

\icmlcorrespondingauthor{Nikita Doikov}{nikita.doikov@uclouvain.be}
\icmlcorrespondingauthor{Yurii Nesterov}{yurii.nesterov@uclouvain.be}

\icmlkeywords{Machine Learning, ICML}

\vskip 0.3in
]



\printAffiliationsAndNotice{}  

\begin{abstract}
	In this paper, we study inexact high-order Tensor Methods
	for solving convex optimization problems with composite
	objective. At every step of such methods, we use
	approximate solution of the auxiliary problem, defined by
	the bound for the residual in function value. We propose
	two dynamic strategies for choosing the inner accuracy:
	the first one is decreasing as $1/k^{p + 1}$, where $p
	\geq 1$ is the order of the method and $k$ is the
	iteration counter, and the second approach is using for
	the inner accuracy the last progress in the target
	objective. We show that inexact Tensor Methods with these
	strategies achieve the same global convergence rate as in
	the error-free case. For the second approach we also
	establish local superlinear rates (for $p \geq 2$), and
	propose the accelerated scheme. Lastly, we present
	computational results on a variety of machine learning
	problems for several methods and different accuracy
	policies.
\end{abstract}

\section{Introduction}
\label{sec:Intro}

\subsection{Motivation}

With the growth of computing power, high-order
optimization methods are becoming more and more popular in
machine learning, due to their ability to tackle the
ill-conditioning and to improve the rate of convergence.
Based on the work~\cite{nesterov2006cubic}, where global
complexity guarantees for the Cubic regularization of
Newton method were established, a significant leap in the
development of second-order optimization algorithms was
made, discovering stochastic and randomized
methods~\cite{kohler2017sub,tripuraneni2018stochastic,cartis2018global,doikov2018randomized,wang2018stochastic,zhou2019stochastic},
which have better convergence rate, than the corresponding
first-order analogues. The main weakness, though, is that
every step of Newton method is much more expensive. It
requires to solve the subproblem, which is a minimization
of quadratic function with a regularizer, and possibly
with  some additional nondifferentiable components.
Therefore, the idea to employ \textit{higher} derivatives
into optimization schemes was remaining questionable,
because of the high cost of the computations. However,
recently~\cite{nesterov2019implementable} it was shown,
that the third-order Tensor Method for convex minimization
problems admits very effective implementation, with the
cost, which is comparable to that of the Newton step.

Now we have a family of methods (starting from the methods
of order one), for each iteration of which may need to
call some auxiliary subsolver. Thus, it becomes important
to study: \textit{which level of exactness we need to
	ensure at the step for not loosing the fast convergence of
	the initial method}. In this work, we suggest to describe
approximate solution of the subproblem in terms of the
residual in function value. We propose two strategies for
the inner accuracies, which are \textit{dynamic} (changing
with iterations). Indeed, there is no need to have a very
precise solution of the subproblem at the first
iterations, but we reasonably ask for higher precision
closer to the end of the optimization process.

\subsection{Related Work}

Global convergence of the first-order methods with inexact
proximal-gradient steps was studied
in~\cite{schmidt2011convergence}. The authors considered
the errors in the residual in function value of the
subproblem, and require them to decrease with iterations
at an appropriate rate. This setting is the most similar
to the current work.

In~\cite{cartis2011adaptive1,cartis2011adaptive2}, adaptive second-order
methods with cubic regularization
and inexact steps were proposed.
High-order inexact tensor methods were considered in
\cite{birgin2017worst,jiang2018unified,grapiglia2019tensor,grapiglia2019tensor2,cartis2019universal,lucchi2019stochastic}.
In all of these works, the authors
describe approximate solution of the subproblem
in terms of the corresponding first-order optimality condition (using the gradients).
This can be difficult to achieve by the current optimization schemes,
since it is more often that we have a better (or the only)
guarantees for the decrease of the residual in function value.
The latter one is used as a measure of inaccuracy in the recent
work~\cite{nesterov2019inexact} on the inexact Basic Tensor Methods.
However, only the constant choice of the accuracy level is considered
there.

\subsection{Contributions}

We propose new dynamic strategies for choosing the inner accuracy
for the general Tensor Methods, and several inexact algorithms based on it,
with proven complexity guarantees,
summarized next (we denote by $\delta_k$ the required precision for the
residual in function value of the auxiliary problem):

\begin{itemize}
	\item The rule $\delta_k := 1 / k^{p + 1}$, where $p \geq 1$ is
	the order of the method, and $k$ is the iteration counter.
	
	Using this strategy, we propose two optimization schemes:
	Monotone Inexact Tensor Method I (Algorithm~\ref{alg:Monotone1})
	and Inexact Tensor Method with Averaging
	(Algorithm~\ref{alg:Averaging}). Both of them have the global complexity
	estimates $O({1 /
		\varepsilon}^{\frac{1}{p}})$
	iterations for minimizing the convex function
	up to $\varepsilon$ -accuracy (see Theorem~\ref{th:GlobalConvAlg1}
	and Theorem~\ref{th:GlobalConvAlg3}). The latter method seems to be the first
	\textit{primal} high-order  scheme (aggregating the points from the primal space only),
	having the explicit distance between the starting point
	and the solution, in the complexity bound.
	
	\item The rule $\delta_k := c \cdot (F(x_{k - 2}) - F(x_{k - 1}))$,
	where $F(x_i)$ are the values of the target objective during
	the iterations, and $c \geq 0$ is a constant.
	
	We incorporate this
	strategy into our Monotone Inexact Tensor Method II
	(Algorithm~\ref{alg:Monotone2}).
	For this scheme, for minimizing convex functions up
	to $\varepsilon$-accuracy by the methods of order $p \geq
	1$, we prove the global complexity
	proportional to $O(1 / \varepsilon^{\frac{1}{p}})$
	(Theorem~\ref{th:GlobalConvAlg2}).
	The global rate becomes linear, if the objective is uniformly convex
	(Theorem~\ref{th:GlobalStrongAlg2}).
	
	Assuming that $\delta_k := c \cdot (F(x_{k - 2}) - F(x_{k
		- 1}))^{p + 1 \over 2}$, for the methods of order $p \geq
	2$ as applied to minimization of strongly convex
	objective, we also establish the local superlinear rate of
	convergence (see Theorem~\ref{th:LocalAlg2}).
	
	\item Using the technique of Contracting Proximal
	iterations~\cite{doikov2019contracting}, we
	propose inexact Accelerated Scheme~(Algorithm~\ref{alg:Accelerated}),
	where at each iteration $k$, we solve the corresponding
	subproblem with the precision $\zeta_k := 1 / k^{p + 2}$
	in the residual of the function value, by inexact Tensor Methods
	of order $p \geq 1$. The resulting complexity bound is
	$\tilde{O}(1 / \varepsilon^{\frac{1}{p + 1}})$
	inexact tensor steps for minimizing the convex function
	up to $\varepsilon$ accuracy (Theorem~\ref{th:GlobalAlg4}).

	\item Numerical results with empirical study of the methods
	for different accuracy policies are provided.
	
\end{itemize}

\subsection{Contents}

The rest of the paper is organized as follows. Section~\ref{sec:Notation}
contains notation which we use throughout the paper, and
declares our problem of interest in the composite form.
In Section~\ref{sec:Tensor} we introduce high-order model of the objective
and describe a general optimization scheme using this model.
Then, we summarize some known techniques for computing a step
for the methods of different order.
In Section~\ref{subsec:Monotone}
we study monotone inexact methods, for which we guarantee
the decrease of the objective function for every iteration, and
in Section~\ref{subsec:Averaging} we study the methods with averaging.
In Section~\ref{sec:Acceleration} we present our accelerated scheme.
Section~\ref{sec:Experiments} contains numerical results.
Missing proofs are provided in the supplementary material.

\section{Notation}
\label{sec:Notation}

In what follows, we denote by $\E$ a finite-dimensional real vector
space and by $\E^{*}$ its dual space, which is a space of linear
functions on $\E$. The value of function $s \in \E^{*}$ on $x \in \E$
is denoted by $\la s, x \ra$.
One can always identify $\E$ and $\E^{*}$ with $\R^n$,
when some basis is fixed, but often it is useful to seperate these
spaces, in order to avoid ambiguities.

Let us fix some symmetric positive definite
linear operator $B: \E \to \E^{*}$
and use it to define Euclidean norm for the primal variables:
$
\|x\| \Def  \la Bx, x \ra^{1/2}, \; x \in \E.
$
Then, the norm for the dual space is defined as:
$$
\ba{rcl}
\|s\|_{*} \Def \max\limits_{h \in \E: \, \|h\| \leq 1} \la s, x \ra
\; = \; \la s, B^{-1}s \ra^{1/2}, \; s \in \E^{*}.
\ea
$$
For a smooth function $f$, its gradient at point $x$ is denoted
by $\nabla f(x)$, and its Hessian is $\nabla^2 f(x)$. Note that
for $x \in \dom f \subseteq \E$ we have
$\nabla f(x) \in \E^{*}$, and $\nabla^2 f(x)h \in \E^{*}$
for $h \in \E$.

For $p \geq 1$, we denote by
$D^p f(x)[h_1, \cdots, h_p]$
$p$th directional derivative of $f$ along
directions $h_1, \dots, h_p \in \E$. If $h_i = h$
for all $1 \leq i \leq p$, the shorter notation $D^p f(x)[h]^p$
is used. The norm of $D^p f(x)$, which is
$p$-linear symmetric form on $\E$, is induced in the standard way:
$$
\ba{rcl}
\| D^p f(x) \| & \Def & \max\limits_{\substack{h_1, \dots, h_p \in \E: \\[1pt]
		\|h_i\| \leq 1, \; 1 \leq i \leq p}}
D^p f(x)[h_1, \dots, h_p] \\
\\
& = & \max\limits_{h \in \E: \, \|h\| \leq 1} \big| D^p f(x)[h]^p \big|.
\ea
$$
See Appendix~1 in \cite{nesterov1994interior} for the proof of the last equation.

We are interested to solve convex optimization problem in the
composite form:
\beq \label{MainProblem}
\ba{rcl}
\min\limits_{x \in \E} \Bigl\{ F(x) & \equiv & f(x) + \psi(x) \Bigr\},
\ea
\eeq
where $\psi: \E \to \R \cup \{ +\infty \}$
is a \textit{simple} proper closed convex function,
and $f : \dom \psi \to \R$ is several times differentiable and convex.
Basic examples of $\psi$ which we keep in mind are:
$\{0, +\infty\}$-indicator of a simple closed convex set
and $\ell_1$-regularization.
We always assume, that solution $x^{*} \in \dom \psi$
of problem~\eqref{MainProblem} does exist, denoting $F^{*} = F(x^{*})$.

\section{Inexact  Tensor Methods}
\label{sec:Tensor}

\subsection{High-Order Model of the Objective}

We assume, that for some $p \geq 1$,
the $p$th derivative of the smooth component
of our objective~\eqref{MainProblem} is Lipschitz continuous.
\BAS \label{as:Lip}
For all $x, y \in \dom \psi$
\beq \label{LipDef}
\ba{rcl}
\| D^p f(x) - D^p f(y) \| & \leq & L_p \|x - y\|.
\ea
\eeq
\EAS

Examples of convex functions with known Lipschitz
constants are as follows.

\BE
For the power of Euclidean norm $f(x) = \frac{1}{p + 1}\|x - x_0\|^{p + 1}$, $p \geq 1$,
\eqref{LipDef} holds with $L_p = p!$
(see Theorem~7.1 in \cite{rodomanov2019smoothness}).
\EE

\BE
For a given $a_i \in \E^{*}, 1 \leq i \leq m$, consider the log-sum-exp function:
$$
\ba{rcl}
f(x) & = & \log \left( \sum\limits_{i = 1}^m e^{\la a_i, x \ra}  \right), \quad x \in \E.
\ea
$$
Then, for $B := \sum_{i = 1}^m a_i a_i^{*}$ (assuming $B \succ 0$,
otherwise we can reduce dimensionality of the problem), \eqref{LipDef} holds
with $L_1 = 1$, $L_2 = 2$ and $L_3 = 4$
(see Lemma~\ref{LemmaLogSumExp} 
in the supplementary material).
\EE

\BE
Using $\E = \R$, and $a_1 = 0, a_2 = 1$ in the previous example,
we obtain the logistic regression loss:
$f(x) = \log(1 + e^x)$.
\EE

Let us consider Taylor's model of $f$ around a fixed point
$x$:
$$
\ba{rcl}
f(y) \; \approx \; f_{p, x}(y) & \Def &
f(x) + \sum\limits_{k = 1}^p \frac{1}{k!} D^k f(x)[y - x]^k.
\ea
$$
Then, from~\eqref{LipDef} we have a global bound
for this approximation. It holds, for all $x, y \in \dom \psi$
\beq \label{LipFuncBound}
\ba{rcl}
|f(y) - f_{p, x}(y)| & \leq & \frac{L_p}{(p + 1)!} \|y - x\|^{p + 1}.
\ea
\eeq
Denote by $\Omega_{H}(x; y)$ the following
regularized model of our objective:
\beq \label{ModelDef}
\ba{rcl}
\Omega_{H}(x; y) & \Def & f_{p, x}(y) + \frac{H\|y - x\|^{p + 1}}{(p + 1)!} + \psi(y),
\ea
\eeq
which serves as the global upper bound: $F(y) \leq \Omega_{H}(x; y)$,
for $H$ big enough (at least, for $H \geq L_p$).
This property suggests us to use the minimizer of~\eqref{ModelDef}
in $y$, as the next point of a hypothetical optimization scheme,
while $x$ being equal to a current iterate:
\beq \label{GenScheme}
\ba{rcl}
x_{k + 1} & \in & \Argmin_y \Omega_{H}(x_k; y), \quad k \geq 0.
\ea
\eeq
The approach of using high-order Taylor model $f_{p, x}(y)$
with its regularization was investigated first in~\cite{baes2009estimate}.

Note, that for $p = 1$, iterations~\eqref{GenScheme}
gives the Gradient Method (see \cite{nesterov2013gradient}
as a modern reference),
and for $p = 2$ it corresponds to the
Newton method with Cubic regularization~\cite{nesterov2006cubic}
(see also~\cite{doikov2018randomized} and~\cite{grapiglia2019accelerated}
for extensions to the composite setting).

Recently it was shown in~\cite{nesterov2019implementable},
that for $H \geq pL_p$, function $\Omega_H(x; \cdot)$ is
\textit{always convex} (despite the Taylor's
polynomial~$f_{p, x}(y)$ is nonconvex for $p \geq 3$).
Thus, computation~\eqref{GenScheme}  of the next point can
be done by powerful tools of Convex Optimization. Let us
summarize some known techniques for computing a step of
the general method~\eqref{GenScheme}, for different $p$:
\begin{itemize}
	\item $p = 1$. When there is no composite part: $\psi(x) = 0$,
	iteration~\eqref{GenScheme} can be represented as the
	Gradient Step with \textit{preconditioning}:
	$x_{k + 1} = x_k - \frac{1}{H}B^{-1} \nabla f(x_k)$.
	One can precompute inverse of $B$ in advance, or
	use some numerical subroutine at every step,
	solving the linear system
	(for example, by Conjugate Gradient method,
	see~\cite{nocedal2006numerical}).
	If $\psi(x) \not= 0$, computing the corresponding
	prox-operator is required~(see~\cite{beck2017first}).
	
	\item $p = 2$. One approach consists in diagonalizing
	the quadratic part by using eigenvalue-
	or tridiagonal-decomposiition of the Hessian matrix.
	Then, computation of the model minimizer~\eqref{GenScheme}
	can be done very efficiently by solving
	some one-dimensional
	nonlinear equation~\cite{nesterov2006cubic,gould2010solving}.
	The usage of fast approximate eigenvalue computation
	was considered in~\cite{agarwal2017finding}.
	Another way to compute the (inexact) second-order
	step~\eqref{GenScheme} is to launch
	the Gradient Method~\cite{carmon2019gradient}.
	Recently, the subsolver based on the Fast Gradient Method
	with restarts was proposed in~\cite{nesterov2019inexact},
	which convergence rate is $O(\frac{1}{t^6})$, where $t$
	is the iteration counter.
	
	\item $p = 3$. In~\cite{nesterov2019implementable},
	an efficient method for computing the inexact
	third-order step~\eqref{GenScheme} was outlined,
	which is based on the notion of
	\textit{relative smoothness}~\cite{van2017forward,bauschke2016descent,lu2018relatively}.
	Thus, every iteration of the subsolver
	can be represented as the Gradient Step for the auxiliary problem,
	with a specific choice of prox-function, formed by the
	second derivative of the initial objective and augmented
	by the fourth power of Euclidean norm.
	The methods of this type were studied in~\cite{grapiglia2019inexact,bullins2018fast,nesterov2019inexact}.
	
\end{itemize}

Up to our knowledge, the efficient implementation of the
tensor step of degree $p \geq 4$ remains to be an open
question.

\subsection{Monotone Inexact Methods}
\label{subsec:Monotone}

Let us assume that at every step of our method, we
minimize the model~\eqref{ModelDef} inexactly by an
auxiliary subroutine, up to some given accuracy $\delta
\geq 0$. We use the following definition of
\textit{inexact $\delta$-step.}

\BD
\label{def:InexStep}
Denote by $T_{H, \delta}(x)$ a point
$T \equiv T_{H, \delta}(x) \in \dom \psi$, satisfying
\beq \label{InexStep}
\ba{rcl}
\Omega_{H}(x; T) - \min\limits_{y} \Omega_{H}(x; y) & \leq & \delta.
\ea
\eeq
\ED

The main property of this point is given by the next
lemma.
\BL \label{lemma:Global}
Let $H = \alpha L_p$ for some $\alpha \geq p$. Then,
for every $y \in \dom \psi$
\beq \label{OneStep}
\ba{rcl}
F(T_{H, \delta}(x)) & \leq & F(y) + \frac{(\alpha + 1)L_p \|y - x\|^{p + 1}}{(p + 1)!} + \delta.
\ea
\eeq
\EL
\proof
Indeed, denoting $T \equiv T_{H, \delta}(x)$, we have
$$
\ba{rcl}
F(T) & \refLE{LipFuncBound} & \Omega_{H}(x; T)
\; \refLE{InexStep} \; \Omega_{H}(x; y) + \delta, \\
\\
& \refLE{LipFuncBound} & F(y) + \frac{(\alpha + 1) L_p \|y - x\|^{p + 1}}{(p + 1)!} + \delta.
\ea
$$
\qed
Now, if we plug $y := x$ (a current iterate) into~\eqref{OneStep}, we obtain
$F(T_{H, \delta}(x)) \leq F(x) + \delta$.
So in the case $\delta = 0$ (exact tensor step),
we would have nonincreasing sequence~$\{ F(x_k) \}_{k \geq 0}$
of test points of the method. However, this is not the case for $\delta > 0$
(inexact tensor step).
Therefore we propose the following minimization scheme \textit{with correction}.

\begin{algorithm}[!h]
	\caption{Monotone Inexact Tensor Method, I}
	\label{alg:Monotone1}
	\begin{algorithmic}
		\STATE {\bfseries Initialization:} Choose $x_0 \in \dom \psi$. Fix $H := p L_p$.
		\FOR{$k = 0, 1, 2, \dots$}
		\STATE Pick up $\delta_{k + 1} \geq 0$
		\STATE Compute inexact tensor step $T_{k + 1} := T_{H, \delta_{k + 1}}(x_k)$
		\IF{$F(T_{k + 1}) < F(x_k)$}
		\STATE $x_{k + 1} := T_{k + 1}$
		\ELSE
		\STATE $x_{k + 1} := x_k$
		\ENDIF
		\ENDFOR
	\end{algorithmic}
\end{algorithm}

If at some step $k \geq 0$ of this algorithm we get $x_{k
	+ 1} = x_k$, then we need to decrease inner accuracy for
the next step. From the practical point of view, an
efficient implementation of this algorithm should include
a possibility of improving accuracy of the previously
computed point.

Denote by $D$ the radius of the initial level set of the objective:
\beq \label{eq-DRad}
\ba{rcl}
D & \Def & \sup\limits_{x}
\Bigl\{  \|x - x^{*}\| \; : \;  F(x) \leq F(x_0)  \Bigr \}.
\ea
\eeq

For Algorithm~\ref{alg:Monotone1}, we can prove the
following convergence result, which uses a simple strategy
for choosing $\delta_{k + 1}$.

\BT \label{th:GlobalConvAlg1}
Let $D < +\infty$. Let the sequence of inner accuracies
$\{ \delta_k \}_{k \geq 1}$ be chosen according to the
rule
\beq \label{DeltaKChoice}
\boxed{
	\ba{rcl}
	\delta_k & := & \frac{c}{k^{p + 1}}
	\ea
}
\eeq
with some $c \geq 0$. Then for the sequence $\{ x_k \}_{k
	\geq 1}$ produced by Algorithm~\ref{alg:Monotone1}, we
have
\beq \label{eq:GlobalConvAlg1}
\ba{rcl}
F(x_k) - F^{*} & \leq &
\frac{(p + 1)^{p + 1} L_p D^{p + 1}}{p! \, k^{p}} + \frac{c}{k^p}.
\ea
\eeq
\ET

We see, that the global convergence rate of the inexact
Tensor Method remains on the same level, as of the exact
one. Namely, in order to achieve $F(x_K) - F^{*} \leq
\varepsilon$, we need to perform $K = O({1 /
	\varepsilon^{1 \over p}})$ iterations of the algorithm.
According to these estimates, at the last iteration $K$,
the rule~\eqref{DeltaKChoice} requires to solve the
subproblem up to accuracy
\beq \label{DeltaLastIterChoice}
\ba{rcl}
\delta_K & = & O(c \varepsilon^{p + 1 \over p}).
\ea
\eeq
This is intriguing, since for bigger $p$ (order of the
method) we need less accurate solutions. Note, that the
estimate~\eqref{DeltaLastIterChoice} coincides with the
constant choice of inner accuracy
in~\cite{nesterov2019inexact}. However, the dynamic
strategy~\eqref{DeltaKChoice} provides a significant
decrease of the computational time on the first iterations
of the method, which is also confirmed by our numerical
results (see Section~\ref{sec:Experiments}).

Now, looking at Algorithm~\ref{alg:Monotone1},
one may think that we are forgetting the points $T_{k + 1}$ such
that $F(T_{k + 1}) \geq F(x_k)$, and thus we are loosing some computations.
However, this is not true: even if point $T_{k + 1}$ has not been taken as $x_{k + 1}$,
we shall use it internally as a starting point for computing the next $T_{k + 2}$.
To support this concept, we introduce \textit{inexact $\delta$-step}
with an additional condition of \textit{monotonicity}.

\BD
\label{def:InexStepMonotone}
Denote by $M_{H, \delta}(x)$ a point
$M \equiv M_{H, \delta}(x) \in \dom \psi$, satisfying
the following two conditions.
\begin{eqnarray}
& \; \Omega_H(x; M) - \min\limits_{y} \Omega_H(x; y)
\; \leq \; \delta, \label{eq-MonDef1} \\[7pt]
& \; \qquad F(M) \; < \; F(x). \label{eq-MonDef2}
\end{eqnarray}
\ED

It is clear, that point $M$ from Definition~\ref{def:InexStepMonotone}
satisfies Definition~\ref{def:InexStep} as well
(while the opposite is not always the case).
Therefore, we can also use Lemma~\ref{lemma:Global} for the
\textit{monotone} inexact tensor step.
Using this definition, we simplify Algorithm~\ref{alg:Monotone1}
and present the following scheme.

\begin{algorithm}[!h]
	\caption{Monotone Inexact Tensor Method, II}
	\label{alg:Monotone2}
	\begin{algorithmic}
		\STATE {\bfseries Initialization:} Choose $x_0 \in \dom \psi$. Fix $H := p L_p$.
		\FOR{$k = 0, 1, 2, \dots$}
		\STATE Pick up $\delta_{k + 1} \geq 0$
		\STATE Compute inexact monotone tensor step $x_{k + 1} := M_{H, \delta_{k + 1}}(x_k)$
		\ENDFOR
	\end{algorithmic}
\end{algorithm}

When our method is strictly monotone, we guarantee that
$F(x_{k + 1}) < F(x_k)$ for all $k \geq 0$, and we propose
to use the following \textit{adaptive strategy} of
defining the inner accuracies.

\BT \label{th:GlobalConvAlg2}
Let $D < +\infty$. Let sequence of inner accuracies $\{
\delta_k \}_{k \geq 1}$ be chosen in accordance to the
rule
\beq \label{DeltaDynamic}
\boxed{
	\ba{rcl}
	\delta_k & := &
	c \cdot  \bigl( F(x_{k - 2}) - F(x_{k - 1}) \bigr),
	\quad k \geq 2
	\ea
}
\eeq
for some fixed $0 \leq c < \frac{1}{(p + 2)3^{p + 1} - 1}$ and $\delta_1 \geq 0$.
Then for the sequence $\{ x_k \}_{k \geq 1}$
produced by Algorithm~\ref{alg:Monotone2}, we have
\beq \label{eq:GlobalConvAlg2}
\ba{rcl}
F(x_k) - F^{*} & \leq &
\frac{\gamma L_p D^{p + 1}}{p!\, k^p}
+ \frac{\beta}{k^{p + 2}},
\ea
\eeq
where $\gamma$ and $\beta$ are the constants:
$$
\ba{rcl}
\gamma & \Def &
\frac{(p + 2)^{p + 1}}{1 - c( (p + 2)3^{p + 1} - 1 )},  \;\;
\beta \; \Def \;
\frac{\delta_1 + c 2^{p + 2}(F(x_0) - F^{*})}{1 - c( (p + 2)^2 / (p + 1) - 1)}.
\ea
$$
\ET

The rule~\eqref{DeltaDynamic} is surprisingly simple and
natural: while the method is approaching the optimum, it
becomes more and more difficult to optimize the function.
Consequently, the progress in the function value at every
step is decreasing. Therefore, we need to solve the
auxiliary problem more accurately, and this is exactly
what we are doing in accordance to this rule.

It is also notable, that the rule~\eqref{DeltaDynamic}
is \textit{universal}, in a sense that it remains
the same (up to a constant factor)
for the methods of \textit{any order}, starting from $p = 1$.

This strategy also works for the nondegenerate case. Let
us assume that our objective is \textit{uniformly convex}
of degree $p + 1$ with constant $\sigma_{p + 1}$. Thus,
for all $x, y \in \dom \psi$ and $F'(x) \in \partial F(x)$
it holds
\beq \label{eq:UnifConv}
\ba{cl}
& F(y) - F(x) + \la F'(x), y - x \ra \\
\\
& \quad \geq \quad \frac{\sigma_{p + 1}}{p + 1}\|y - x\|^{p + 1}.
\ea
\eeq
For $p = 1$ this definition corresponds to the standard class
of \textit{strongly convex} functions.
One of the main sources of uniform convexity is a
regularization by power of Euclidean norm (we
use this construction in Section~\ref{sec:Acceleration},
where we accelerate our methods):

\BE \label{ExampleUnifConv}
Let $\psi(x) = \frac{\mu}{p + 1}\|x - x_0\|^{p + 1}$,
$\mu \geq 0$.
Then~\eqref{eq:UnifConv}
holds with $\sigma_{p + 1} = \mu 2^{1 - p}$
(see Lemma 5 in~\cite{doikov2019minimizing}).
\EE

\BE
Let $\psi(x) = \frac{\mu}{2}\|x - x_0\|^2$, $\mu \geq 0$.
Consider the ball of radius $D$ around the optimum:
$\mathcal{B} = \{ x : \|x - x^{*}\| \leq D  \}.$
Then~\eqref{eq:UnifConv} holds
for all $x, y \in \mathcal{B}$ with
$\sigma_{p + 1} = \frac{(p + 1)\mu}{2^p D^{p - 1}}$
(see Lemma~\ref{LemmaStrongUniformConvex}
in the supplementary material).
\EE

Denote by $\omega_p$ the \textit{condition number} of degree $p$:
\beq \label{eq-CondNumber}
\ba{rcl}
\omega_p & \Def & \max\{
\frac{(p + 1)^2  L_p}{p! \, \sigma_{p + 1}}, 1 \}.
\ea
\eeq

The next theorem shows, that $\omega_p$
serves as the main factor in the complexity of solving the
uniformly convex problems by inexact Tensor Methods.

\BT \label{th:GlobalStrongAlg2}
Let $\sigma_{p + 1} > 0$. Let sequence of inner accuracies
$\{ \delta_k \}_{k \geq 1}$ be chosen in accordance to the
rule
\beq \label{DeltaDynamic2}
\boxed{
	\ba{rcl}
	\delta_k & := &
	c \cdot  \bigl( F(x_{k - 2}) - F(x_{k - 1}) \bigr),
	\quad k \geq 2
	\ea
}
\eeq
for some fixed $0 \leq c < \frac{p}{p + 1}\omega_p^{-1 /
	p}$ and $\delta_1 \geq 0$. Then for the sequence $\{ x_k
\}_{k \geq 1}$ produced by Algorithm~\ref{alg:Monotone2},
we have the following \underline{linear} rate of
convergence:
\beq \label{eq:GlobalStrongAlg2}
\ba{cl}
& F(x_{k + 1}) - F^{*} \\
\\
& \;\; \leq \;\;
\Bigl(1 - \frac{p}{p + 1} \omega_p^{-1/p} + c \Bigr) (F(x_{k - 1}) - F^{*}).
\ea
\eeq
\ET
Let us pick $c := \frac{p}{2(p + 1)}\omega_p^{-1/p}$.
Then, according to estimate~\eqref{eq:GlobalStrongAlg2},
in order solve the problem up to $\varepsilon$ accuracy:
$F(x_K) - F^{*} \leq \varepsilon$, we need to perform
\beq \label{LinearComplexity}
\ba{rcl}
K & = & O\left(
\omega_p^{1 / p} \log \frac{F(x_0) - F^{*}}{\varepsilon}
\right)
\ea
\eeq
iterations of the algorithm.

Finally, we study the local behavior of the method for
strongly convex objective.
\BT \label{th:LocalAlg2}
Let $\sigma_2 > 0$. Let sequence of inner accuracies $\{
\delta_k \}_{k \geq 1}$ be chosen in accordance to the
rule
\beq \label{DeltaDynamic3}
\boxed{
	\ba{rcl}
	\delta_k & := &
	c \cdot \bigl( F(x_{k - 2}) - F(x_{k - 1}) \bigr)^{p + 1 \over 2},
	\, k \geq 2
	\ea
}
\eeq
with some fixed $c \geq 0$ and $\delta_1 \geq 0$. Then for
$p \geq 2$ the sequence $\{ x_k \}_{k \geq 1}$ produced by
Algorithm~\ref{alg:Monotone2} has the \underline{local
	superlinear rate of convergence}:
\beq \label{eq:LocalAlg2}
\ba{cl}
& F(x_{k + 1}) - F^{*} \\
\\
& \;\; \leq \;\;
\Bigl( \frac{L_p}{p!}\bigl( \frac{2}{\sigma_2} \bigr)^{p + 1 \over 2} + c \Bigr)
(F(x_{k - 1}) - F^{*})^{p + 1 \over 2}.
\ea
\eeq
\ET

Let us assume for simplicity, that the constant $c$ is chosen to be small enough:
$c \leq \frac{L_p}{p!}\bigl( \frac{2}{\sigma_2} \bigr)^{(p + 1) / 2}$.
Then, we are able to describe the region of superlinear convergence as
$$
\ba{cl}
\mathcal{Q} = \left\{
x \in \dom \psi : F(x) - F^{*} \leq
\Bigl(\frac{\sigma_2^{p + 1}}{2^{p + 3}} \bigl( \frac{p!}{L_p} \bigr)^2 \Bigr)^{1 \over (p - 1)}
\right\}.
\ea
$$
After reaching it, the method becomes very fast: we need
to perform no more than $O(\log \log
\frac{1}{\varepsilon})$ additional iterations to solve the
problem.

Note, that estimate~\eqref{eq:LocalAlg2} of the local
convergence is slightly weaker, than the corresponding one
for \textit{exact} Tensor Methods~\cite{doikov2019local}.
For example, for $p = 2$ (Cubic regularization of Newton
Method) we obtain the convergence of order~${3 \over 2}$,
not the quadratic, which affects only a constant factor in
the complexity estimate. The region~$\mathcal{Q}$ of the
superlinear convergence is remaining the same.

\subsection{Inexact Methods with Averaging}
\label{subsec:Averaging}

Methods from the previous section were developed by
forcing the monotonicity of the sequence of function
values $\{ F(x_k) \}_{k \geq 0}$ into the scheme. As a
byproduct, we get the radius of the initial level set $D$
(see definition~\eqref{eq-DRad}) in the right-hand side of
our complexity estimates \eqref{eq:GlobalConvAlg1}
and~\eqref{eq:GlobalConvAlg2}. Note, that $D$ may be
significantly bigger than the distance $\|x_0 - x^{*}\|$
from the initial point to the solution.
\BE \label{ex:BadFunction}
Consider the following function, for $x \in \R^n$:
$$
\ba{rcl}
f(x) & = & |x^{(1)}|^{p + 1}
+ \sum\limits_{i = 2}^n |x^{(i)} - 2x^{(i - 1)}|^{p + 1},
\ea
$$
where $x^{(i)}$ indicates $i$th coordinate of $x$. Clearly,
the minimum of $f$ is at the origin: $x^{*} = (0, \dots, 0)^T$.
Let us take two points: $x_0 = (1, \dots, 1)^T$ and
$x_1$, such that $x_1^{(i)} = 2^i - 1$.
It holds, $f(x_0) = f(x_1) = n$, so they belong to the same level set.
However, we have (for the standard Euclidean norm): $\|x_0 - x^{*}\| = \sqrt{n}$,
while
$D \geq \|x_1 - x^{*}\| \geq 2^{n - 1}$.
\EE

Here we present an alternative approach, Tensor Methods
with \textit{Averaging}. In this scheme, we perform a step
not from the previous point $x_k$, but from a point $y_k$,
which is a convex combination of the previous point and
the starting point:
$$
\ba{rcl}
y_{k}  & = & \lambda_{k} x_k
+ (1 - \lambda_{k}) x_0,
\ea
$$
where $\lambda_{k} \equiv \bigl(\frac{k}{k + 1}\bigr)^{p + 1}$.
The whole optimization scheme remains very simple.

\begin{algorithm}[!h]
	\caption{Inexact Tensor Method with Averaging}
	\label{alg:Averaging}
	\begin{algorithmic}
		\STATE {\bfseries Initialization:} Choose $x_0 \in \dom \psi$. Fix $H := pL_p$.
		\FOR{$k = 0, 1, 2, \dots$}
		\STATE Set $\lambda_{k} := \bigl( \frac{k}{k + 1} \bigr)^{p + 1}$,
		$\; y_{k} := \lambda_{k} x_k + (1 - \lambda_{k}) x_0$
		\STATE Pick up $\delta_{k + 1} \geq 0$
		\STATE Compute inexact tensor step $x_{k + 1} := T_{H, \delta_{k + 1}}(y_{k})$
		\ENDFOR
	\end{algorithmic}
\end{algorithm}

For this method, we are able to prove a similar
convergence result as that of
Algorithm~\ref{alg:Monotone1}. However, now we have the
explicit distance $\|x_0 - x^{*}\|$ in the right hand side
of our bound for the convergence rate (compare with
Theorem~\ref{th:GlobalConvAlg1}).

\BT \label{th:GlobalConvAlg3}
Let sequence of inner accuracies $\{ \delta_k \}_{k \geq 1}$
be chosen according the rule
\beq \label{DeltaKChoice3}
\boxed{
	\ba{rcl}
	\delta_k & := & \frac{c}{k^{p + 1}}
	\ea
}
\eeq
for some $c \geq 0$. Then for the sequence $\{ x_k \}_{k \geq 1}$
produced by Algorithm~\ref{alg:Averaging}, we have
\beq \label{eq:GlobalConvAlg3}
\ba{rcl}
F(x_k) - F^{*} & \leq &
\frac{(p + 1)^{p + 1} L_p \|x_0 - x^{*}\|^{p + 1}}{p! \, k^{p}}
+ \frac{c}{k^p}.
\ea
\eeq
\ET

Thus, Algorithm~\ref{alg:Averaging} seems
to be the first \textit{Primal} Tensor method
(aggregating only the points from the primal space $\E$),
which admits the explicit initial distance
in the global convergence estimate~\eqref{eq:GlobalConvAlg3}.
Table~\ref{table:Methods} contains a short overview
of the inexact Tensor methods from this section.

\begin{table}[!h]
	\caption{Summary on the methods.}
	\label{table:Methods}
	\vskip 0.15in
	\begin{center}
		\scriptsize
		\setlength{\tabcolsep}{1pt}
		\begin{tabular}{cccc}
			\toprule
			\textbf{Algorithm}
			&
			\begin{tabular}{c}
				\textbf{The rule for} $\boldsymbol{\delta_k}$
			\end{tabular}
			&
			\begin{tabular}{c}
				\textbf{Global} \\
				\textbf{rate}
			\end{tabular}
			&
			\begin{tabular}{c}
				\textbf{Local} \\
				\textbf{superl.}
			\end{tabular}
			\\
			\midrule
			\begin{tabular}{c}
				Tensor Method \\
				\cite{nesterov2019implementable}
			\end{tabular}
			&
			0
			&
			$O\Bigl(\frac{L_p D^{p + 1}}{k^p}\Bigr)$
			&
			Yes
			\\
			\midrule
			\begin{tabular}{c}
				Monotone Inexact \\
				Tensor Method, I \\
				(Algorithm~\ref{alg:Monotone1})
			\end{tabular}
			&
			$1 / k^{p + 1}$
			&
			$O\Bigl(\frac{L_p D^{p + 1}}{k^p}\Bigr)$
			&
			No \\
			\midrule
			\begin{tabular}{c}
				Monotone Inexact \\
				Tensor Method, II \\
				(Algorithm~\ref{alg:Monotone2})
			\end{tabular}
			&
			$\, (F(x_{k - 1}) - F(x_{k}))^{\alpha}$
			&
			\begin{tabular}{c}
				$O\Bigl(\frac{L_p D^{p + 1}}{k^p}\Bigr),$ \\[3pt]
				$\; \alpha = 1$
			\end{tabular}
			&
			\begin{tabular}{c}
				Yes, \\[3pt]
				$\alpha = \frac{p + 1}{2}$
			\end{tabular}
			\\
			\midrule
			\begin{tabular}{c}
				Inexact \\
				Tensor Method \\
				with Averaging \\
				(Algorithm~\ref{alg:Averaging})
			\end{tabular}
			&
			$1 / k^{p + 1}$
			&
			$O\Bigl(\frac{L_p \|x_0 - x^{*}\|^{p + 1}}{k^p}\Bigr)$
			&
			No \\
			\bottomrule
		\end{tabular}
	\end{center}
	\vskip -0.1in
\end{table}

\section{Acceleration}
\label{sec:Acceleration}

After the Fast Gradient Method
had been discovered in~\cite{nesterov1983method},
there were made huge efforts to develop accelerated
second-order~\cite{nesterov2008accelerating,monteiro2013accelerated,grapiglia2019accelerated}
and high-order~\cite{baes2009estimate,nesterov2019implementable,gasnikov2019near,grapiglia2019tensor,song2019towards}
optimization algorithms.
Most of these schemes are based on the notion of Estimating sequences~(see \cite{nesterov2018lectures}).

An alternative approach of using the Proximal Point
iterations with the acceleration was studied first
in~\cite{guler1992new}. It became very popular recently,
in the context of machine learning
applications~\cite{lin2015universal,lin2018catalyst,kulunchakov2019generic,ivanova2019adaptive}.
In this section, we use the technique of Contracting
Proximal iterations~\cite{doikov2019contracting}, to
accelerate our inexact tensor methods.

In the accelerated scheme, two sequences of points are used:
the main sequence $\{ x_k \}_{k \geq 0}$,
for which we are able to
guarantee the convergence in function residuals,
and auxiliary sequence $\{ v_k \}_{k \geq 0}$
of prox-centers,
starting from the same initial point: $v_0 = x_0$.
Also, we use the sequence $\{ A_k \}_{k \geq 0}$
of scaling coefficients. Denote
$
a_k \Def A_k - A_{k - 1}, k \geq 1.
$

Then, at every iteration, we apply
Monotone Inexact Tensor Method, II (Algorithm~\ref{alg:Monotone2})
to minimize the following \textit{contracted}
objective with regularization:
$$
\ba{rcl}
h_{k + 1}(x) & \Def &
A_{k + 1}f\bigl(\frac{a_{k + 1}x + A_k x_k}{A_{k + 1}}\bigr)
+ a_{k + 1}\psi(x) \\[7pt]
& \; & \; + \; \beta_d(v_k; x).
\ea
$$
Here
$\beta_d(v_k; x) \Def d(x) - d(v_k) - \la \nabla d(v_k), x - v_k \ra$
is \textit{Bregman divergence} centered at $v_k$,
for the following choice of prox-function:
$$
\ba{rcl}
d(x) & := & \frac{1}{p + 1}\|x - x_0\|^{p + 1},
\ea
$$
which is uniformly convex of degree $p + 1$
(Example~\ref{ExampleUnifConv}). Therefore, Tensor Method
achieves fast linear rate of convergence
(Theorem~\ref{th:GlobalStrongAlg2}). By an appropriate
choice of scaling coefficients $\{ A_k \}_{k \geq 1}$, we
are able to make the condition number of the subproblem
being an absolute constant. This means that only
$\tilde{O}(1)$ steps of Algorithm~\ref{alg:Monotone2} are
needed to find an approximate minimizer of $h_{k +
	1}(\cdot)$:
\beq \label{StopCondition}
\ba{rcl}
h_{k + 1}(v_{k + 1}) - h_{k + 1}^{*} & \leq & \zeta_{k + 1}.
\ea
\eeq
Note that inexact condition~\eqref{StopCondition} was
considered first in~\cite{lin2015universal}, in a general
algorithmic framework for accelerating first-order
methods. It differs from the corresponding one
from~\cite{doikov2019contracting}, where a bound for the
(sub)gradient norm was used.

\begin{algorithm}[!h]
	\caption{Accelerated Scheme}
	\label{alg:Accelerated}
	\begin{algorithmic}
		\STATE {\bfseries Initialization:} Choose $x_0 \in \dom \psi$.
		Set $v_0 := x_0$, $A_0 := 0$.
		\FOR{$k = 0, 1, 2, \dots$}
		\STATE Set $A_{k + 1} := \frac{(k + 1)^{p + 1}}{L_p}$
		\STATE Pick up $\zeta_{k + 1} \geq 0$
		\STATE Find $v_{k + 1}$ such that~\eqref{StopCondition} holds
		\STATE Set $x_{k + 1} := \frac{a_{k + 1} v_{k + 1} + A_k x_k}{A_{k + 1}}$
		\ENDFOR
	\end{algorithmic}
\end{algorithm}

Therefore, for accelerating inexact Tensor Methods, we
propose a multi-level approach. On the upper level, we run
Algorithm~\ref{alg:Accelerated}. At each iteration of this
method, we call Algorithm~\ref{alg:Monotone2} (with
dynamic rule~\eqref{DeltaDynamic2} for inner accuracies)
to find $v_{k + 1}$. For this optimization scheme, we
obtain the following global convergence guarantee.
\BT\label{th:GlobalAlg4}
Let sequence $\{ \zeta_k \}_{k \geq 1}$ be chosen
according to the rule
\beq \label{ZetaKChoice}
\boxed{
	\ba{rcl}
	\zeta_k & := & \frac{c}{k^{p + 2}}
	\ea
}
\eeq
with some $c \geq 0$. Then for the iterations $\{x_k\}_{k
	\geq 1}$ produced by Algorithm~\ref{alg:Accelerated}, it
holds:
\beq \label{eq-GlobalAccelerated}
\ba{rcl}
F(x_k) - F^{*} & \leq & O\Bigl(\frac{L_p(\|x_0 - x^{*}\|^{p + 1} + c)}{k^{p + 1}}\Bigr).
\ea
\eeq
For every $k \geq 0$, in order to find $v_{k + 1}$ by
Algorithm~\ref{alg:Monotone2} (for minimizing
$h_{k+1}(\cdot)$, starting from $v_k$), it is enough to
perform no more than
\beq \label{eq-NStepsBound}
\ba{rcl}
O\Bigl( \log \frac{(k + 1)(\|x_0 - x^{*}\|^{p + 1} + c)}{c} \Bigr).
\ea
\eeq
inexact monotone tensor steps.
\ET

Therefore, the total number of the inexact tensor steps
for finding $\varepsilon$-solution of~\eqref{MainProblem}
is bounded by $\widetilde{O}\bigl( 1 / \varepsilon^{1
	\over {p + 1}} \bigr)$. One theoretical question remains
open: is it possible to construct in the framework of
inexact tensor steps, the \textit{optimal methods} with
the complexity estimate $O\bigl(1 / \varepsilon^{2 \over
	{3p + 1}} \bigr)$ having no hidden logarithmic factors.
This would match the existing lower
bound~\cite{arjevani2019oracle,nesterov2019implementable}.

\section{Experiments}
\label{sec:Experiments}

Let us demonstrate computational results with 
empirical study of different accuracy policies.
We consider inexact methods of order $p = 2$ 
(Cubic regularization of Newton method),
and to solve the corresponding subproblem 
we call the Fast Gradient Method with restarts from~\cite{nesterov2019inexact}.
To estimate the residual in function value 
of the subproblem, we use a simple stopping criterion, given by uniform convexity of the model $g(y) = \Omega_H(x; y)$:
\beq \label{UnifConvBound}
\ba{rcl}
g(y) - \min\limits_{y} g(y) & \leq & \frac{4}{3}\left( \frac{1}{H} \right)^{1/2} \| \nabla g(y) \|_{*}^{3 / 2}.
\ea
\eeq
An alternative approach would be to bound the functional residual by the duality gap\footnote{Note that
	the left hand side in~\eqref{UnifConvBound}
	can be bounded from below by the distance from $y$ to the optimum of the model, using uniform convexity. Therefore,
we have a computable bound for the distance to the solution of the subproblem.}.

We compare the \textit{adaptive} rule for inner accuracies~\eqref{DeltaDynamic}
with dynamic strategies in the form $\delta_k = 1/k^{\alpha}$, for
different~$\alpha$
(left graphs), and with the constant choices (right).

\subsection{Logistic Regression}

First, let us consider the problem of training $\ell_2$-regularized logistic regression
model for classification task with two classes, on several real datasets\footnote{\url{https://www.csie.ntu.edu.tw/~cjlin/libsvmtools/datasets/}}: \textit{mashrooms} $(m=8124, n=112)$, \textit{w8a} $(m=49749, n=300)$, and \textit{a8a} $(m=22696,n=123)$\footnote{$m$ is the number of training examples and $n$ is the dimension of the problem (the number of features).}.

\begin{figure}[h!]
	\begin{minipage}{0.210\textwidth}
		\centering
		\includegraphics[width=\textwidth ]{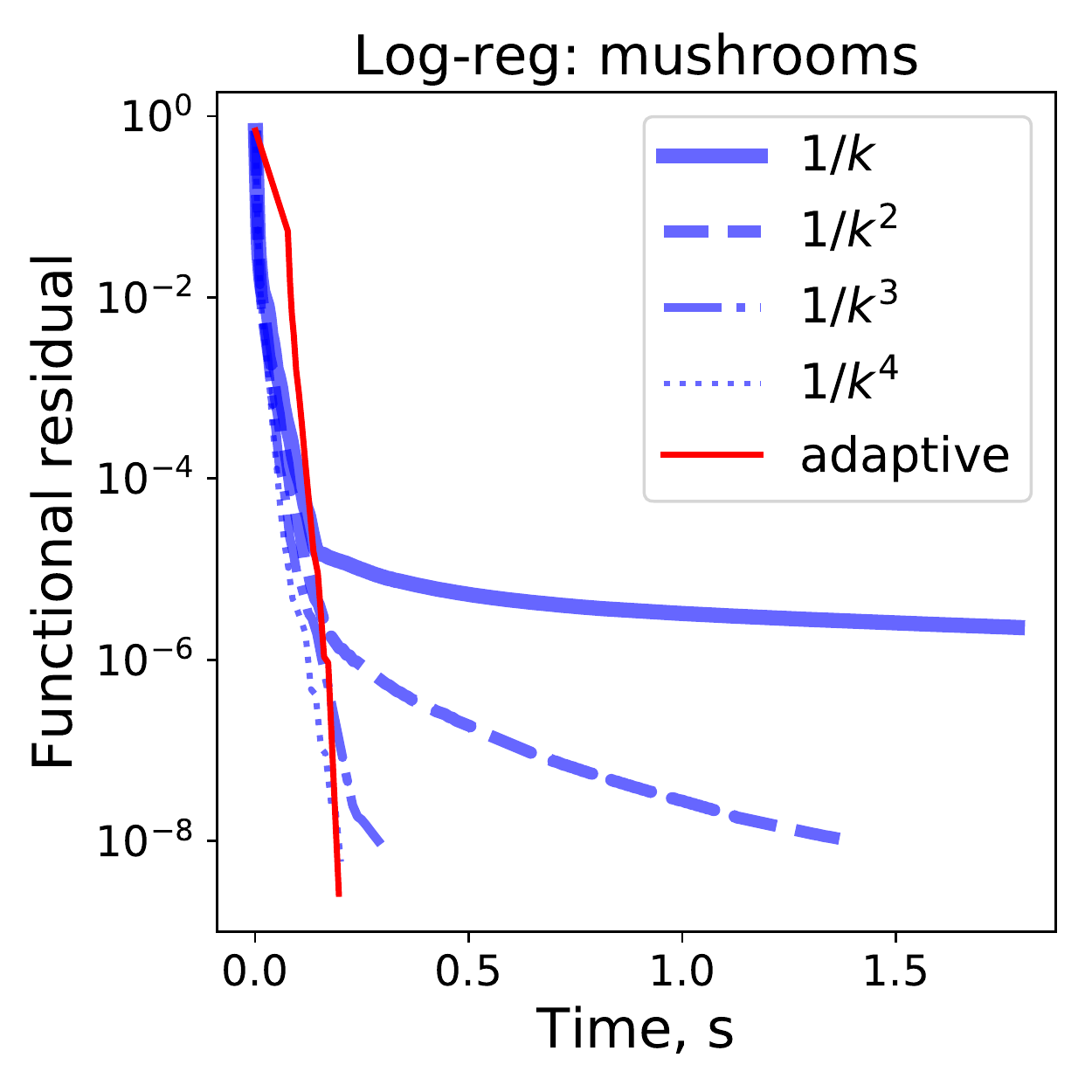}
	\end{minipage}
	\begin{minipage}{0.210\textwidth}
		\centering
		\includegraphics[width=\textwidth ]{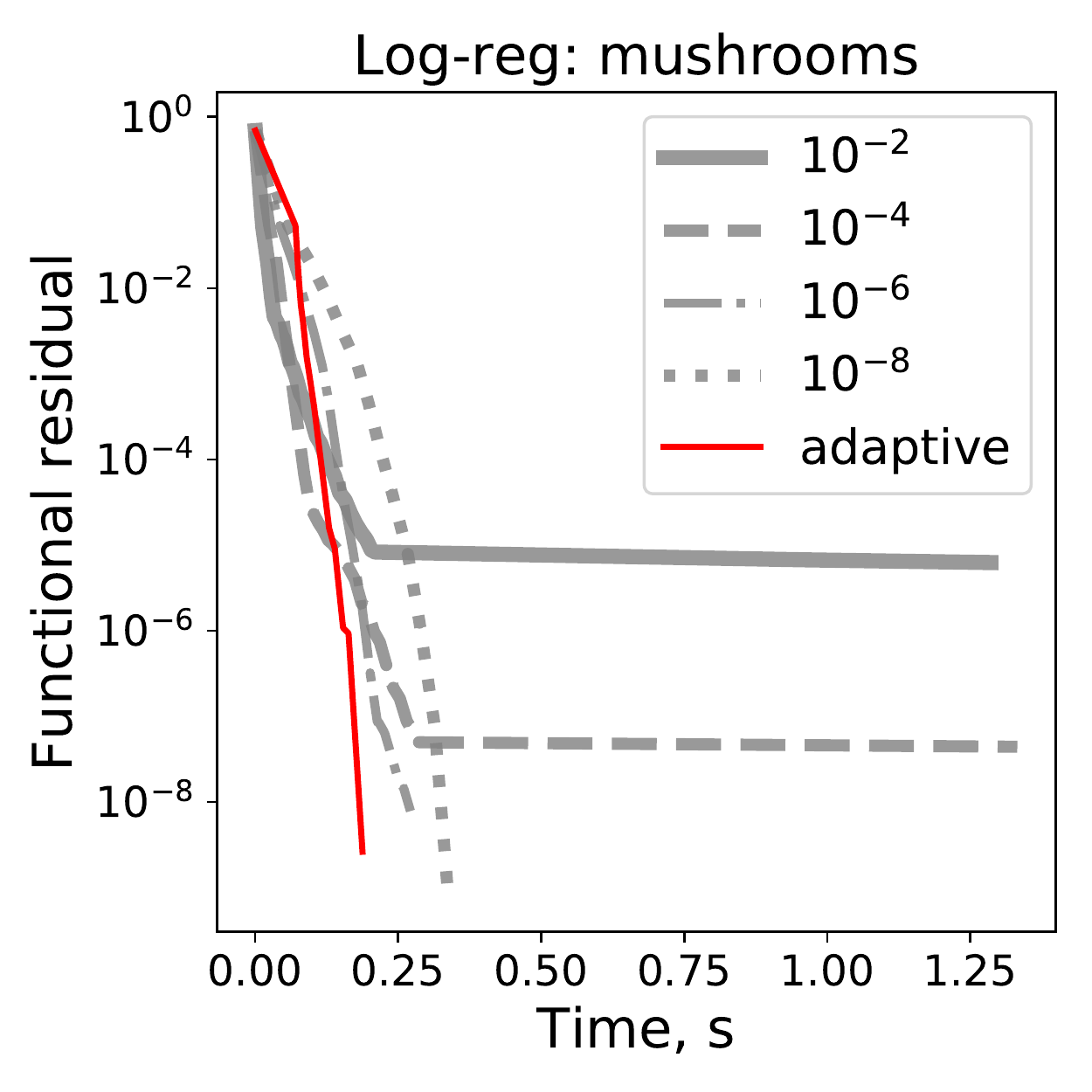}
	\end{minipage}
	
	\begin{minipage}{0.210\textwidth}
		\centering
		\includegraphics[width=\textwidth ]{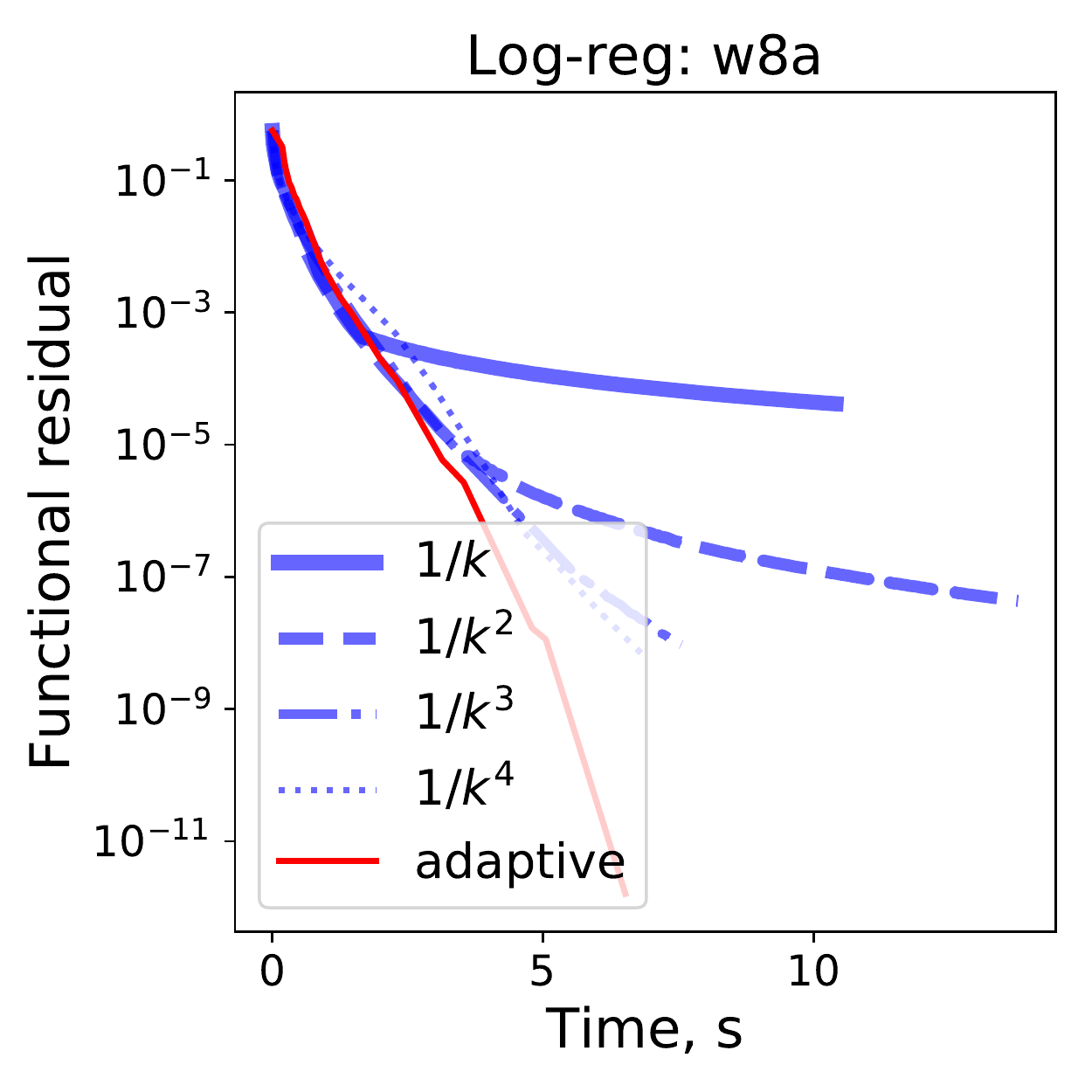}
	\end{minipage}
	\begin{minipage}{0.210\textwidth}
		\centering
		\includegraphics[width=\textwidth ]{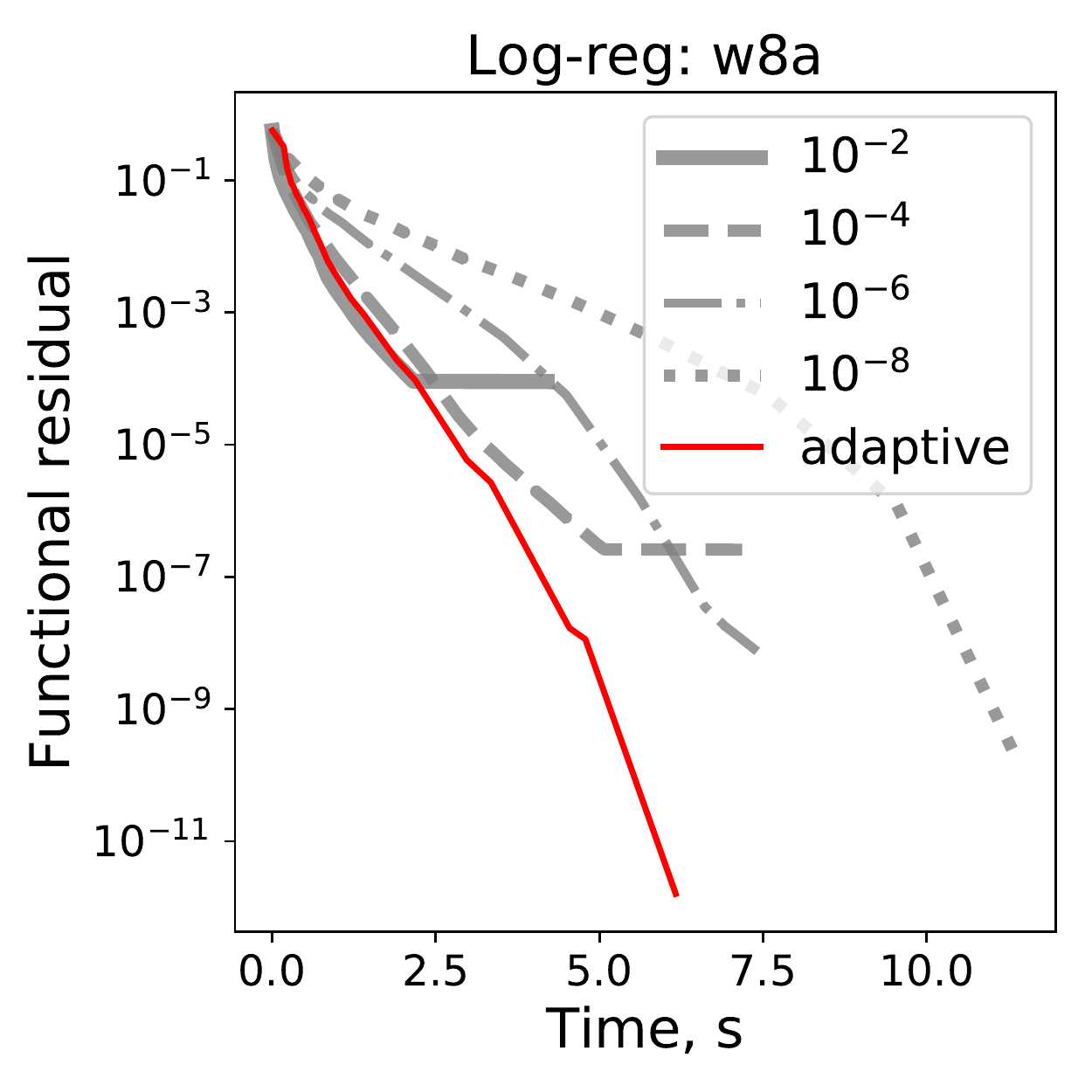}
	\end{minipage}

	\begin{minipage}{0.210\textwidth}
		\centering
		\includegraphics[width=\textwidth ]{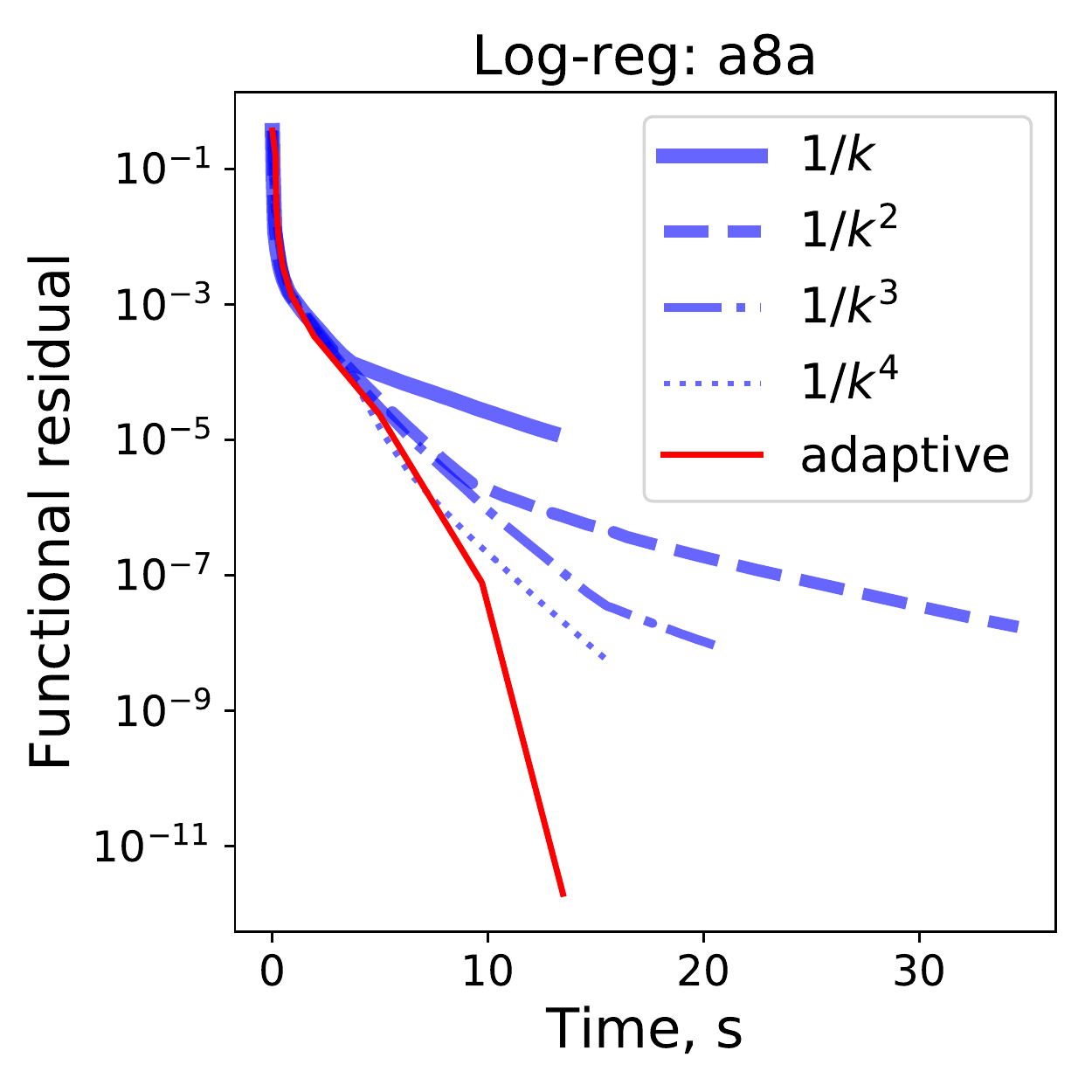}
	\end{minipage}
	\begin{minipage}{0.210\textwidth}
		\centering
		\includegraphics[width=\textwidth ]{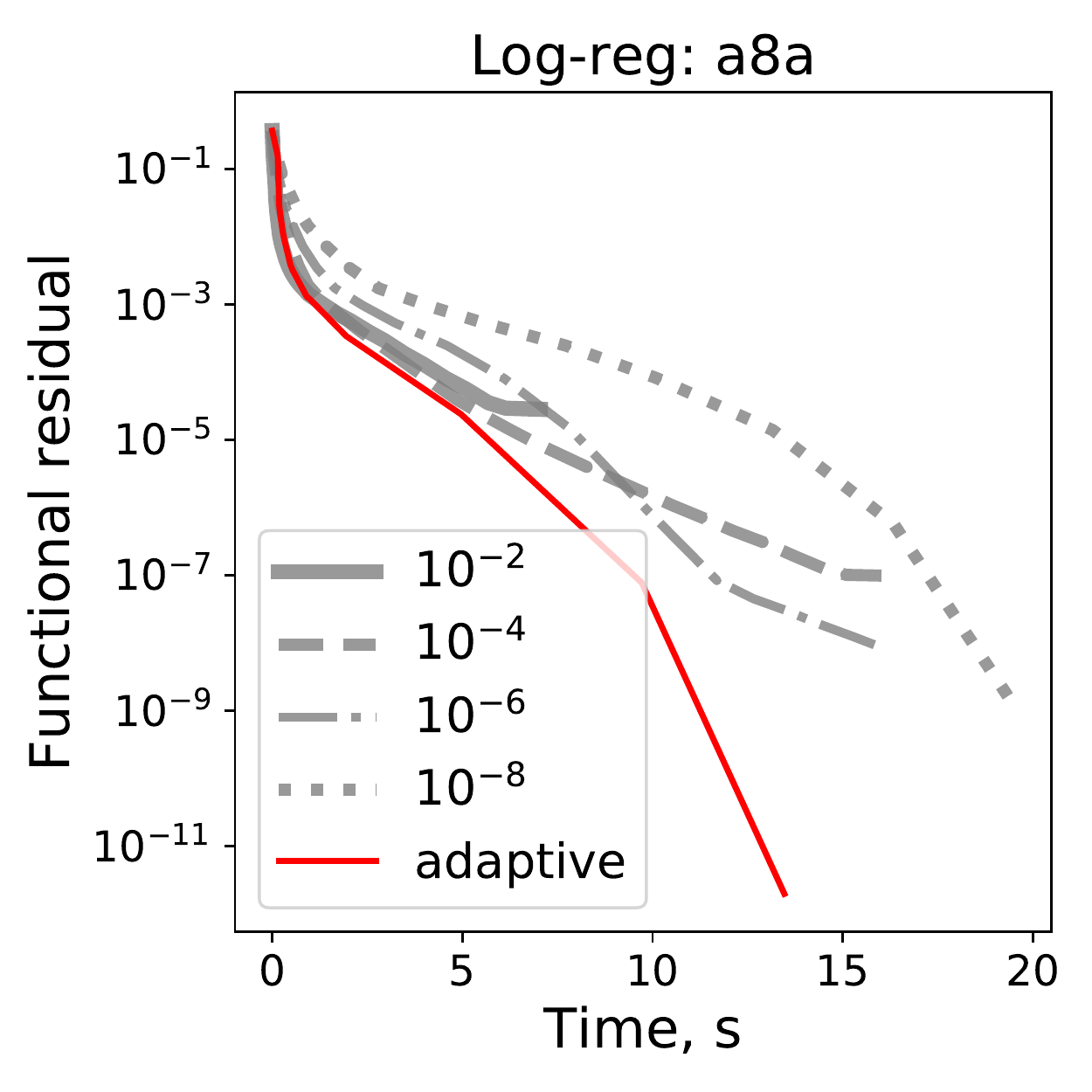}
	\end{minipage}
	\caption{Comparison of different accuracy policies for inexact Cubic Newton, training logistic regression.}\label{fig:CN_logistic_regression}
\end{figure}

We use the standard Euclidean norm for this problem, and simple line search at every iteration, to fit the regularization parameter $H$. The results are shown on Figure~\ref{fig:CN_logistic_regression}.

\subsection{Log-Sum-Exp}

In the next set of experiments, we consider unconstrained minimization
of the following objective:
$$
\ba{rcl}
f_{\mu}(x) & = & \mu \ln \left( 
\sum\limits_{i = 1}^m \exp\Bigl( \frac{\la a_i, x \ra - b_i}{\mu} \Bigr) 
\right), \;\; x \in \R^n,
\ea
$$
where $\mu > 0$ is a \textit{smoothing} parameter. 
To generate the data, we sample coefficients $\{ \tilde{a}_i \}_{i = 1}^m$ and $b$ randomly from the uniform distribution on $[-1, 1]$.
Then, we shift the parameters in a way to have the solution $x^{*}$ in the origin. Namely, using $\{ \tilde{a}_i \}_{i = 1}^m$
we form a preliminary function $\tilde{f}_{\mu}(x)$, and set
$a_i := \tilde{a}_i - \nabla \tilde{f}_{\mu}(0)$. Thus we essentially obtain
$\nabla f_{\mu}(0) = 0$.

\begin{figure}[h!]
	\begin{minipage}{0.235\textwidth}
		\centering
		\includegraphics[width=\textwidth ]{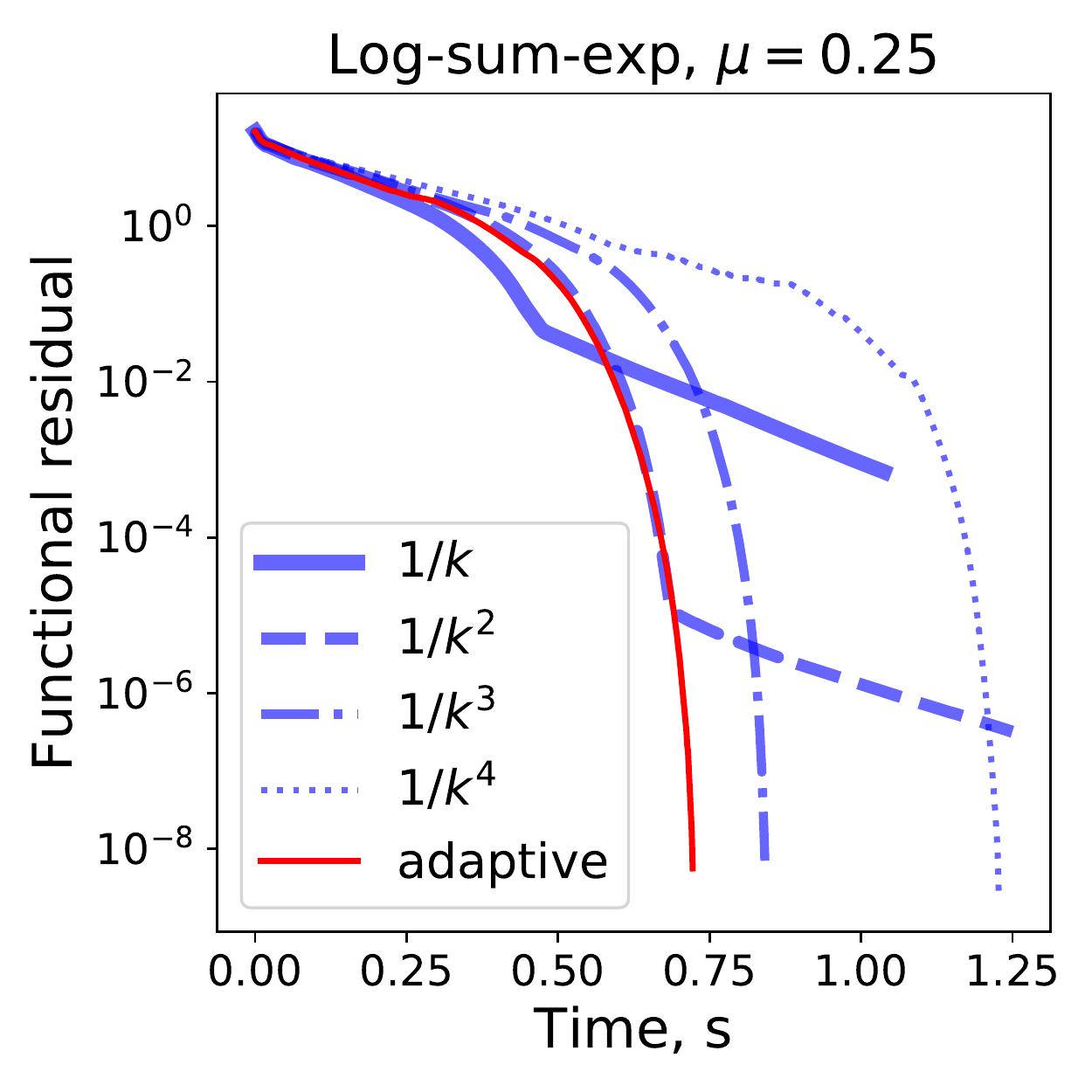}
	\end{minipage}
	\begin{minipage}{0.235\textwidth}
		\centering
		\includegraphics[width=\textwidth ]{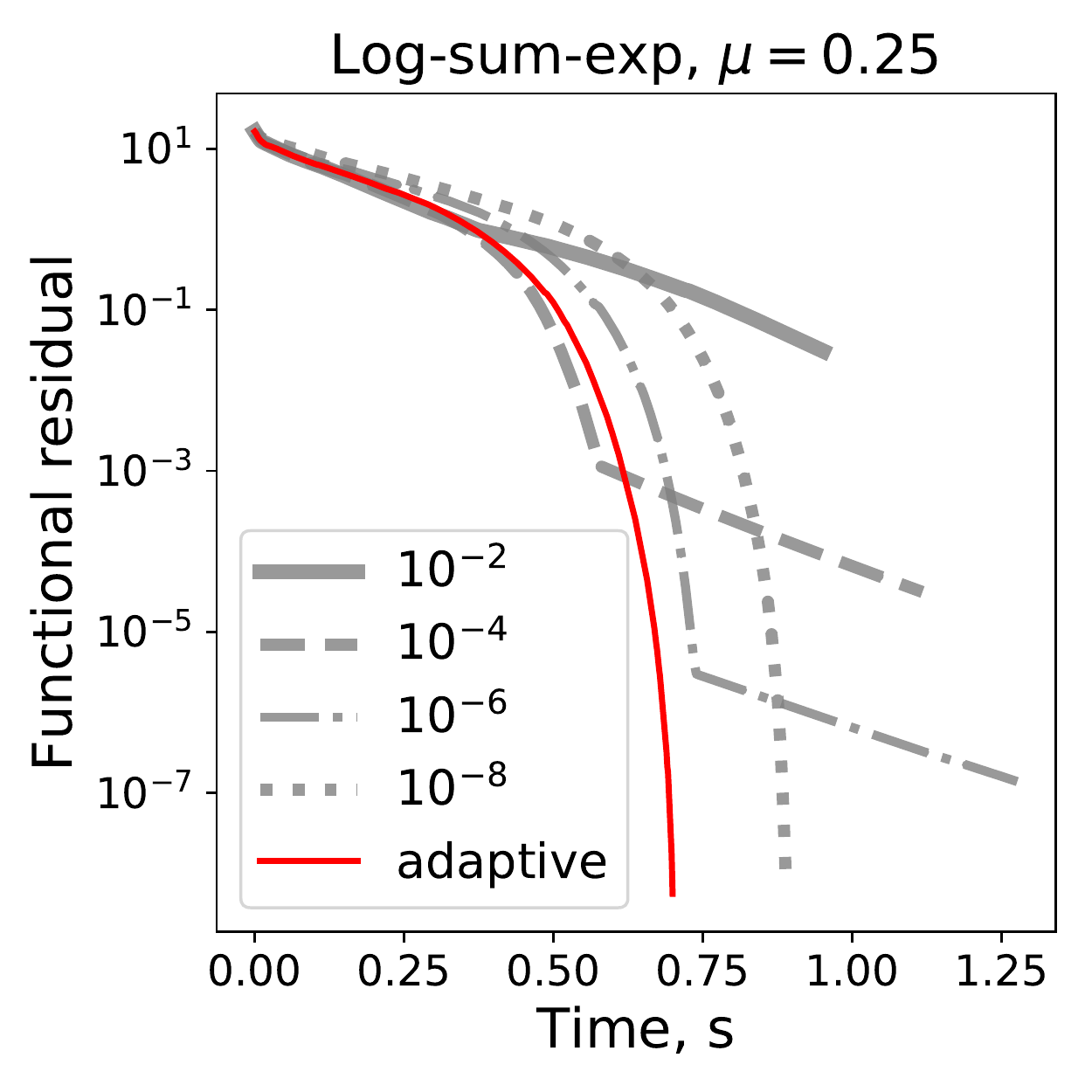}
	\end{minipage}

	\begin{minipage}{0.235\textwidth}
		\centering
		\includegraphics[width=\textwidth ]{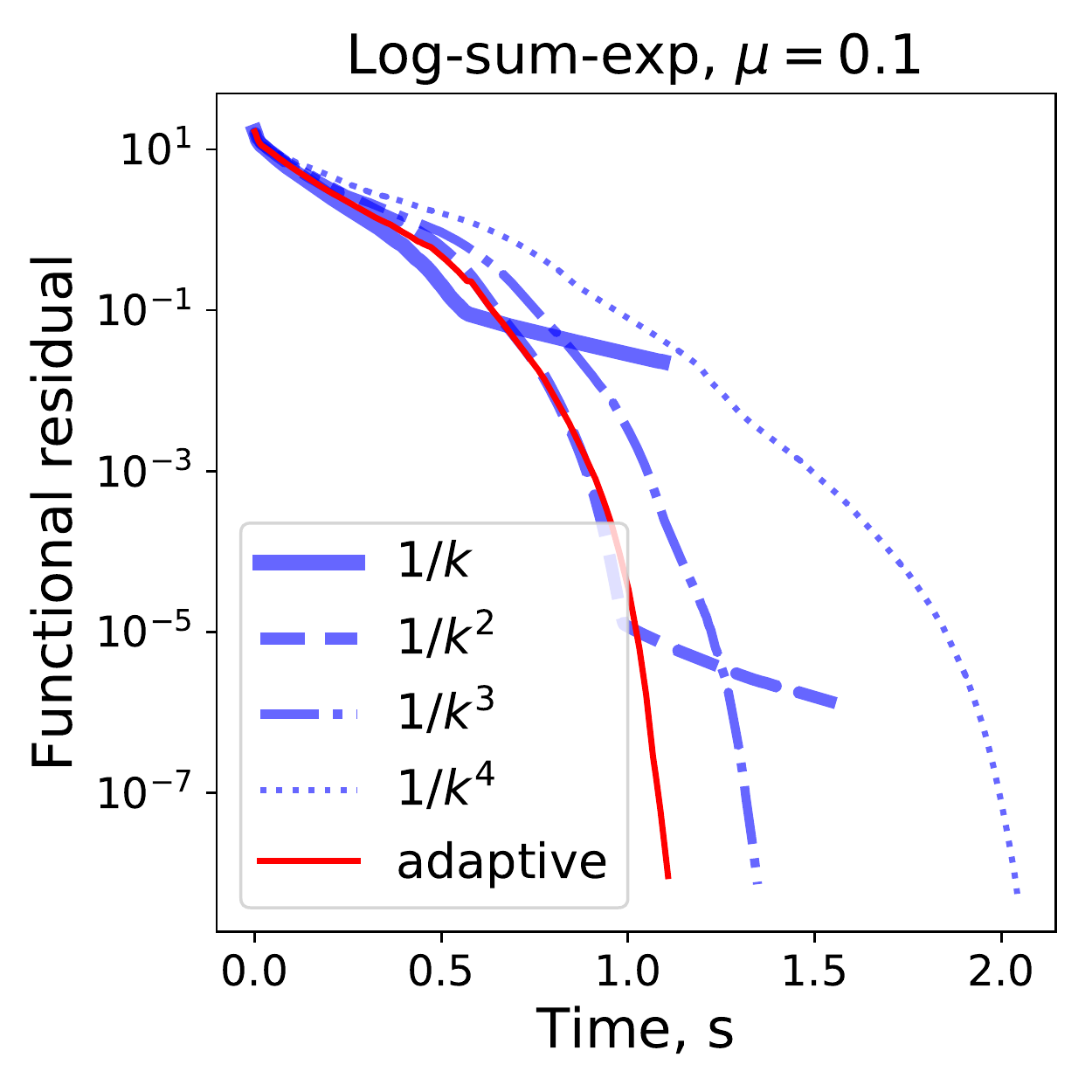}
	\end{minipage}
	\begin{minipage}{0.235\textwidth}
		\centering
		\includegraphics[width=\textwidth ]{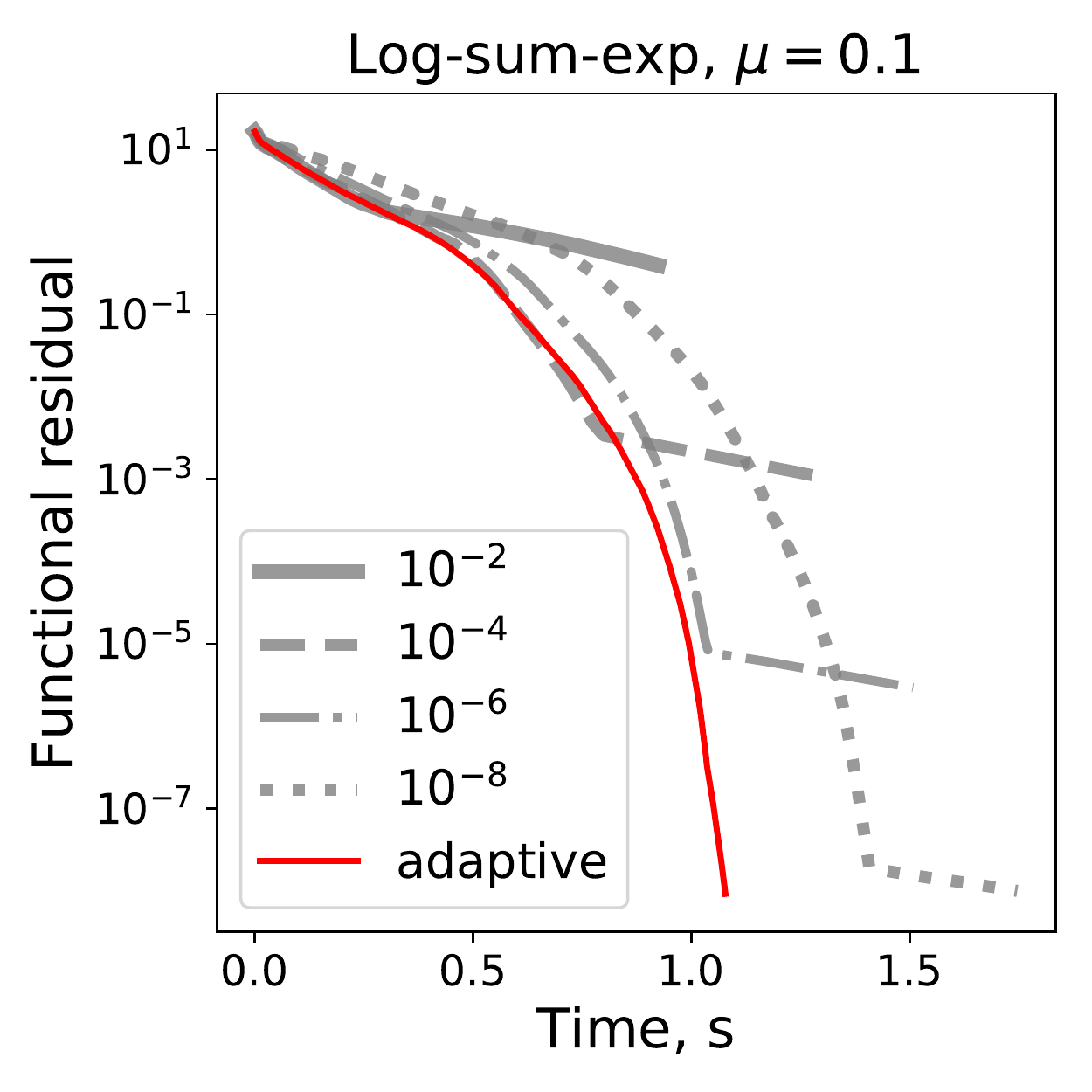}
	\end{minipage}

	\begin{minipage}{0.235\textwidth}
		\centering
		\includegraphics[width=\textwidth ]{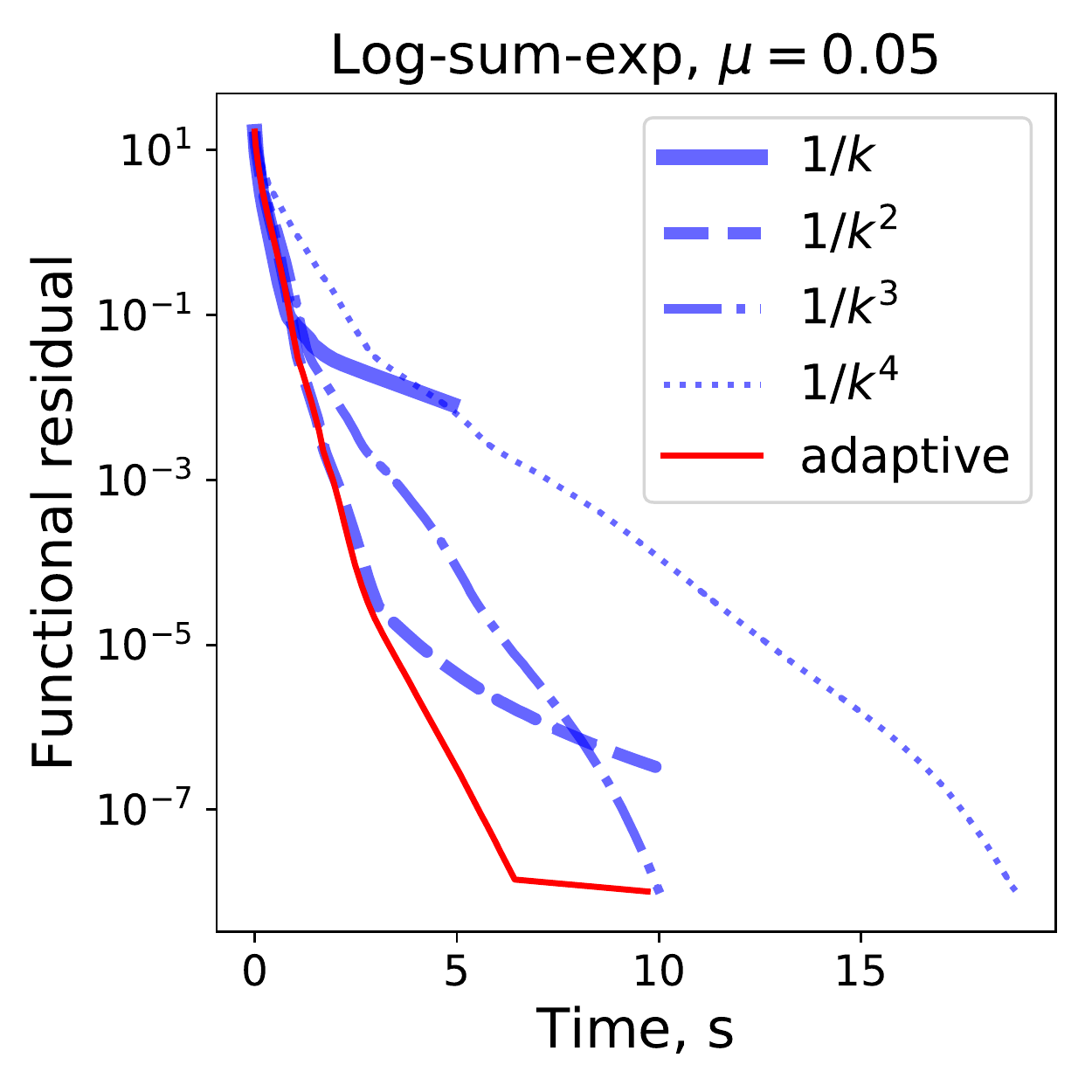}
	\end{minipage}
	\begin{minipage}{0.235\textwidth}
		\centering
		\includegraphics[width=\textwidth ]{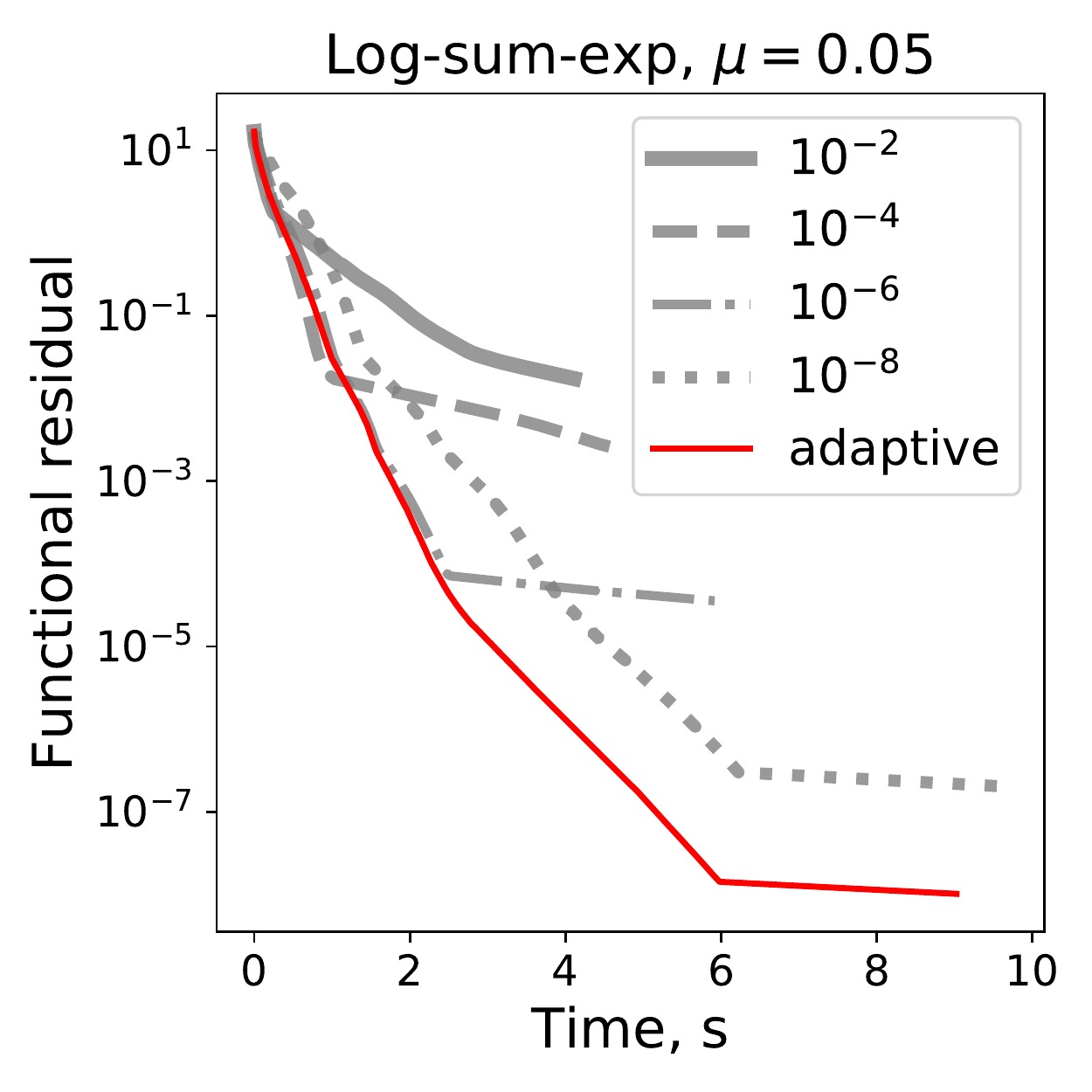}
	\end{minipage}
	\caption{Comparison of different accuracy policies for inexact Cubic Newton, minimizing Log-Sum-Exp function.} \label{fig:CN_log_sum_exp}
\end{figure}

We set $m = 6n$, and $n = 100$.
In the method, we use the following Euclidean norm 
for the primal space: $\|x\| = \la Bx, x \ra^{1/2}$, 
with the matrix $B = \sum_{i = 1}^m a_i a_i^T$,
and fix regularization parameter $H$ being equal $1$.
The results are shown on Figure~\ref{fig:CN_log_sum_exp}.

We see, that the adaptive rule
demonstrates reasonably good performance (in terms of the total computational time\footnote{Clock time was evaluated using the machine with Intel Core i5 CPU, 1.6GHz; 8 GB RAM. All methods were implemented in Python. The source code can be found at \url{https://github.com/doikov/dynamic-accuracies/}}) in all the scenarios.

\section*{Acknowledgements}

The research results of this paper were obtained in the framework of ERC Advanced Grant 788368.

\bibliography{bibliography}

\begin{thebibliography}{48}
\providecommand{\natexlab}[1]{#1}
\providecommand{\url}[1]{\texttt{#1}}
\expandafter\ifx\csname urlstyle\endcsname\relax
  \providecommand{\doi}[1]{doi: #1}\else
  \providecommand{\doi}{doi: \begingroup \urlstyle{rm}\Url}\fi

\bibitem[Agarwal et~al.(2017)Agarwal, Allen-Zhu, Bullins, Hazan, and
  Ma]{agarwal2017finding}
Agarwal, N., Allen-Zhu, Z., Bullins, B., Hazan, E., and Ma, T.
\newblock Finding approximate local minima faster than gradient descent.
\newblock In \emph{Proceedings of the 49th Annual ACM SIGACT Symposium on
  Theory of Computing}, pp.\  1195--1199. ACM, 2017.

\bibitem[Arjevani et~al.(2019)Arjevani, Shamir, and Shiff]{arjevani2019oracle}
Arjevani, Y., Shamir, O., and Shiff, R.
\newblock Oracle complexity of second-order methods for smooth convex
  optimization.
\newblock \emph{Mathematical Programming}, 178\penalty0 (1-2):\penalty0
  327--360, 2019.

\bibitem[Baes(2009)]{baes2009estimate}
Baes, M.
\newblock Estimate sequence methods: extensions and approximations.
\newblock \emph{Institute for Operations Research, ETH, Z{\"u}rich,
  Switzerland}, 2009.

\bibitem[Bauschke et~al.(2016)Bauschke, Bolte, and
  Teboulle]{bauschke2016descent}
Bauschke, H.~H., Bolte, J., and Teboulle, M.
\newblock A descent lemma beyond lipschitz gradient continuity: first-order
  methods revisited and applications.
\newblock \emph{Mathematics of Operations Research}, 42\penalty0 (2):\penalty0
  330--348, 2016.

\bibitem[Beck(2017)]{beck2017first}
Beck, A.
\newblock \emph{First-order methods in optimization}, volume~25.
\newblock SIAM, 2017.

\bibitem[Birgin et~al.(2017)Birgin, Gardenghi, Mart{\'\i}nez, Santos, and
  Toint]{birgin2017worst}
Birgin, E.~G., Gardenghi, J., Mart{\'\i}nez, J.~M., Santos, S.~A., and Toint,
  P.~L.
\newblock Worst-case evaluation complexity for unconstrained nonlinear
  optimization using high-order regularized models.
\newblock \emph{Mathematical Programming}, 163\penalty0 (1-2):\penalty0
  359--368, 2017.

\bibitem[Bullins(2018)]{bullins2018fast}
Bullins, B.
\newblock Fast minimization of structured convex quartics.
\newblock \emph{arXiv preprint arXiv:1812.10349}, 2018.

\bibitem[Carmon \& Duchi(2019)Carmon and Duchi]{carmon2019gradient}
Carmon, Y. and Duchi, J.
\newblock Gradient descent finds the cubic-regularized nonconvex {N}ewton step.
\newblock \emph{SIAM Journal on Optimization}, 29\penalty0 (3):\penalty0
  2146--2178, 2019.

\bibitem[Cartis \& Scheinberg(2018)Cartis and Scheinberg]{cartis2018global}
Cartis, C. and Scheinberg, K.
\newblock Global convergence rate analysis of unconstrained optimization
  methods based on probabilistic models.
\newblock \emph{Mathematical Programming}, 169\penalty0 (2):\penalty0 337--375,
  2018.

\bibitem[Cartis et~al.(2011{\natexlab{a}})Cartis, Gould, and
  Toint]{cartis2011adaptive1}
Cartis, C., Gould, N.~I., and Toint, P.~L.
\newblock Adaptive cubic regularisation methods for unconstrained optimization.
  {P}art {I}: motivation, convergence and numerical results.
\newblock \emph{Mathematical Programming}, 127\penalty0 (2):\penalty0 245--295,
  2011{\natexlab{a}}.

\bibitem[Cartis et~al.(2011{\natexlab{b}})Cartis, Gould, and
  Toint]{cartis2011adaptive2}
Cartis, C., Gould, N.~I., and Toint, P.~L.
\newblock Adaptive cubic regularisation methods for unconstrained optimization.
  {P}art {II}: worst-case function-and derivative-evaluation complexity.
\newblock \emph{Mathematical programming}, 130\penalty0 (2):\penalty0 295--319,
  2011{\natexlab{b}}.

\bibitem[Cartis et~al.(2019)Cartis, Gould, and Toint]{cartis2019universal}
Cartis, C., Gould, N.~I., and Toint, P.~L.
\newblock Universal regularization methods: varying the power, the smoothness
  and the accuracy.
\newblock \emph{SIAM Journal on Optimization}, 29\penalty0 (1):\penalty0
  595--615, 2019.

\bibitem[Doikov \& Nesterov(2019{\natexlab{a}})Doikov and
  Nesterov]{doikov2019contracting}
Doikov, N. and Nesterov, Y.
\newblock Contracting proximal methods for smooth convex optimization.
\newblock \emph{CORE Discussion Papers 2019/27}, 2019{\natexlab{a}}.

\bibitem[Doikov \& Nesterov(2019{\natexlab{b}})Doikov and
  Nesterov]{doikov2019local}
Doikov, N. and Nesterov, Y.
\newblock Local convergence of tensor methods.
\newblock \emph{CORE Discussion Papers 2019/21}, 2019{\natexlab{b}}.

\bibitem[Doikov \& Nesterov(2019{\natexlab{c}})Doikov and
  Nesterov]{doikov2019minimizing}
Doikov, N. and Nesterov, Y.
\newblock Minimizing uniformly convex functions by cubic regularization of
  {N}ewton method.
\newblock \emph{arXiv preprint arXiv:1905.02671}, 2019{\natexlab{c}}.

\bibitem[Doikov \& Richt{\'a}rik(2018)Doikov and
  Richt{\'a}rik]{doikov2018randomized}
Doikov, N. and Richt{\'a}rik, P.
\newblock Randomized block cubic {N}ewton method.
\newblock In \emph{International Conference on Machine Learning}, pp.\
  1289--1297, 2018.

\bibitem[Gasnikov et~al.(2019)Gasnikov, Dvurechensky, Gorbunov, Vorontsova,
  Selikhanovych, Uribe, Jiang, Wang, Zhang, Bubeck, et~al.]{gasnikov2019near}
Gasnikov, A., Dvurechensky, P., Gorbunov, E., Vorontsova, E., Selikhanovych,
  D., Uribe, C.~A., Jiang, B., Wang, H., Zhang, S., Bubeck, S., et~al.
\newblock Near optimal methods for minimizing convex functions with lipschitz $
  p $-th derivatives.
\newblock In \emph{Conference on Learning Theory}, pp.\  1392--1393, 2019.

\bibitem[Gould et~al.(2010)Gould, Robinson, and Thorne]{gould2010solving}
Gould, N.~I., Robinson, D.~P., and Thorne, H.~S.
\newblock On solving trust-region and other regularised subproblems in
  optimization.
\newblock \emph{Mathematical Programming Computation}, 2\penalty0 (1):\penalty0
  21--57, 2010.

\bibitem[Grapiglia \& Nesterov(2019{\natexlab{a}})Grapiglia and
  Nesterov]{grapiglia2019accelerated}
Grapiglia, G.~N. and Nesterov, Y.
\newblock Accelerated regularized {N}ewton methods for minimizing composite
  convex functions.
\newblock \emph{SIAM Journal on Optimization}, 29\penalty0 (1):\penalty0
  77--99, 2019{\natexlab{a}}.

\bibitem[Grapiglia \& Nesterov(2019{\natexlab{b}})Grapiglia and
  Nesterov]{grapiglia2019inexact}
Grapiglia, G.~N. and Nesterov, Y.
\newblock On inexact solution of auxiliary problems in tensor methods for
  convex optimization.
\newblock \emph{arXiv preprint arXiv:1907.13023}, 2019{\natexlab{b}}.

\bibitem[Grapiglia \& Nesterov(2019{\natexlab{c}})Grapiglia and
  Nesterov]{grapiglia2019tensor}
Grapiglia, G.~N. and Nesterov, Y.
\newblock Tensor methods for minimizing functions with {H}\"{o}lder continuous
  higher-order derivatives.
\newblock \emph{arXiv preprint arXiv:1904.12559}, 2019{\natexlab{c}}.

\bibitem[Grapiglia \& Nesterov(2019{\natexlab{d}})Grapiglia and
  Nesterov]{grapiglia2019tensor2}
Grapiglia, G.~N. and Nesterov, Y.
\newblock Tensor methods for finding approximate stationary points of convex
  functions.
\newblock \emph{arXiv preprint arXiv:1907.07053}, 2019{\natexlab{d}}.

\bibitem[G{\"u}ler(1992)]{guler1992new}
G{\"u}ler, O.
\newblock New proximal point algorithms for convex minimization.
\newblock \emph{SIAM Journal on Optimization}, 2\penalty0 (4):\penalty0
  649--664, 1992.

\bibitem[Ivanova et~al.(2019)Ivanova, Grishchenko, Gasnikov, and
  Shulgin]{ivanova2019adaptive}
Ivanova, A., Grishchenko, D., Gasnikov, A., and Shulgin, E.
\newblock Adaptive catalyst for smooth convex optimization.
\newblock \emph{arXiv preprint arXiv:1911.11271}, 2019.

\bibitem[Jiang et~al.(2018)Jiang, Lin, and Zhang]{jiang2018unified}
Jiang, B., Lin, T., and Zhang, S.
\newblock A unified adaptive tensor approximation scheme to accelerate
  composite convex optimization.
\newblock \emph{arXiv preprint arXiv:1811.02427}, 2018.

\bibitem[Kohler \& Lucchi(2017)Kohler and Lucchi]{kohler2017sub}
Kohler, J.~M. and Lucchi, A.
\newblock Sub-sampled cubic regularization for non-convex optimization.
\newblock In \emph{International Conference on Machine Learning}, pp.\
  1895--1904, 2017.

\bibitem[Kulunchakov \& Mairal(2019)Kulunchakov and
  Mairal]{kulunchakov2019generic}
Kulunchakov, A. and Mairal, J.
\newblock A generic acceleration framework for stochastic composite
  optimization.
\newblock In \emph{Advances in Neural Information Processing Systems}, pp.\
  12556--12567, 2019.

\bibitem[Lin et~al.(2015)Lin, Mairal, and Harchaoui]{lin2015universal}
Lin, H., Mairal, J., and Harchaoui, Z.
\newblock A universal catalyst for first-order optimization.
\newblock In \emph{Advances in Neural Information Processing Systems}, pp.\
  3384--3392, 2015.

\bibitem[Lin et~al.(2018)Lin, Mairal, and Harchaoui]{lin2018catalyst}
Lin, H., Mairal, J., and Harchaoui, Z.
\newblock Catalyst acceleration for first-order convex optimization: from
  theory to practice.
\newblock \emph{Journal of Machine Learning Research}, 18\penalty0
  (212):\penalty0 1--54, 2018.

\bibitem[Lu et~al.(2018)Lu, Freund, and Nesterov]{lu2018relatively}
Lu, H., Freund, R.~M., and Nesterov, Y.
\newblock Relatively smooth convex optimization by first-order methods, and
  applications.
\newblock \emph{SIAM Journal on Optimization}, 28\penalty0 (1):\penalty0
  333--354, 2018.

\bibitem[Lucchi \& Kohler(2019)Lucchi and Kohler]{lucchi2019stochastic}
Lucchi, A. and Kohler, J.
\newblock A stochastic tensor method for non-convex optimization.
\newblock \emph{arXiv preprint arXiv:1911.10367}, 2019.

\bibitem[Monteiro \& Svaiter(2013)Monteiro and
  Svaiter]{monteiro2013accelerated}
Monteiro, R.~D. and Svaiter, B.~F.
\newblock An accelerated hybrid proximal extragradient method for convex
  optimization and its implications to second-order methods.
\newblock \emph{SIAM Journal on Optimization}, 23\penalty0 (2):\penalty0
  1092--1125, 2013.

\bibitem[Nesterov(1983)]{nesterov1983method}
Nesterov, Y.
\newblock A method for solving the convex programming problem with convergence
  rate {O}(1/k\^{}2).
\newblock In \emph{Dokl. akad. nauk Sssr}, volume 269, pp.\  543--547, 1983.

\bibitem[Nesterov(2008)]{nesterov2008accelerating}
Nesterov, Y.
\newblock Accelerating the cubic regularization of {N}ewton's method on convex
  problems.
\newblock \emph{Mathematical Programming}, 112\penalty0 (1):\penalty0 159--181,
  2008.

\bibitem[Nesterov(2013)]{nesterov2013gradient}
Nesterov, Y.
\newblock Gradient methods for minimizing composite functions.
\newblock \emph{Mathematical Programming}, 140\penalty0 (1):\penalty0 125--161,
  2013.

\bibitem[Nesterov(2018)]{nesterov2018lectures}
Nesterov, Y.
\newblock \emph{Lectures on convex optimization}, volume 137.
\newblock Springer, 2018.

\bibitem[Nesterov(2019{\natexlab{a}})]{nesterov2019implementable}
Nesterov, Y.
\newblock Implementable tensor methods in unconstrained convex optimization.
\newblock \emph{Mathematical Programming}, pp.\  1--27, 2019{\natexlab{a}}.

\bibitem[Nesterov(2019{\natexlab{b}})]{nesterov2019inexact}
Nesterov, Y.
\newblock Inexact basic tensor methods.
\newblock \emph{CORE Discussion Papers 2019/23}, 2019{\natexlab{b}}.

\bibitem[Nesterov \& Nemirovskii(1994)Nesterov and
  Nemirovskii]{nesterov1994interior}
Nesterov, Y. and Nemirovskii, A.
\newblock \emph{Interior-point polynomial algorithms in convex programming}.
\newblock SIAM, 1994.

\bibitem[Nesterov \& Polyak(2006)Nesterov and Polyak]{nesterov2006cubic}
Nesterov, Y. and Polyak, B.~T.
\newblock Cubic regularization of {N}ewton's method and its global performance.
\newblock \emph{Mathematical Programming}, 108\penalty0 (1):\penalty0 177--205,
  2006.

\bibitem[Nocedal \& Wright(2006)Nocedal and Wright]{nocedal2006numerical}
Nocedal, J. and Wright, S.~J.
\newblock Numerical optimization 2nd, 2006.

\bibitem[Rodomanov \& Nesterov(2019)Rodomanov and
  Nesterov]{rodomanov2019smoothness}
Rodomanov, A. and Nesterov, Y.
\newblock Smoothness parameter of power of euclidean norm.
\newblock \emph{arXiv preprint arXiv:1907.12346}, 2019.

\bibitem[Schmidt et~al.(2011)Schmidt, Roux, and Bach]{schmidt2011convergence}
Schmidt, M., Roux, N.~L., and Bach, F.~R.
\newblock Convergence rates of inexact proximal-gradient methods for convex
  optimization.
\newblock In \emph{Advances in neural information processing systems}, pp.\
  1458--1466, 2011.

\bibitem[Song \& Ma(2019)Song and Ma]{song2019towards}
Song, C. and Ma, Y.
\newblock Towards unified acceleration of high-order algorithms under
  {H}\"{o}lder continuity and uniform convexity.
\newblock \emph{arXiv preprint arXiv:1906.00582}, 2019.

\bibitem[Tripuraneni et~al.(2018)Tripuraneni, Stern, Jin, Regier, and
  Jordan]{tripuraneni2018stochastic}
Tripuraneni, N., Stern, M., Jin, C., Regier, J., and Jordan, M.~I.
\newblock Stochastic cubic regularization for fast nonconvex optimization.
\newblock In \emph{Advances in Neural Information Processing Systems}, pp.\
  2899--2908, 2018.

\bibitem[Van~Nguyen(2017)]{van2017forward}
Van~Nguyen, Q.
\newblock Forward-backward splitting with bregman distances.
\newblock \emph{Vietnam Journal of Mathematics}, 45\penalty0 (3):\penalty0
  519--539, 2017.

\bibitem[Wang et~al.(2018)Wang, Zhou, Liang, and Lan]{wang2018stochastic}
Wang, Z., Zhou, Y., Liang, Y., and Lan, G.
\newblock Stochastic variance-reduced cubic regularization for nonconvex
  optimization.
\newblock \emph{arXiv preprint arXiv:1802.07372}, 2018.

\bibitem[Zhou et~al.(2019)Zhou, Xu, and Gu]{zhou2019stochastic}
Zhou, D., Xu, P., and Gu, Q.
\newblock Stochastic variance-reduced cubic regularization methods.
\newblock \emph{Journal of Machine Learning Research}, 20\penalty0
  (134):\penalty0 1--47, 2019.

\end{thebibliography}
\bibliographystyle{icml2020}


\clearpage
\onecolumn
\appendix

\icmltitle{Supplementary Material}

\section{Extra Experiments}

\subsection{Exact Stopping Criterion}

\begin{figure}[h!]
	
	\begin{minipage}{0.33\textwidth}
		\centering
		\includegraphics[width=\textwidth ]{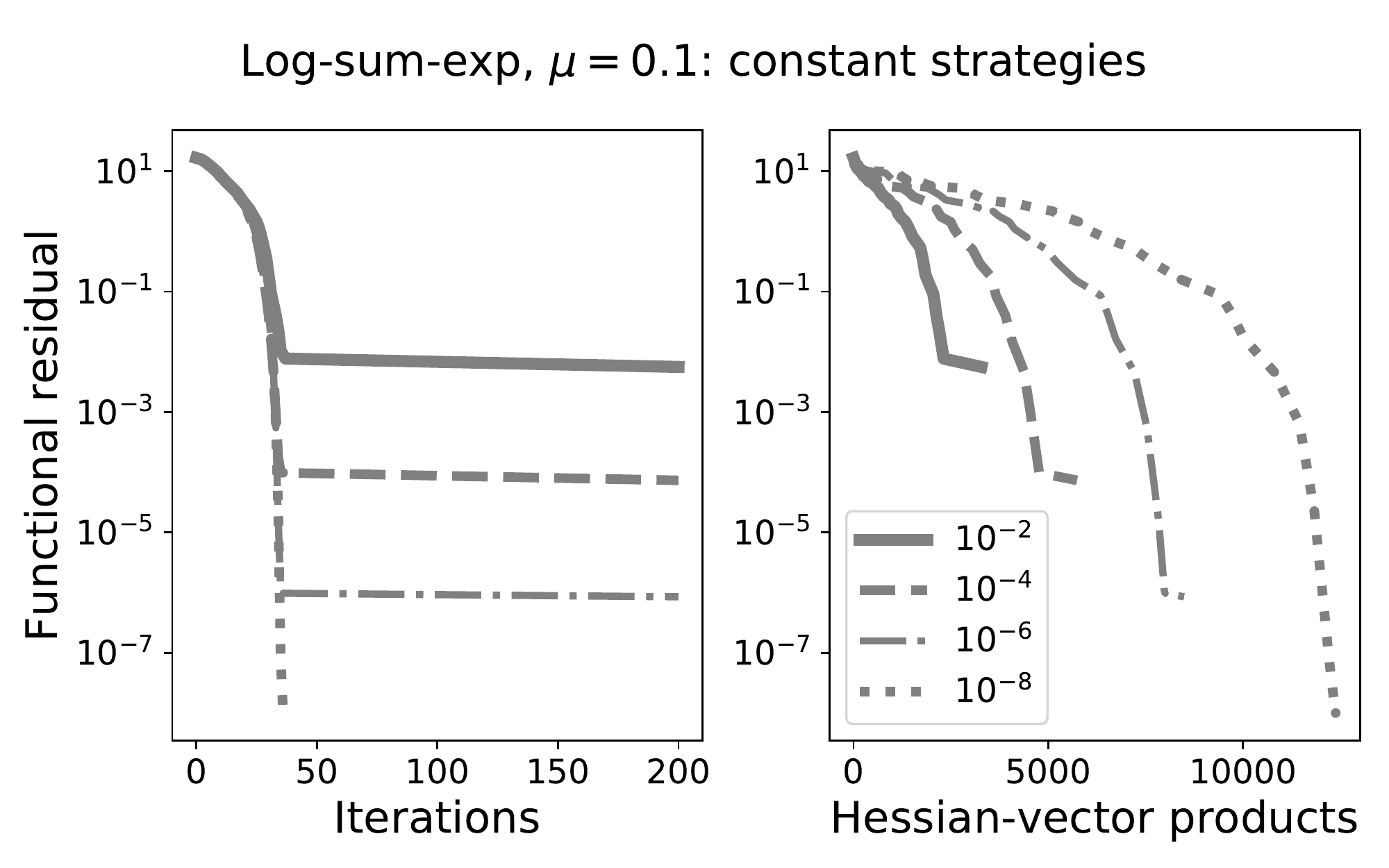}
	\end{minipage}
	\begin{minipage}{0.33\textwidth}
	\centering
	\includegraphics[width=\textwidth ]{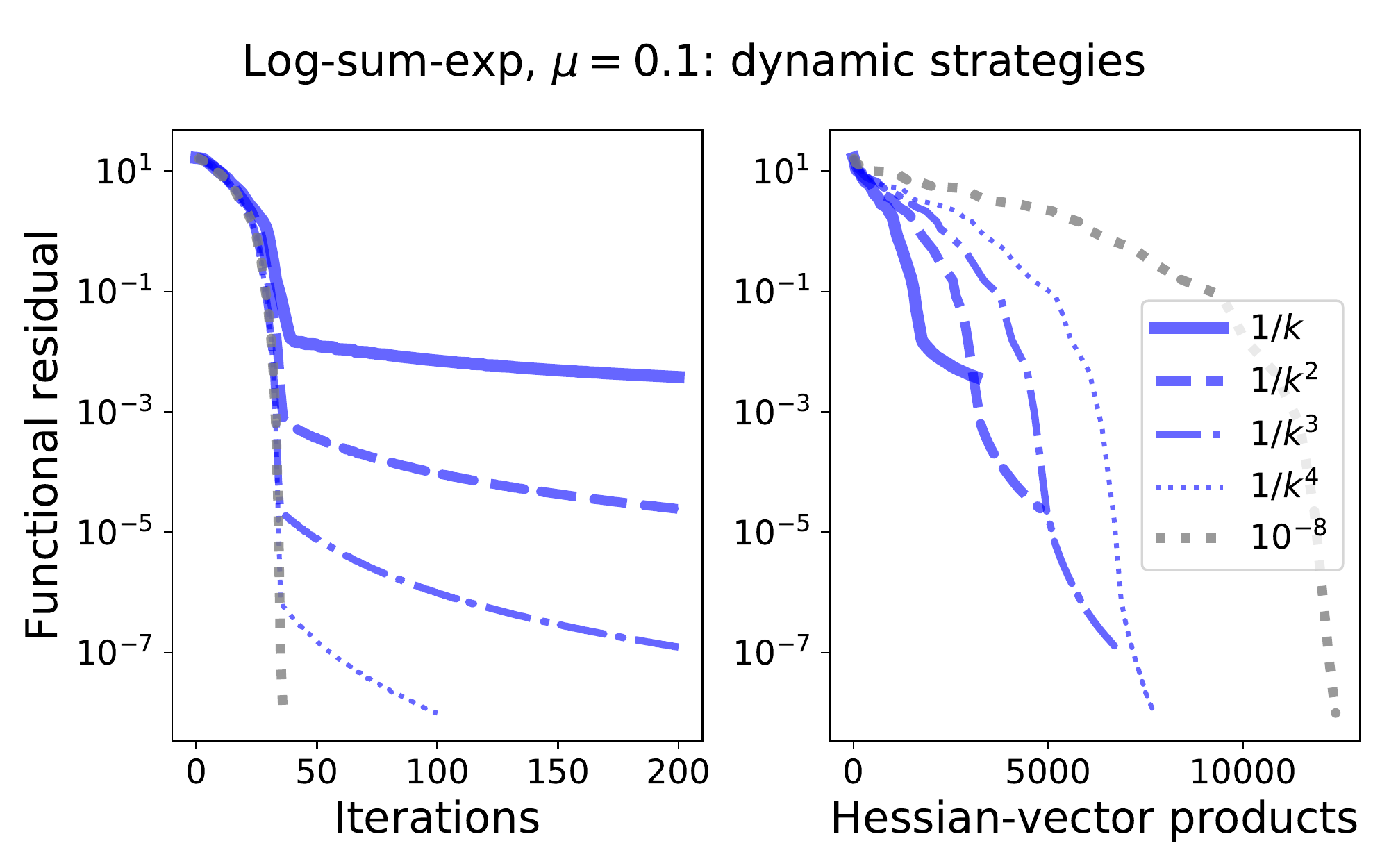}
	\end{minipage}
		\begin{minipage}{0.33\textwidth}
		\centering
		\includegraphics[width=\textwidth ]{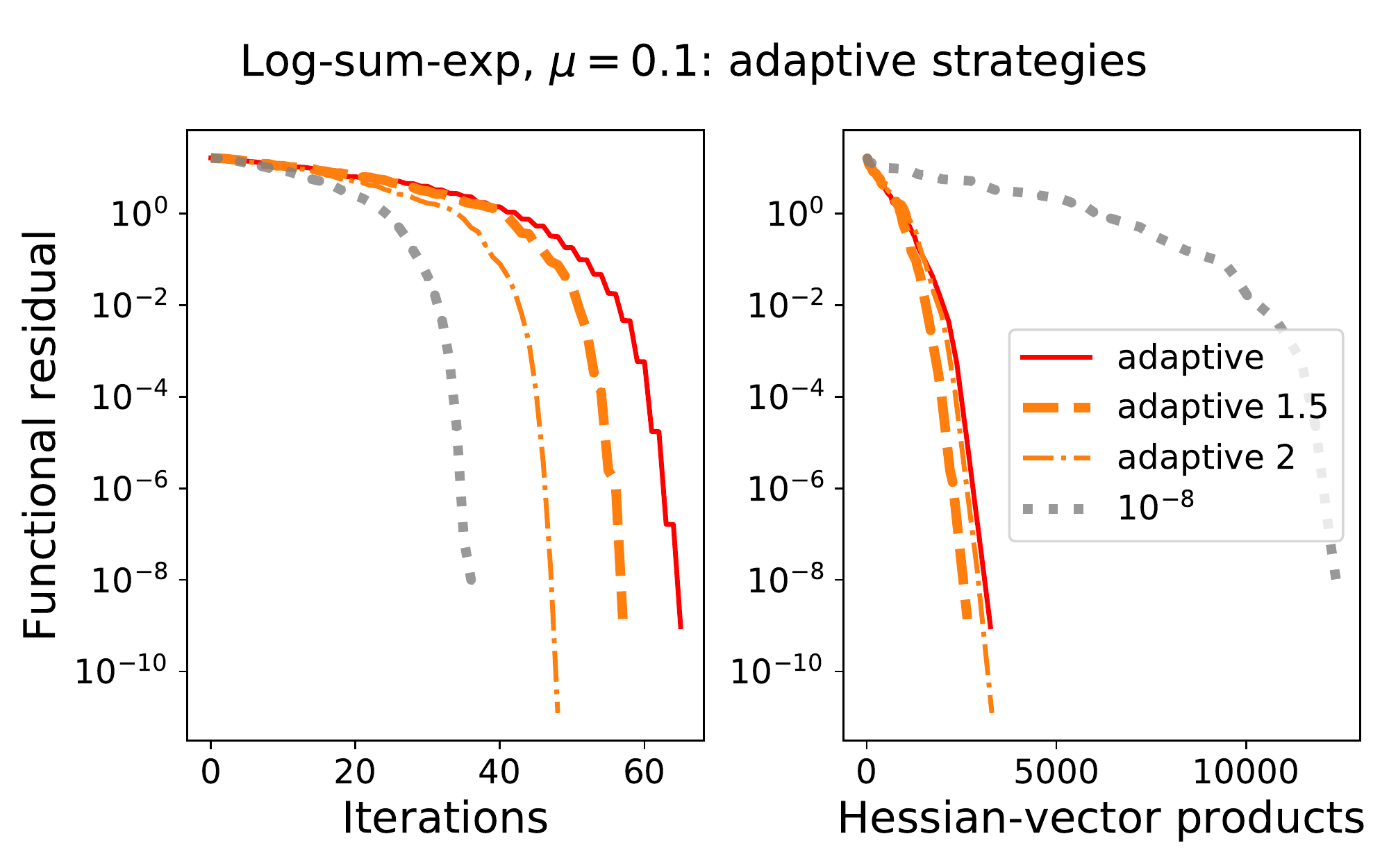}
	\end{minipage}

	\begin{minipage}{0.33\textwidth}
	\centering
	\includegraphics[width=\textwidth ]{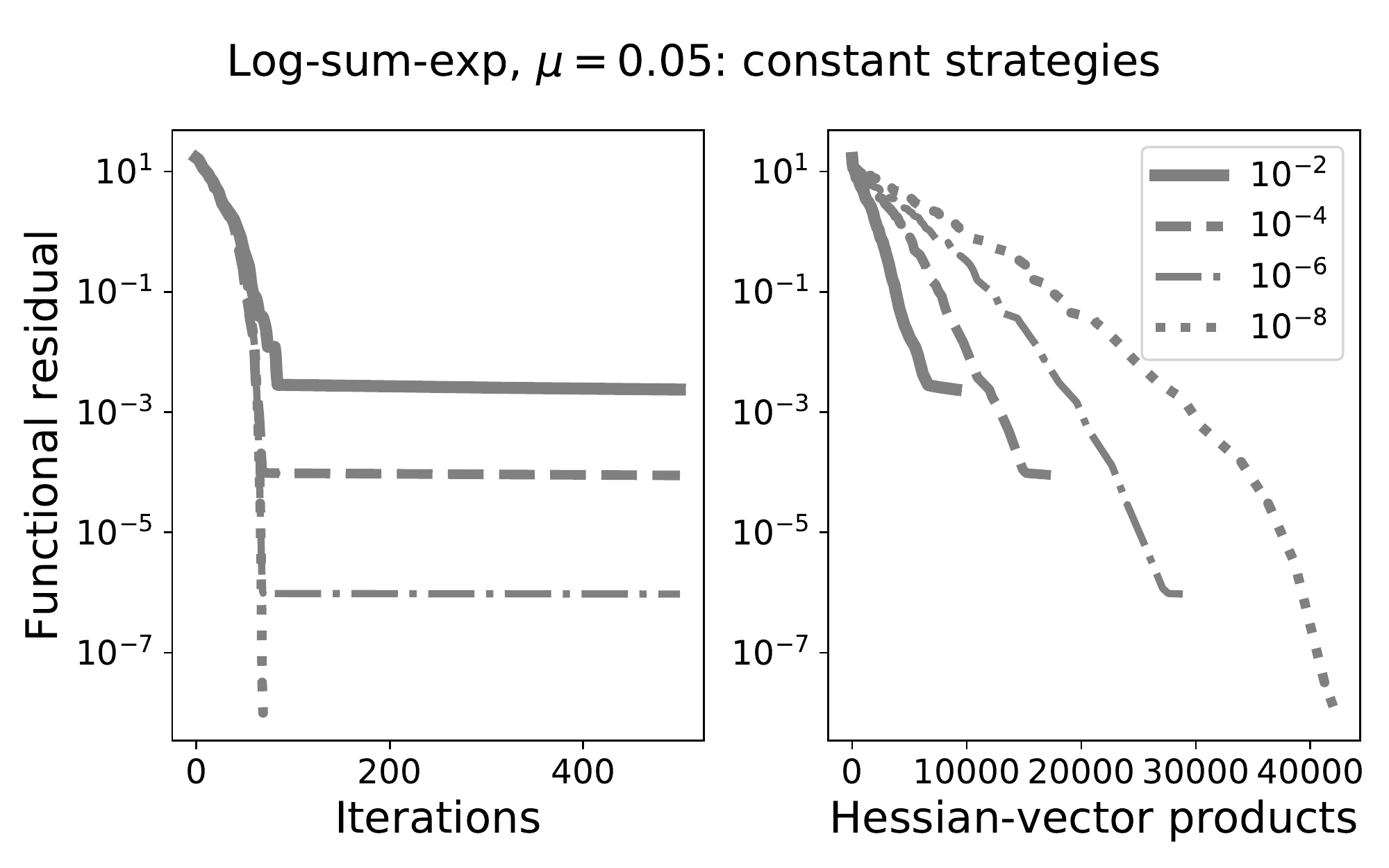}
	\end{minipage}
	\begin{minipage}{0.33\textwidth}
	\centering
	\includegraphics[width=\textwidth ]{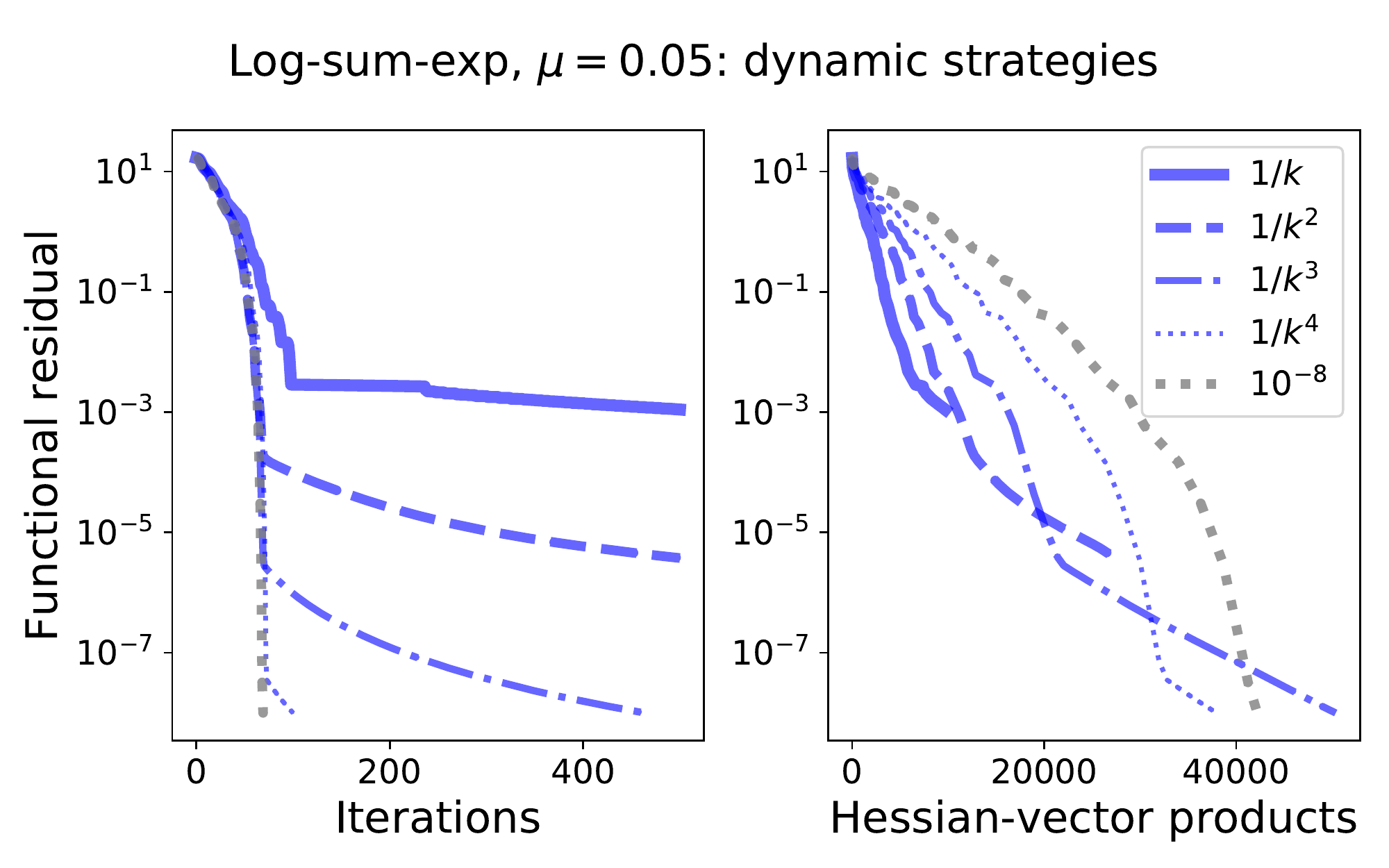}
	\end{minipage}
	\begin{minipage}{0.33\textwidth}
	\centering
	\includegraphics[width=\textwidth ]{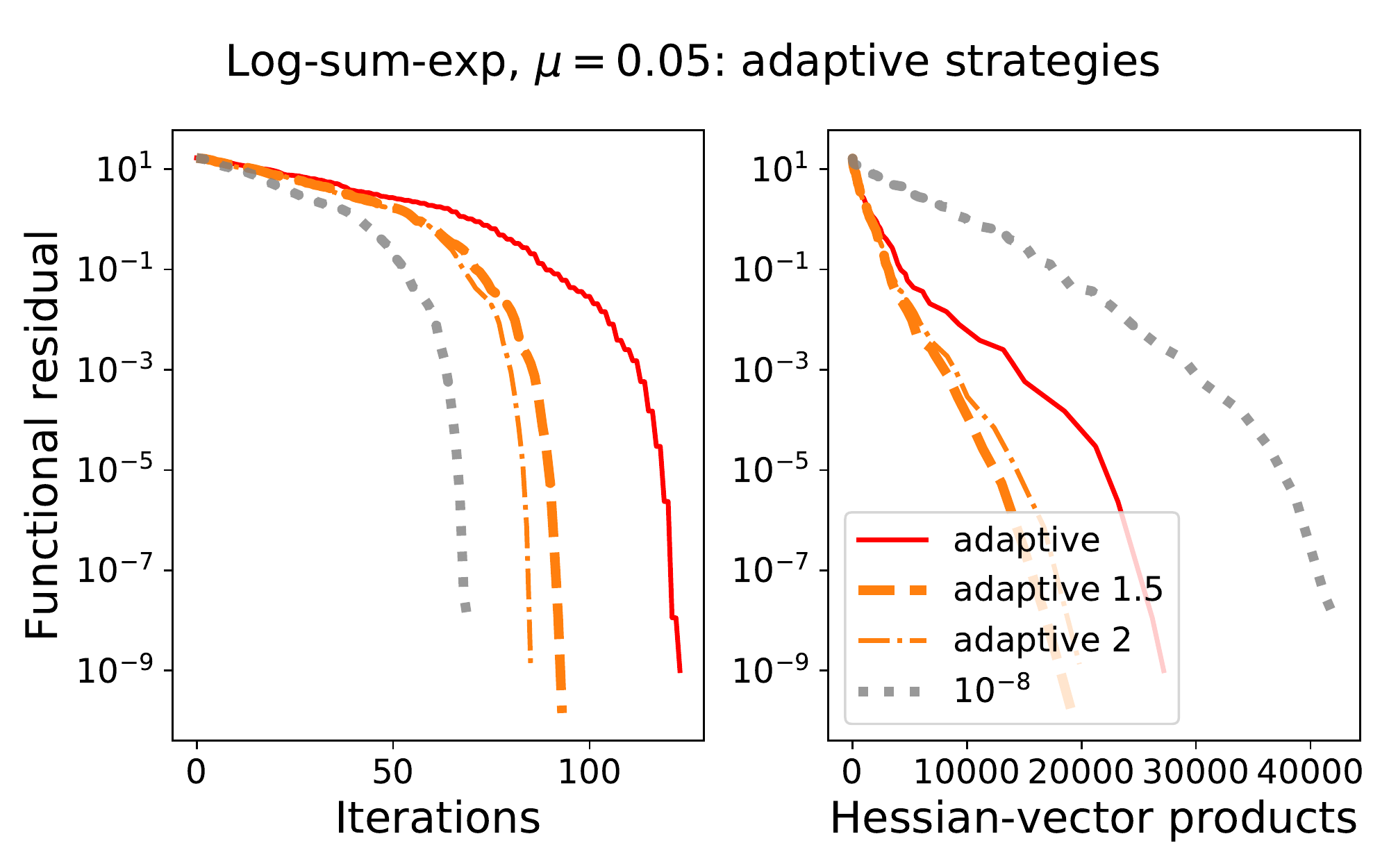}
	\end{minipage}

		\caption{Exact stopping criterion, minimizing Log-Sum-Exp function, $n = 100$.} 
		\label{fig:CN_log_sum_exp_exact_100}
\end{figure}

\begin{figure}[h!]
	
	\begin{minipage}{0.33\textwidth}
		\centering
		\includegraphics[width=\textwidth ]{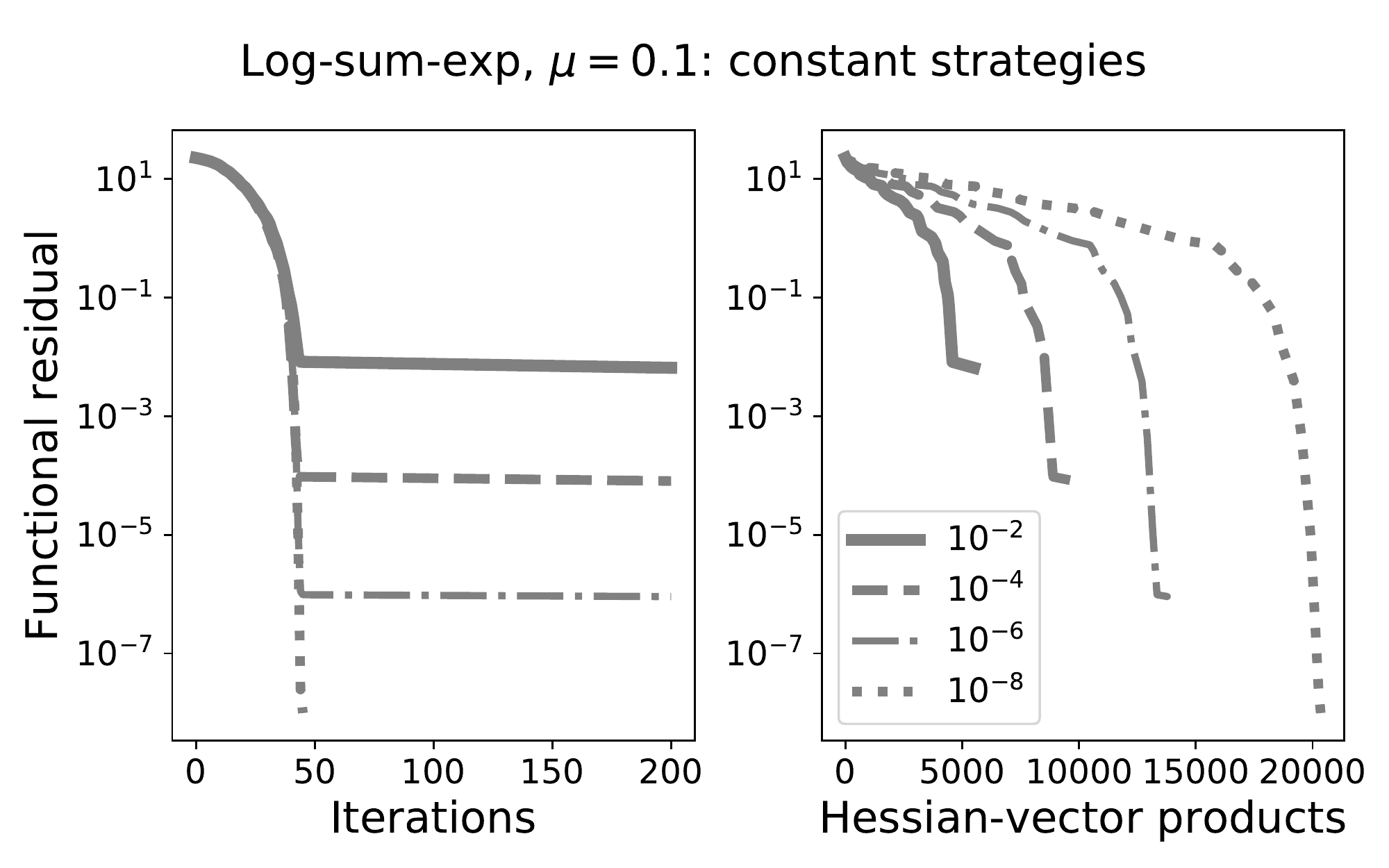}
	\end{minipage}
	\begin{minipage}{0.33\textwidth}
		\centering
		\includegraphics[width=\textwidth ]{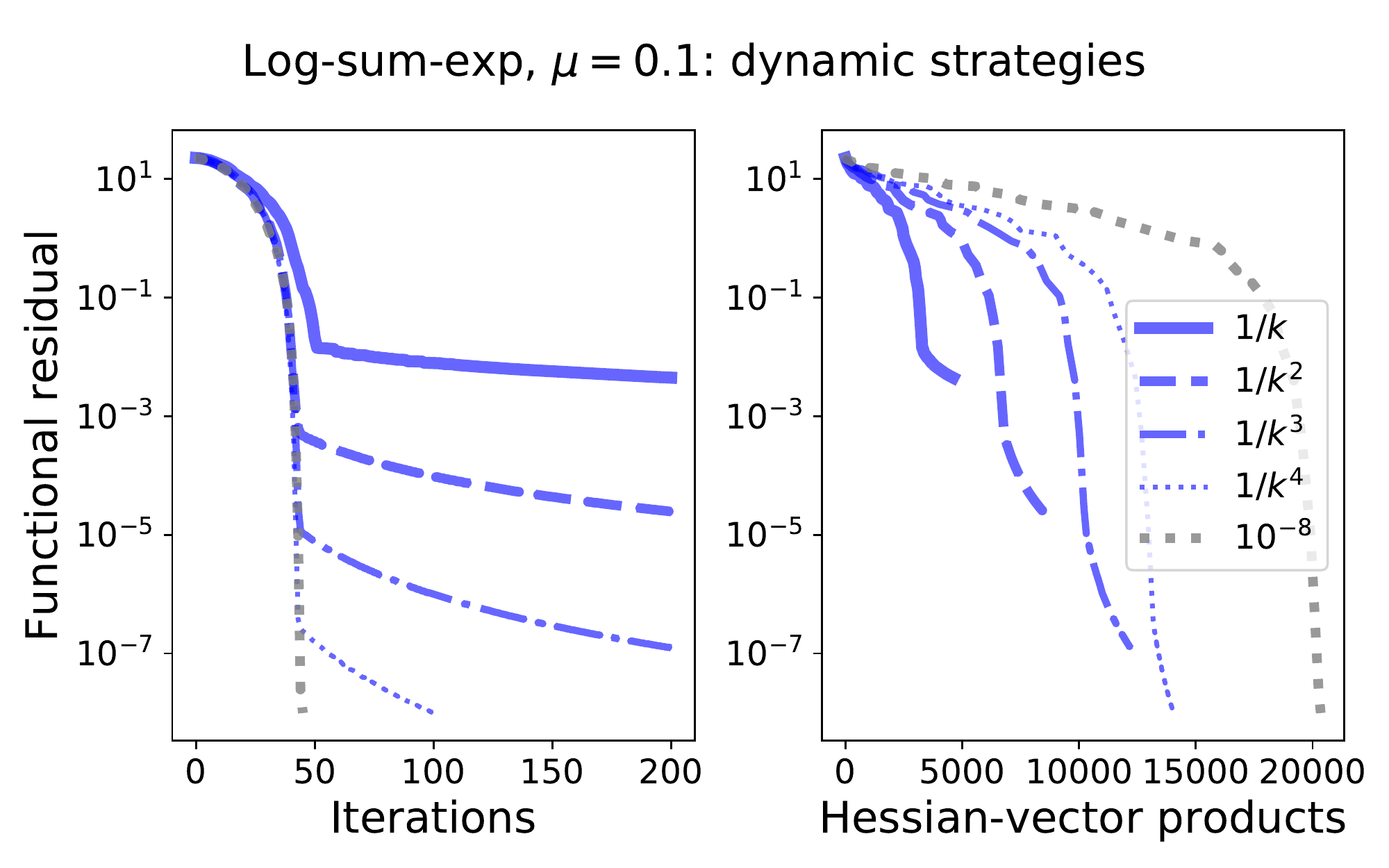}
	\end{minipage}
	\begin{minipage}{0.33\textwidth}
		\centering
		\includegraphics[width=\textwidth ]{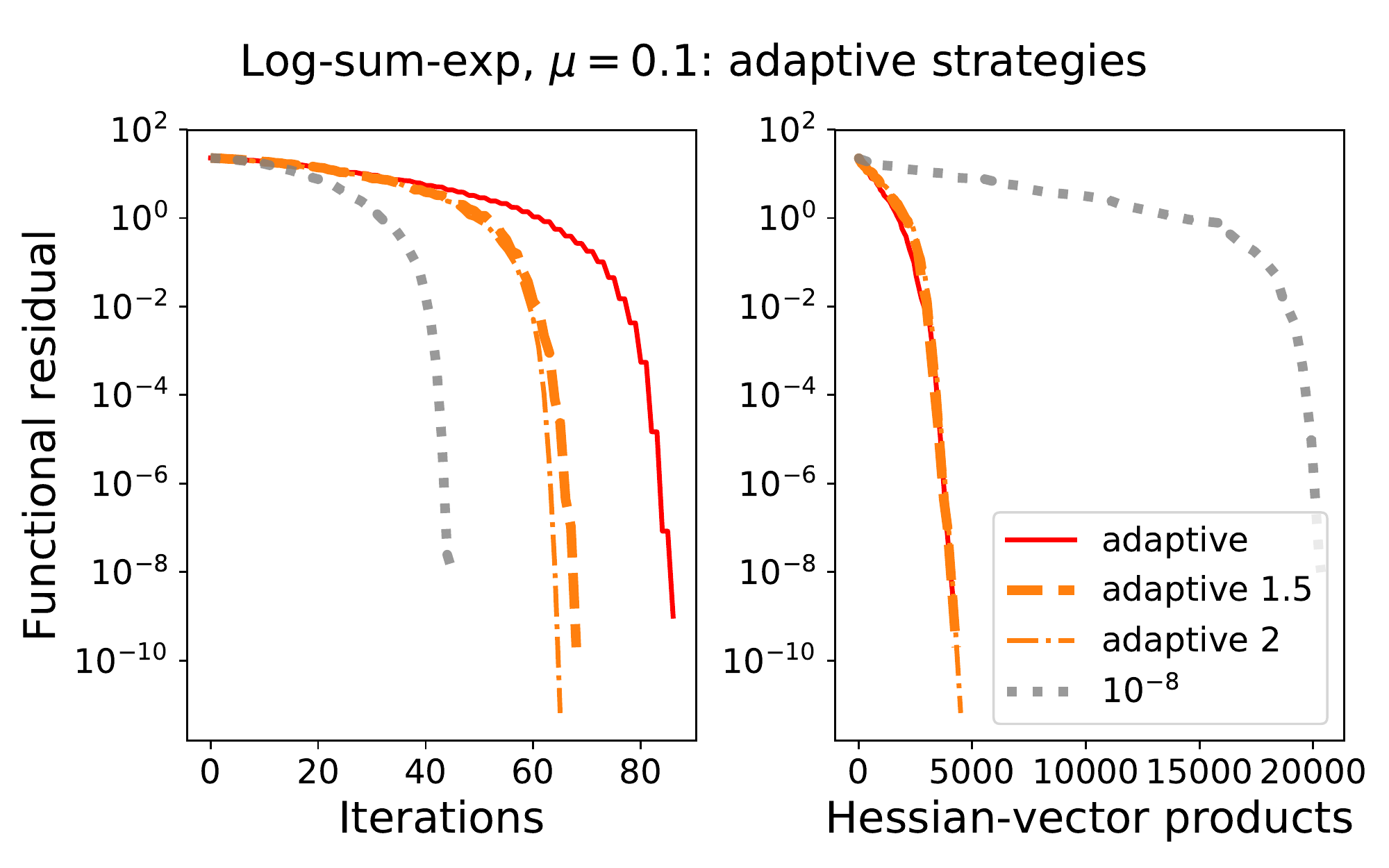}
	\end{minipage}
	
	\begin{minipage}{0.33\textwidth}
		\centering
		\includegraphics[width=\textwidth ]{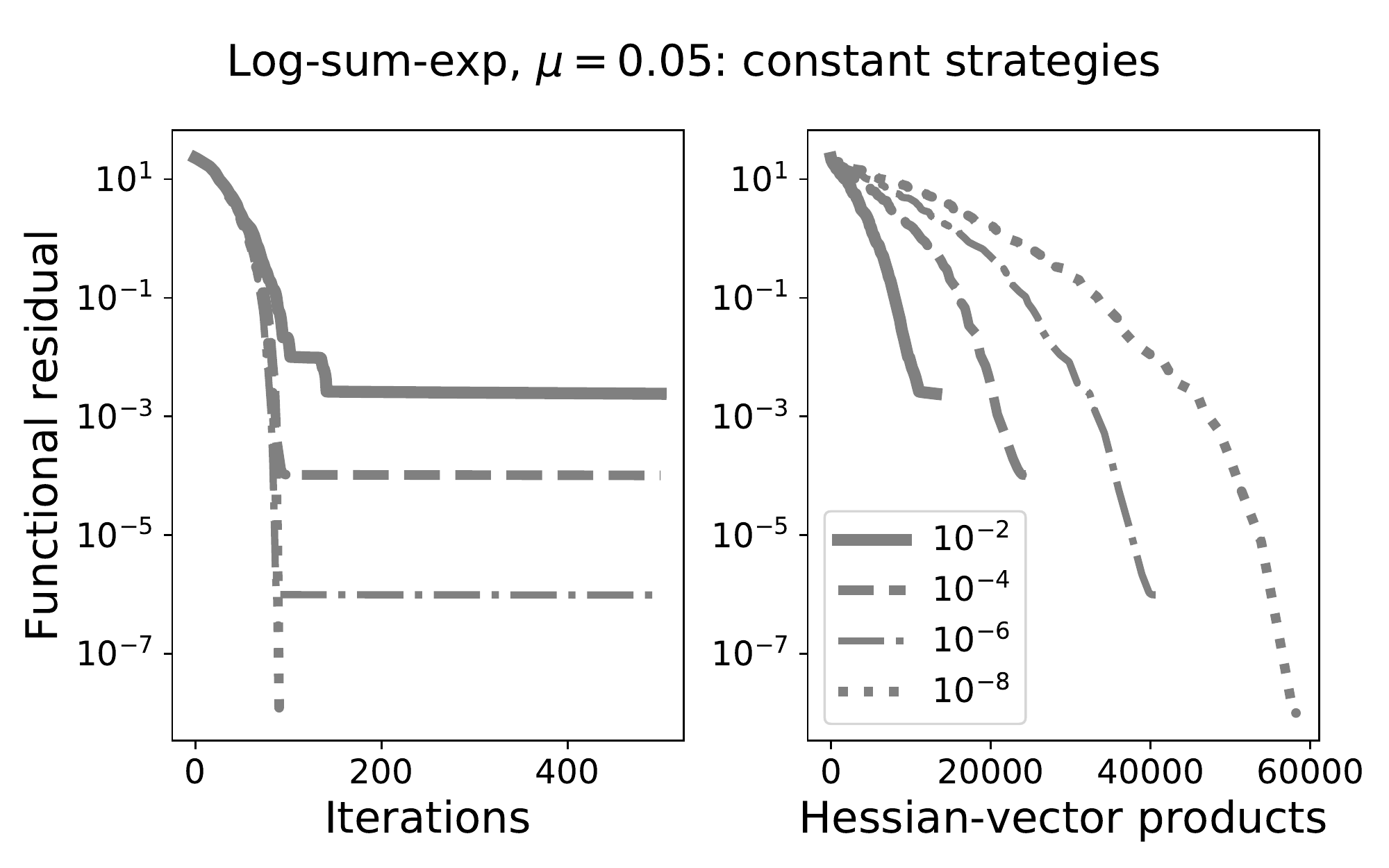}
	\end{minipage}
	\begin{minipage}{0.33\textwidth}
		\centering
		\includegraphics[width=\textwidth ]{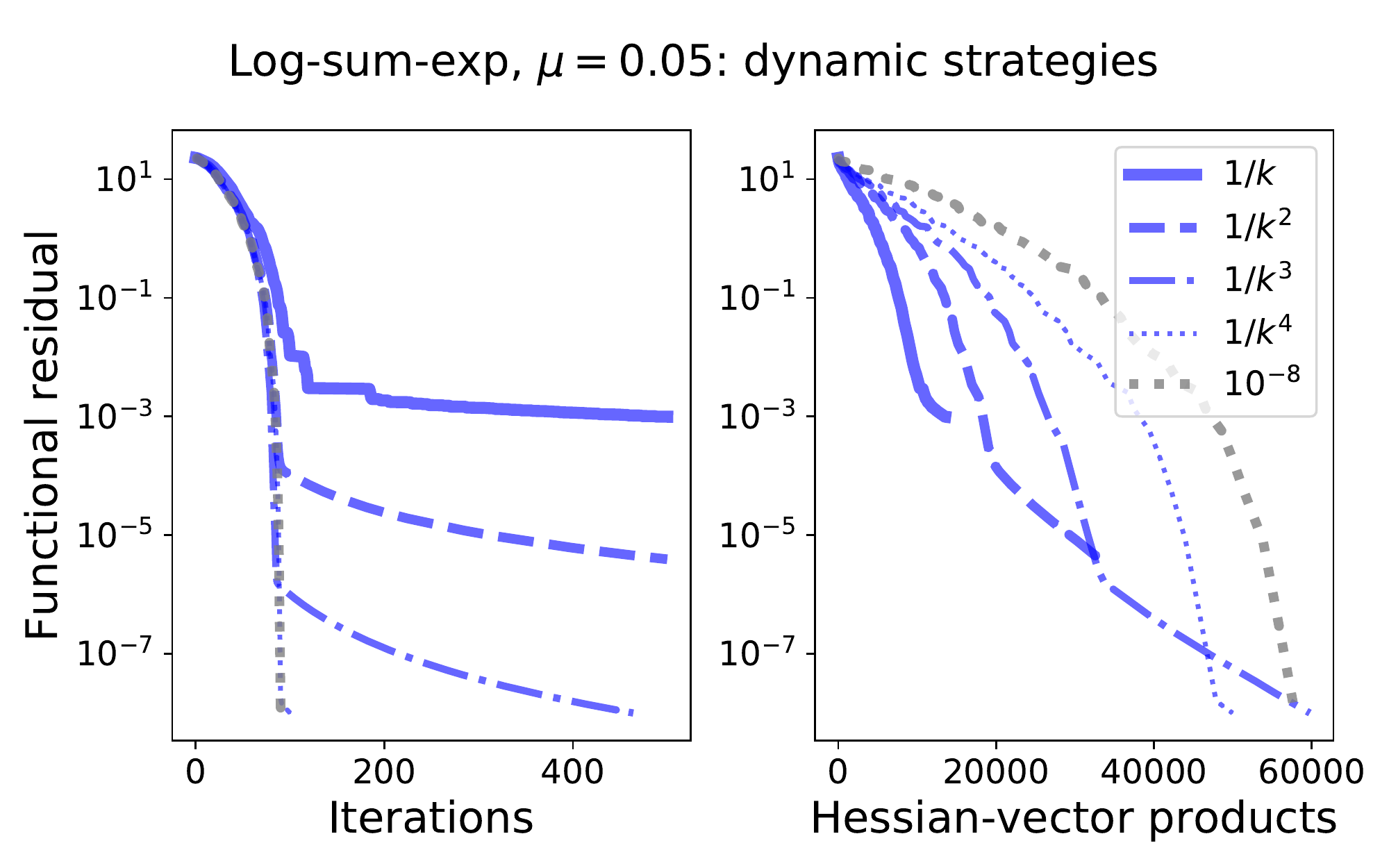}
	\end{minipage}
	\begin{minipage}{0.33\textwidth}
		\centering
		\includegraphics[width=\textwidth ]{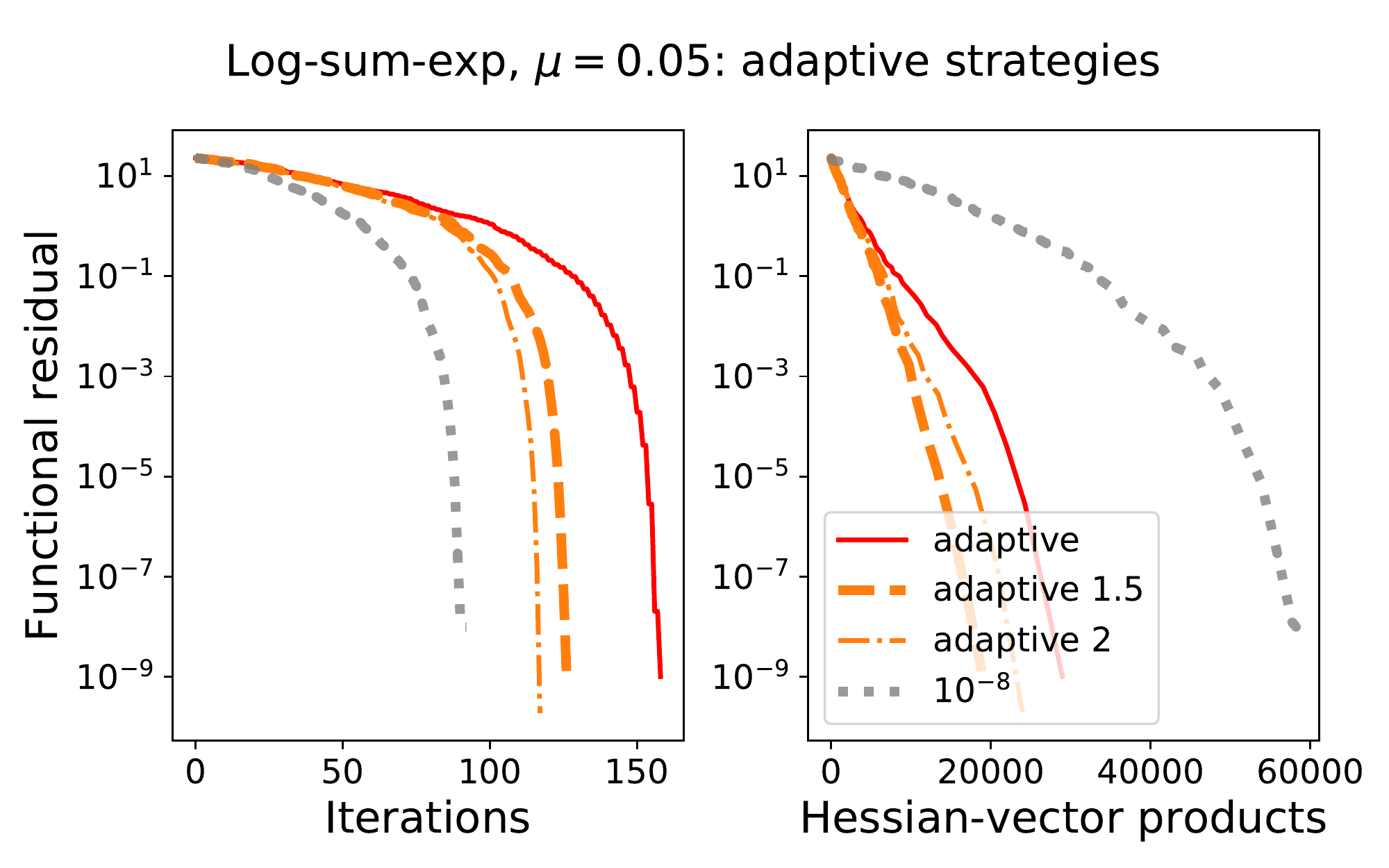}
	\end{minipage}
	
	\caption{Exact stopping criterion, minimizing Log-Sum-Exp function, $n = 200$.} 
	\label{fig:CN_log_sum_exp_exact_200}
\end{figure}

\begin{figure}[h!]	
	\begin{minipage}{0.33\textwidth}
		\centering
		\includegraphics[width=\textwidth ]{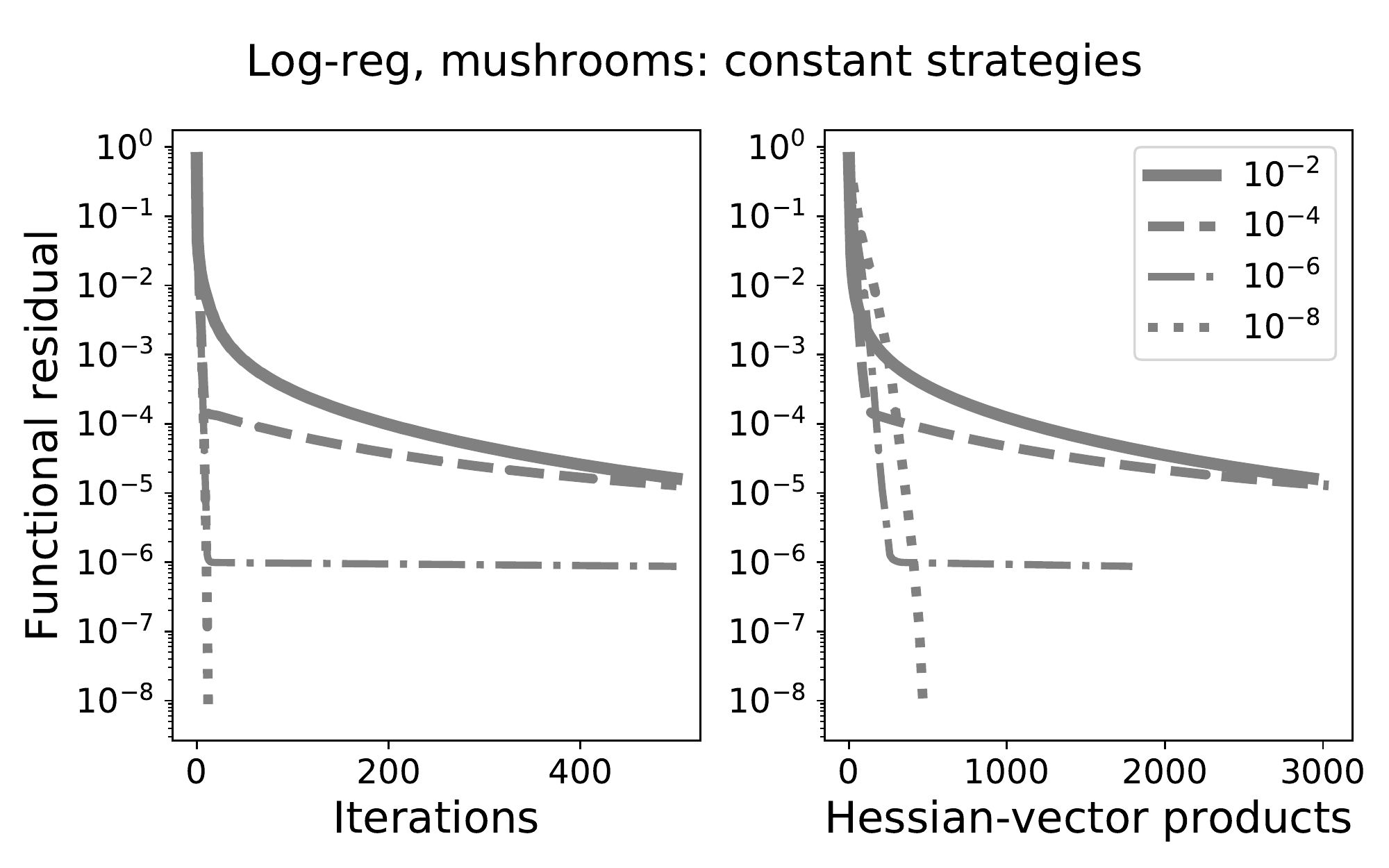}
	\end{minipage}
	\begin{minipage}{0.33\textwidth}
		\centering
		\includegraphics[width=\textwidth ]{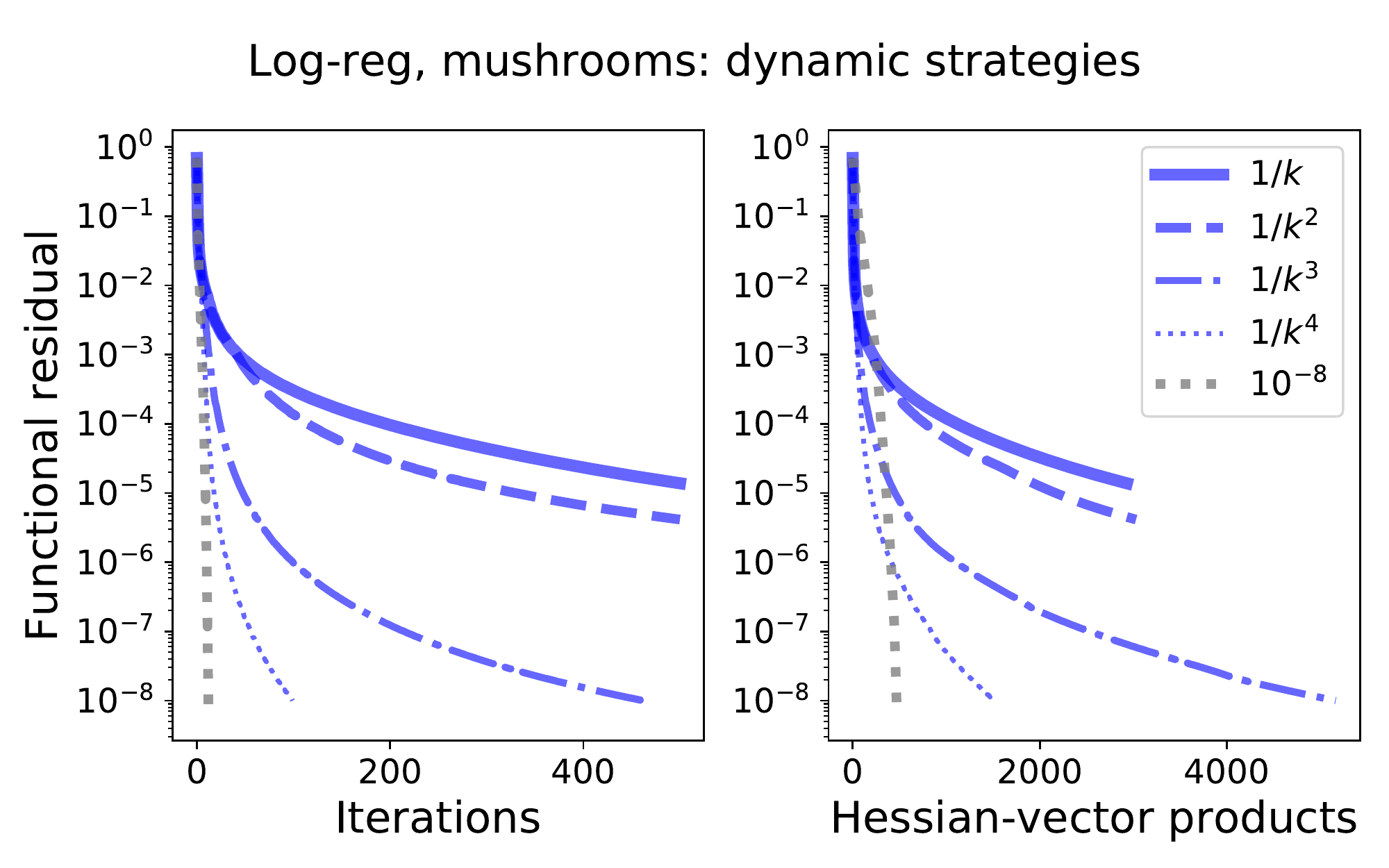}
	\end{minipage}
	\begin{minipage}{0.33\textwidth}
		\centering
		\includegraphics[width=\textwidth ]{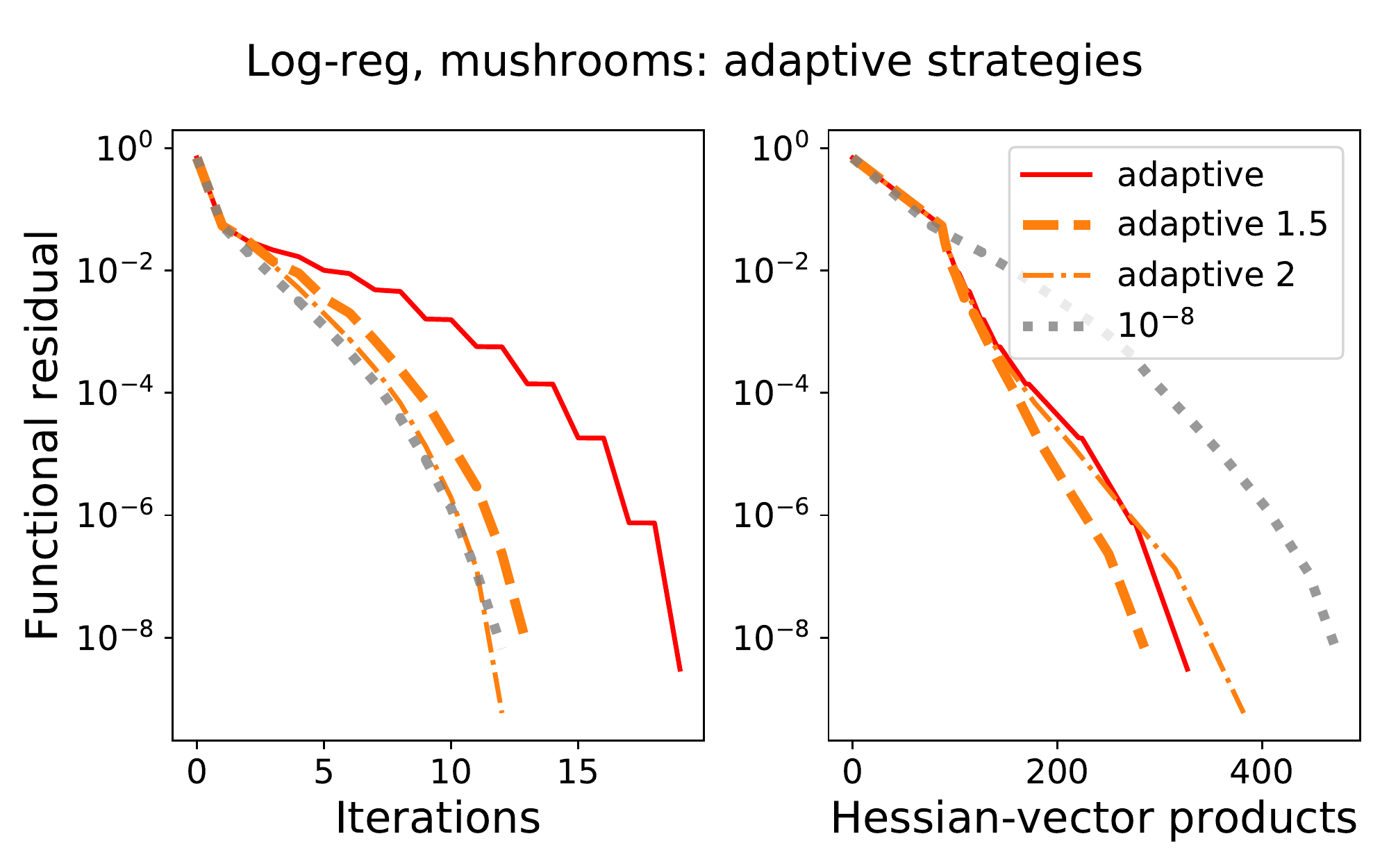}
	\end{minipage}
	
	\begin{minipage}{0.33\textwidth}
		\centering
		\includegraphics[width=\textwidth ]{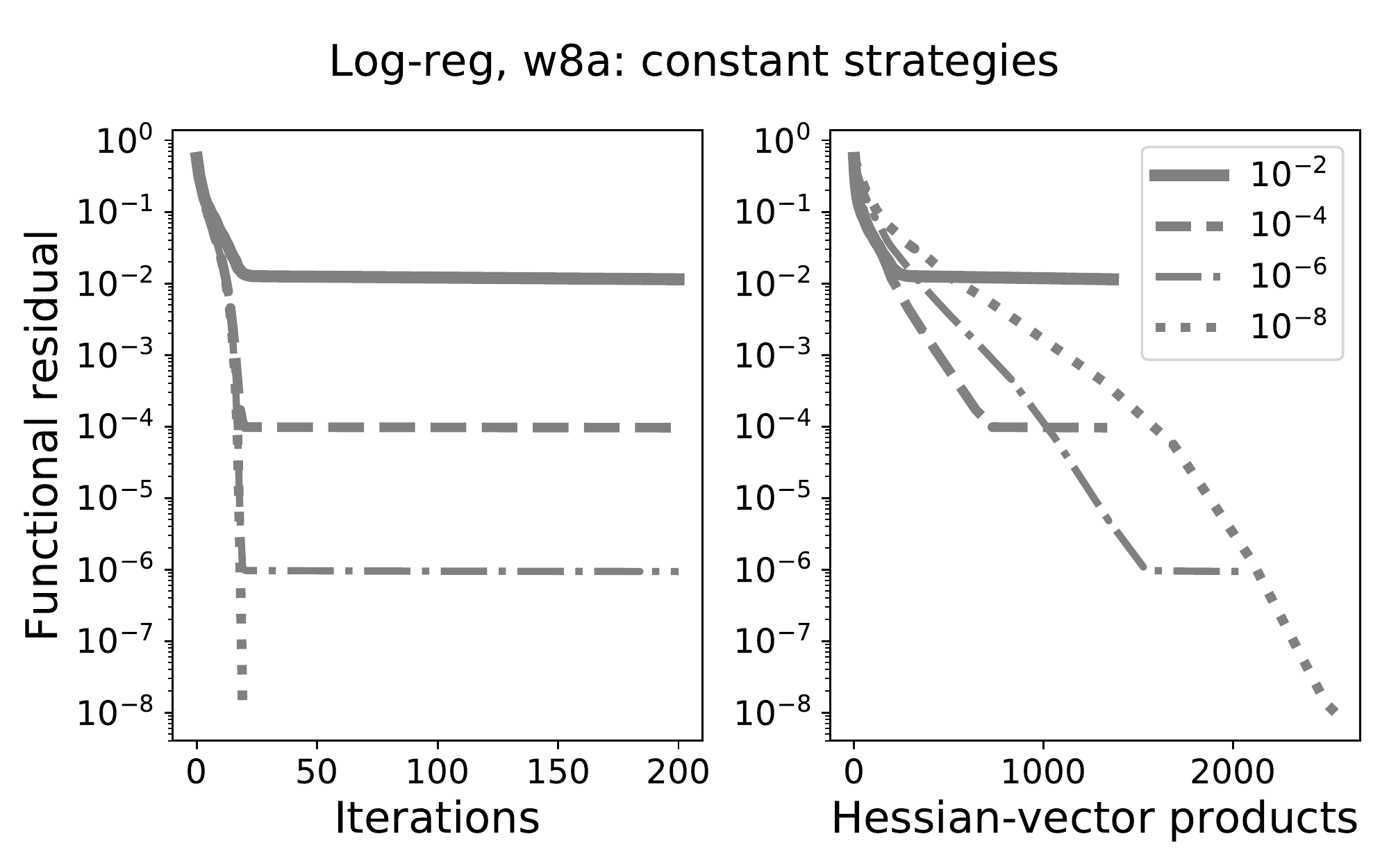}
	\end{minipage}
	\begin{minipage}{0.33\textwidth}
		\centering
		\includegraphics[width=\textwidth ]{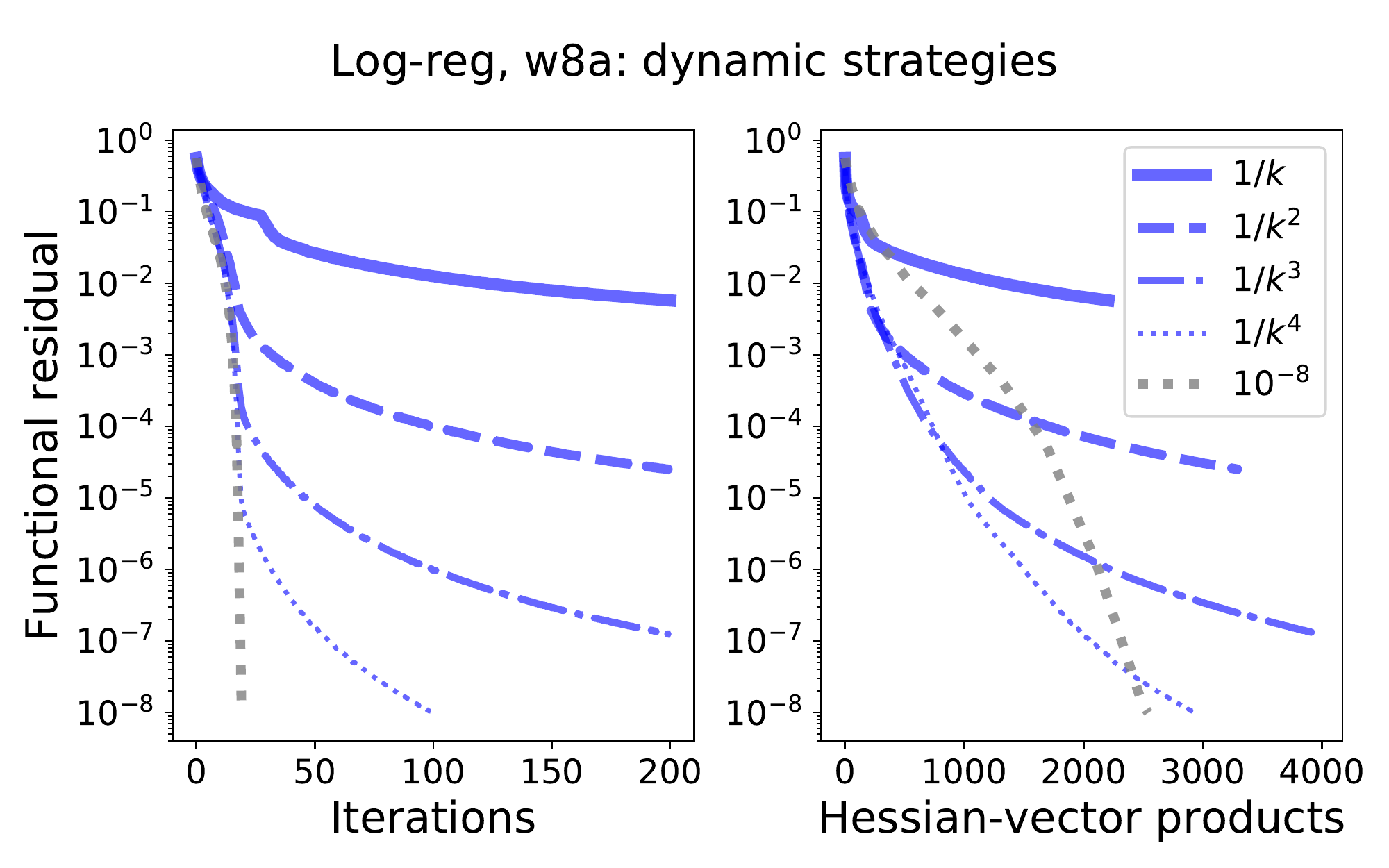}
	\end{minipage}
	\begin{minipage}{0.33\textwidth}
		\centering
		\includegraphics[width=\textwidth ]{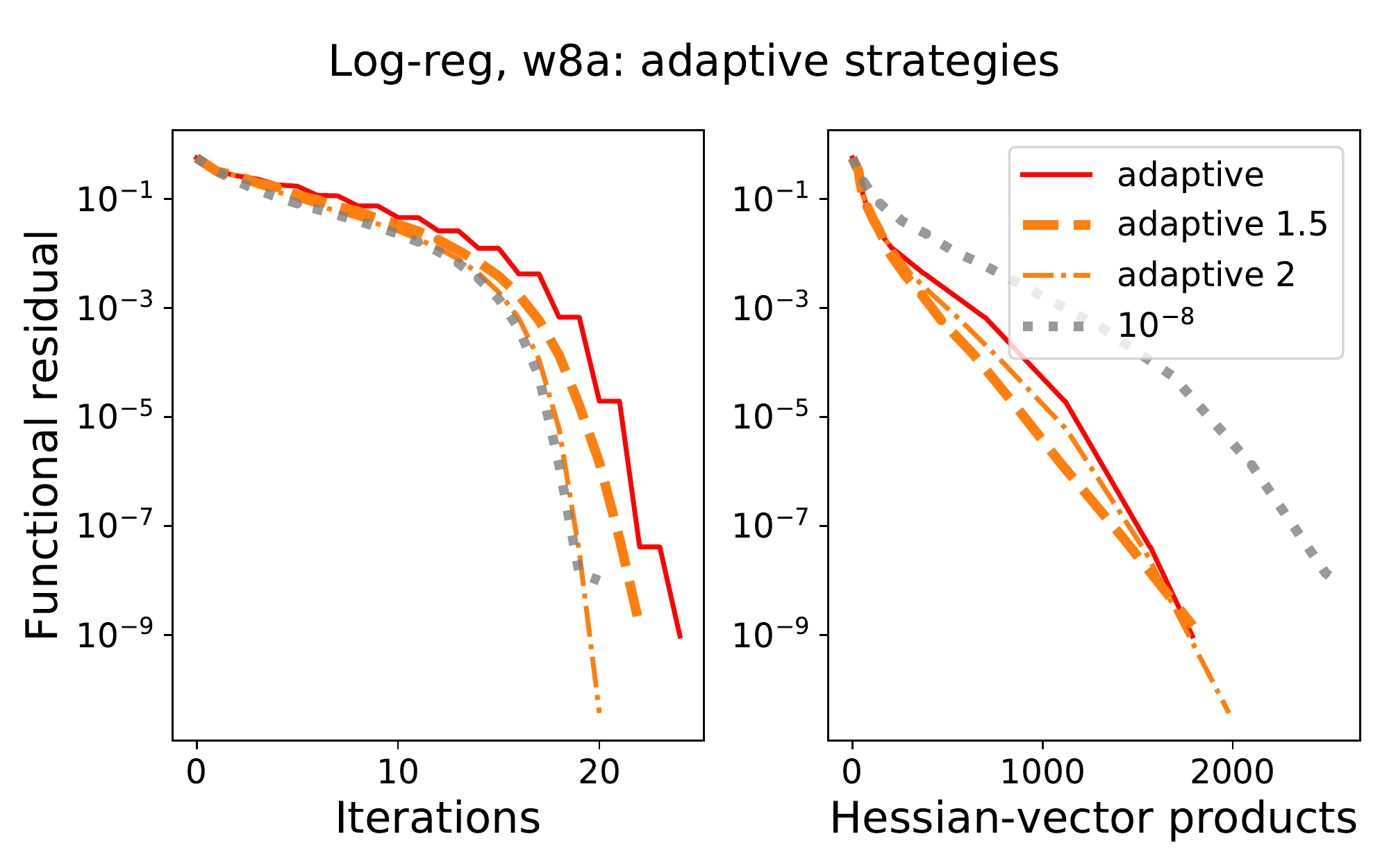}
	\end{minipage}

	\begin{minipage}{0.33\textwidth}
		\centering
		\includegraphics[width=\textwidth ]{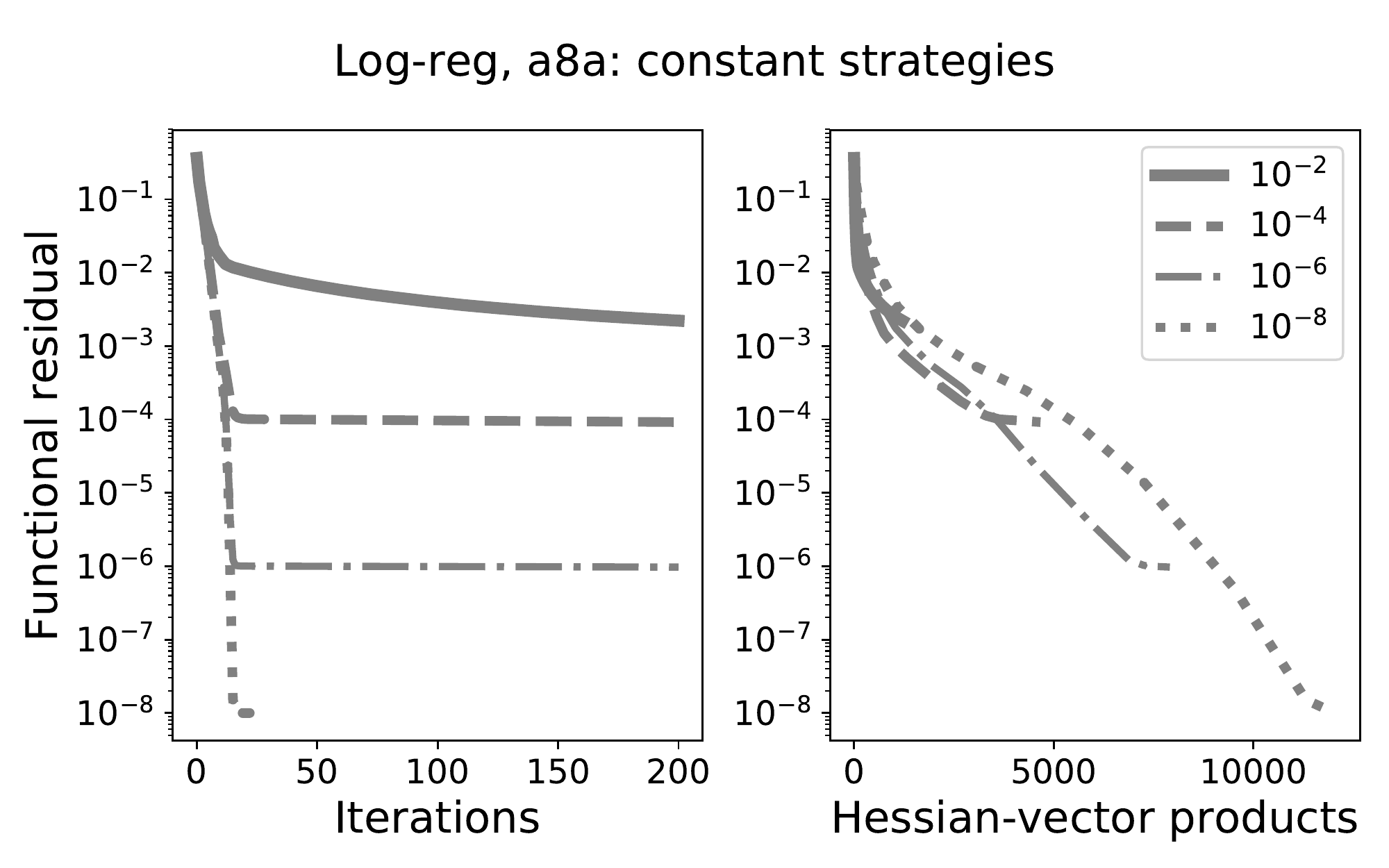}
	\end{minipage}
	\begin{minipage}{0.33\textwidth}
		\centering
		\includegraphics[width=\textwidth ]{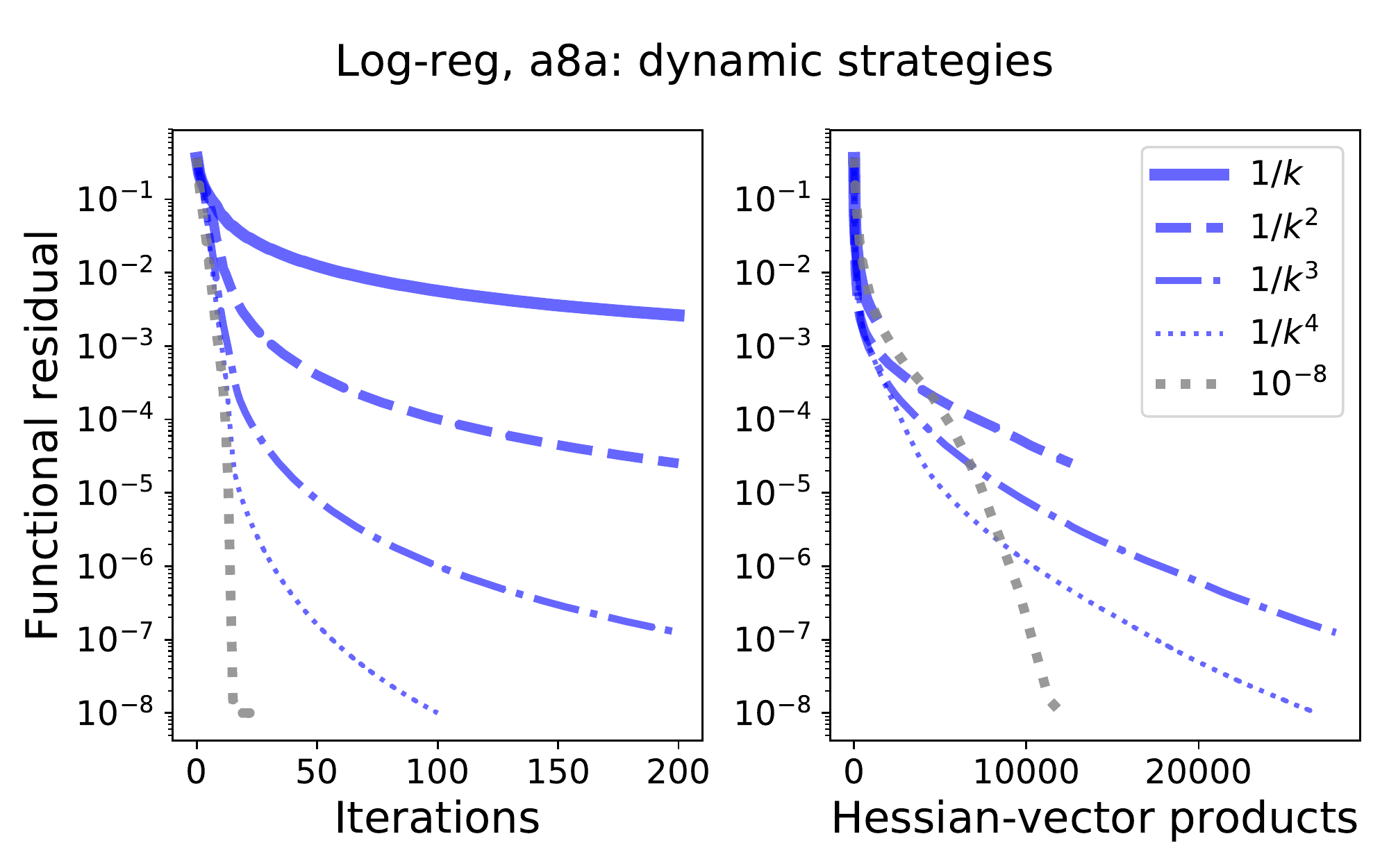}
	\end{minipage}
	\begin{minipage}{0.33\textwidth}
		\centering
		\includegraphics[width=\textwidth ]{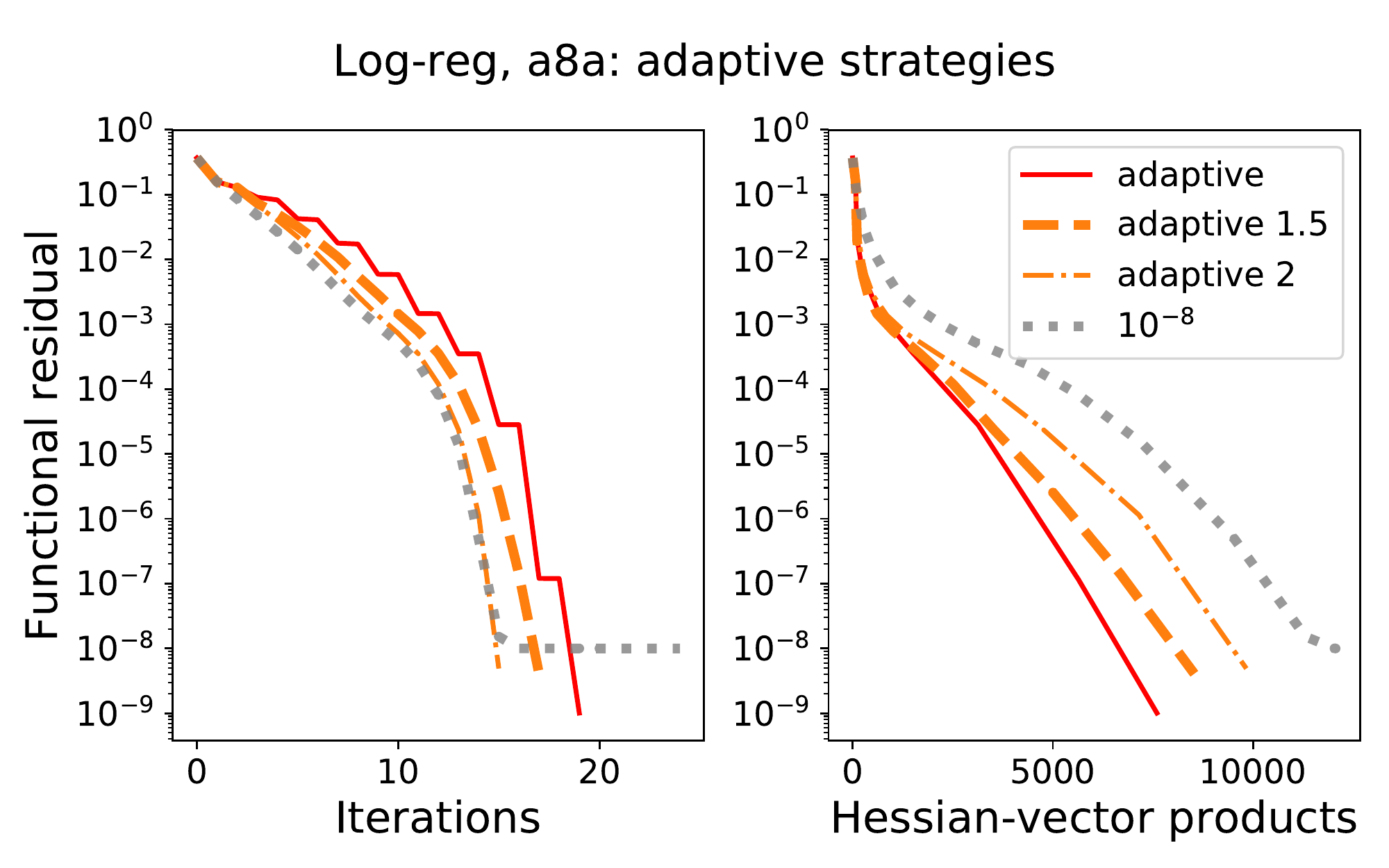}
	\end{minipage}
	
	\begin{minipage}{0.33\textwidth}
		\centering
		\includegraphics[width=\textwidth ]{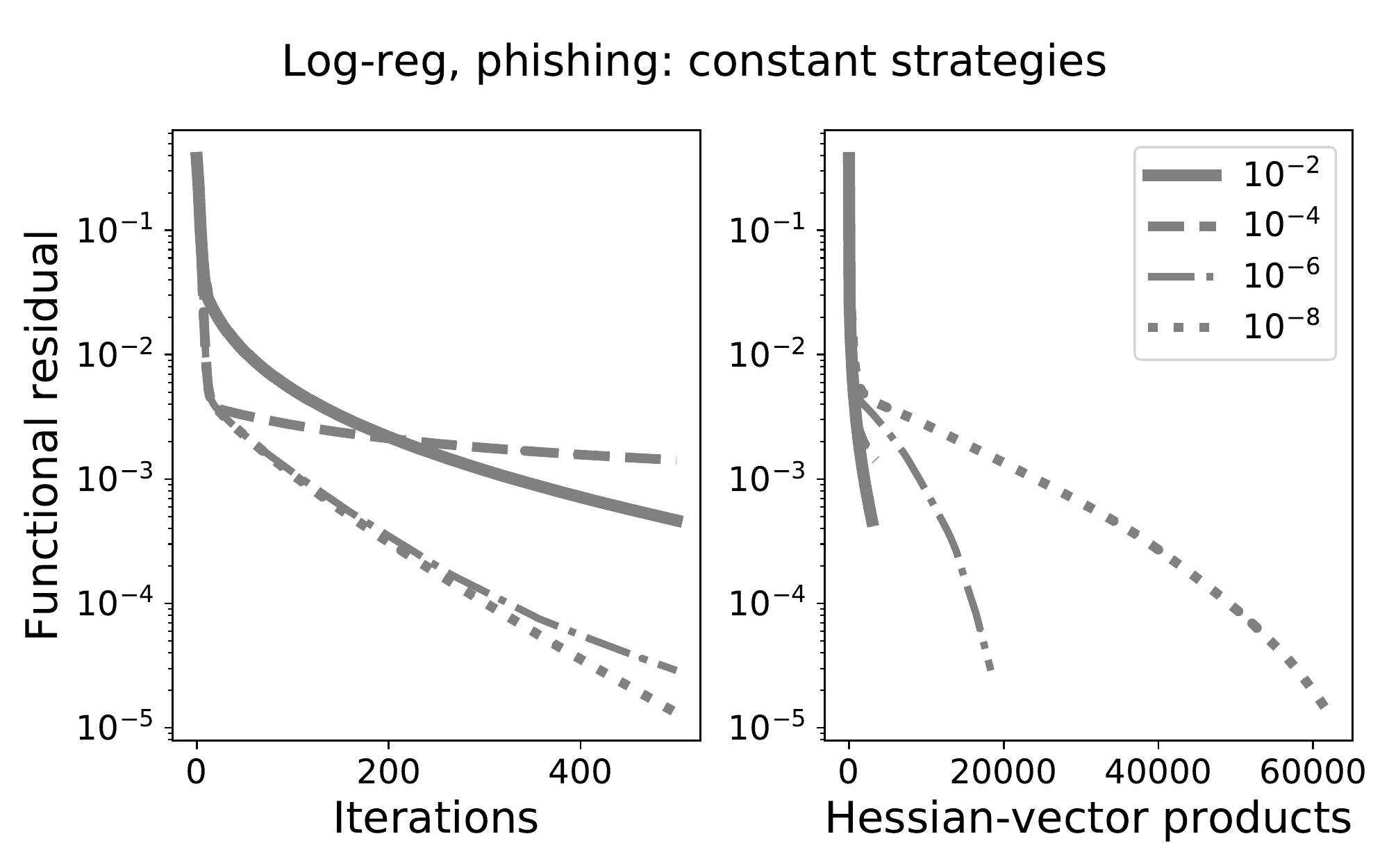}
	\end{minipage}
	\begin{minipage}{0.33\textwidth}
		\centering
		\includegraphics[width=\textwidth ]{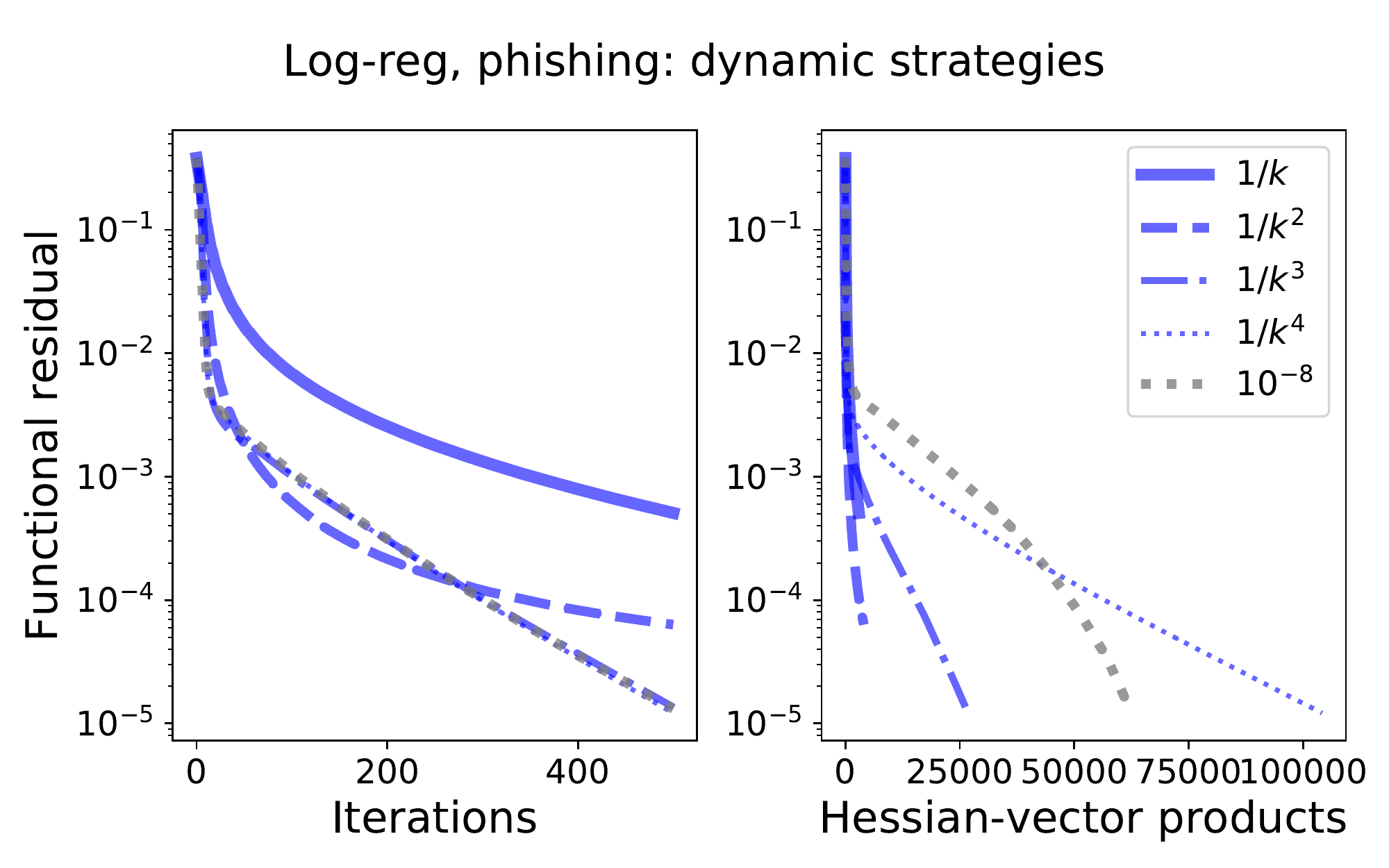}
	\end{minipage}
	\begin{minipage}{0.33\textwidth}
		\centering
		\includegraphics[width=\textwidth ]{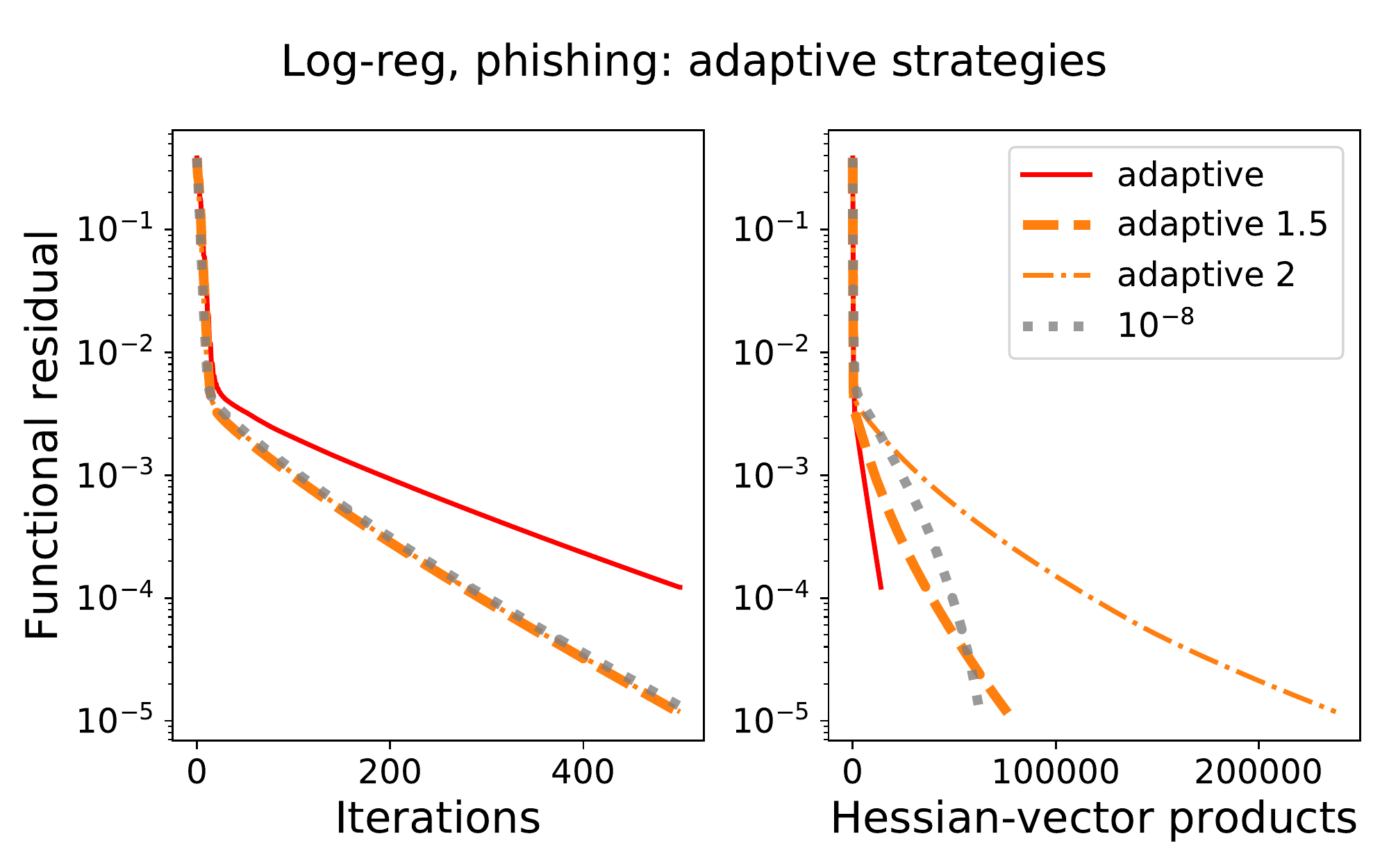}
	\end{minipage}

	\begin{minipage}{0.33\textwidth}
	\centering
	\includegraphics[width=\textwidth ]{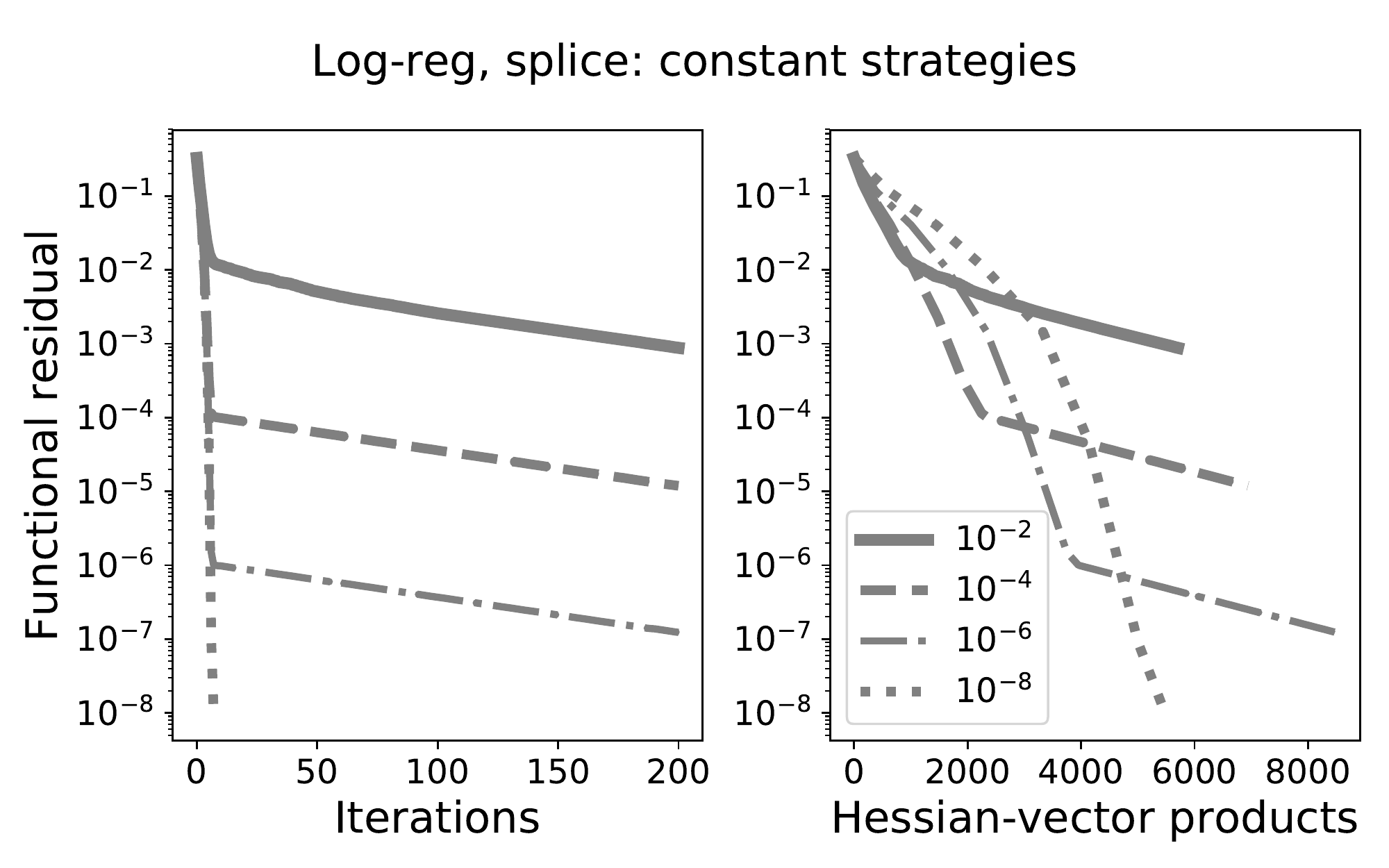}
	\end{minipage}
	\begin{minipage}{0.33\textwidth}
	\centering
	\includegraphics[width=\textwidth ]{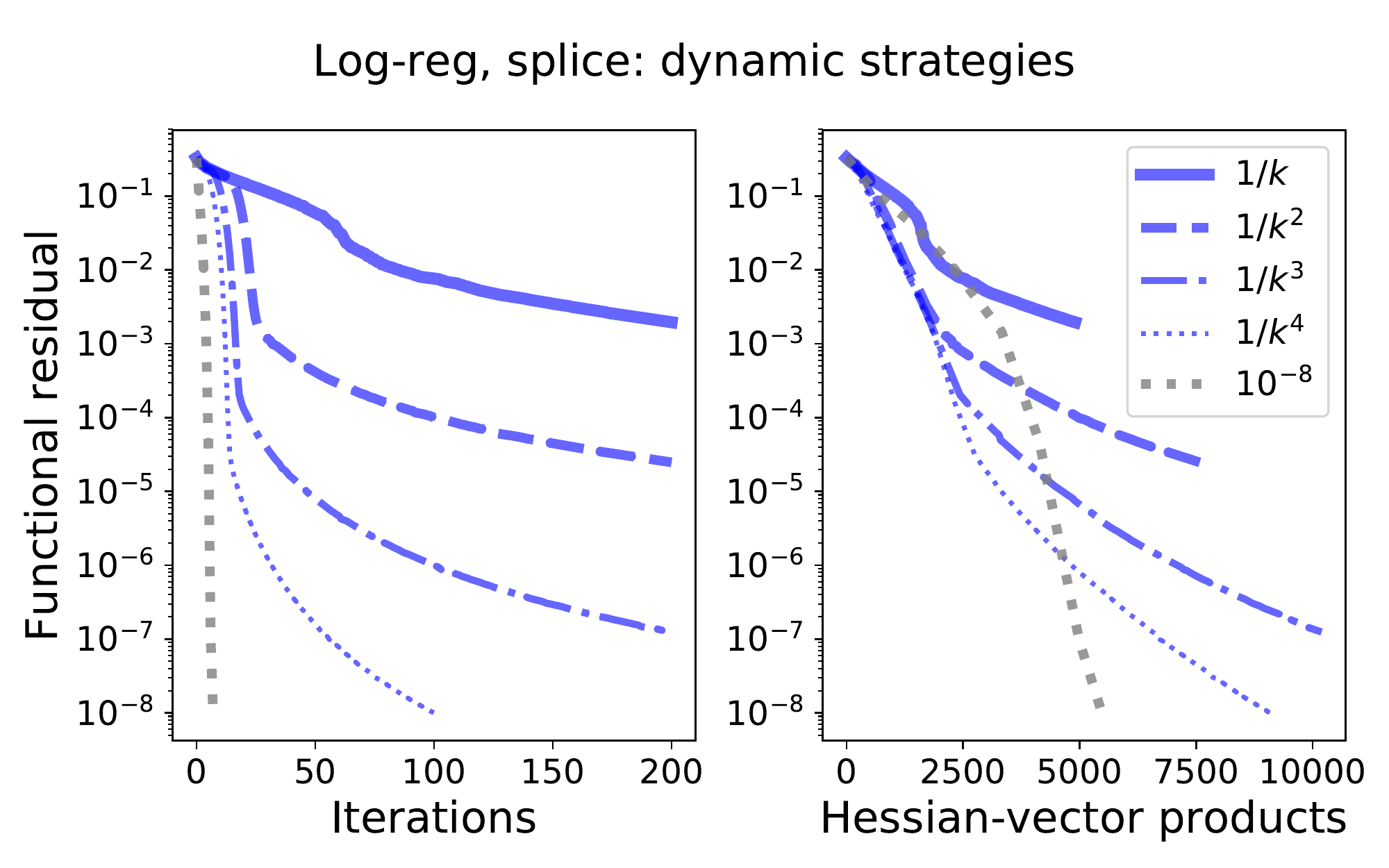}
	\end{minipage}
	\begin{minipage}{0.33\textwidth}
	\centering
	\includegraphics[width=\textwidth ]{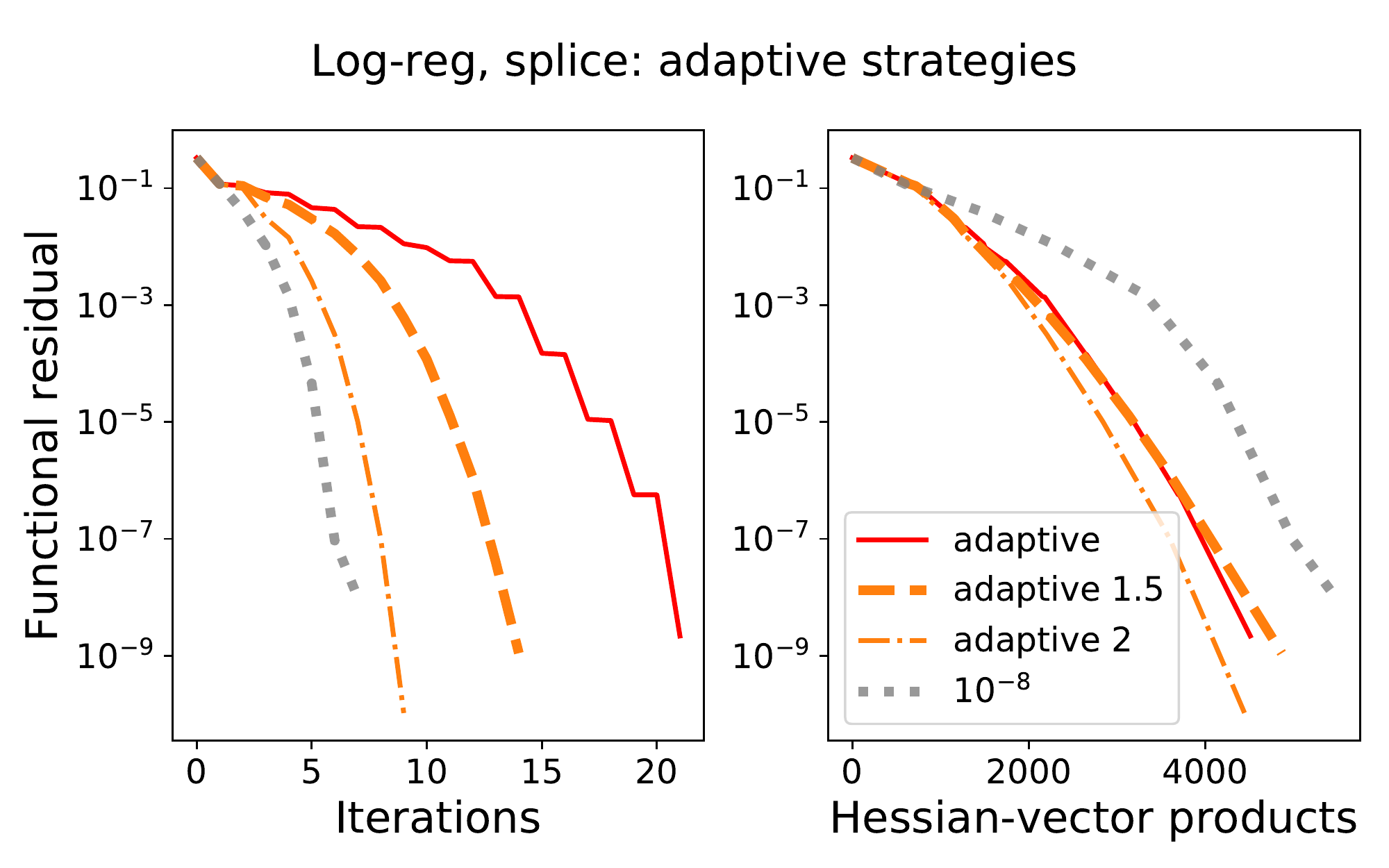}
	\end{minipage}

	\caption{Exact stopping criterion, training Logistic regression.} 
	\label{fig:CN_logreg_exact}
\end{figure}

In the following set of experiments with Cubic Newton method,
we compute the \textit{exact} minimizer of the model~\eqref{GenScheme}, at every iteration.
Then, we use this value to ensure the required precision in function value of the subproblem
for the inexact step (in the previous settings we used the upper bound~\eqref{UnifConvBound} for this purpose).
The results for Log-Sum-Exp function are shown on Figures~\ref{fig:CN_log_sum_exp_exact_100} and~\ref{fig:CN_log_sum_exp_exact_200}. The results for Logistic regression are shown on Figure~\ref{fig:CN_logreg_exact}.

We compare the iteration rate and 
the corresponding number of Hessian-vector products used,
for the constant choice of inner accuracy (left graphs),
dynamic strategies in the form $\delta_k = 1/k^{\alpha}$ (center),
and adaptive strategies $\delta_k = (F(x_{k - 1}) - F(x_k))^{\alpha}$ (right graphs).
We use the names "adaptive", "adaptive 1.5" and "adaptive 2" for
$\alpha = 1$, $\alpha = 3/2$, and $\alpha = 2$, respectively.

We see, that the constant choice of inner accuracy reasonably depends 
on the desired precision for solving the initial problem.
At the same time, the dynamic strategies are adjusting with the iterations.
The best performance is achieved by the use of the adaptive policies.
It is also notable, that in some cases, we need to use "adaptive 1.5" or "adaptive 2" strategy,
to have the local superlinear convergence. This confirms our theory.

\newpage

\subsection{Averaging and Acceleration}

In this experiment, we consider unconstrained minimization 
of the following objective ($x^{(i)}$ indicates $i$th coordinate of $x$)
\beq \label{PoweredDiff}
\ba{rcl}
f(x) & = & |x^{(1)}|^3 + \sum\limits_{i = 2}^n |x^{(i)} - x^{(i - 1)}|^3, 
\qquad x \in \R^n,
\ea
\eeq
by different inexact Newton methods.
Note, that the structure of~\eqref{PoweredDiff} is similar to that one
of the worst function for the second-order methods
(see Chapter~4.3.1 in~\cite{nesterov2018lectures}).
It is also similar to the function from Example~\ref{ex:BadFunction}.

We compare iteration rates of the following algorithms: 
Cubic Newton (CN) with dynamic rule $\delta_k = 1/k^3$,
Cubic Newton with adaptive rule~\eqref{DeltaDynamic},
the method with Averaging (Algorithm~\ref{alg:Averaging})
with $\delta_k = 1/k^3$,
and the accelerated method with Contracting proximal iterations 
(Algorithm~\ref{alg:Accelerated}).
For the latter one we use $\zeta_k = 1/k$ and $\delta_k = 1/k$,
as the rules for choosing the accuracy of inexact (outer) proximal steps,
and inexact (inner) Newton steps, respectively.\footnote{In our experiments, there is no need of high precision for the inexact contracting proximal steps. A faster decrease of $\delta_k$ did not improve the rate of convergence.}

\begin{figure}[h!]	
	\centering 
	\begin{minipage}{0.30\textwidth}
		\centering
		\includegraphics[width=\textwidth ]{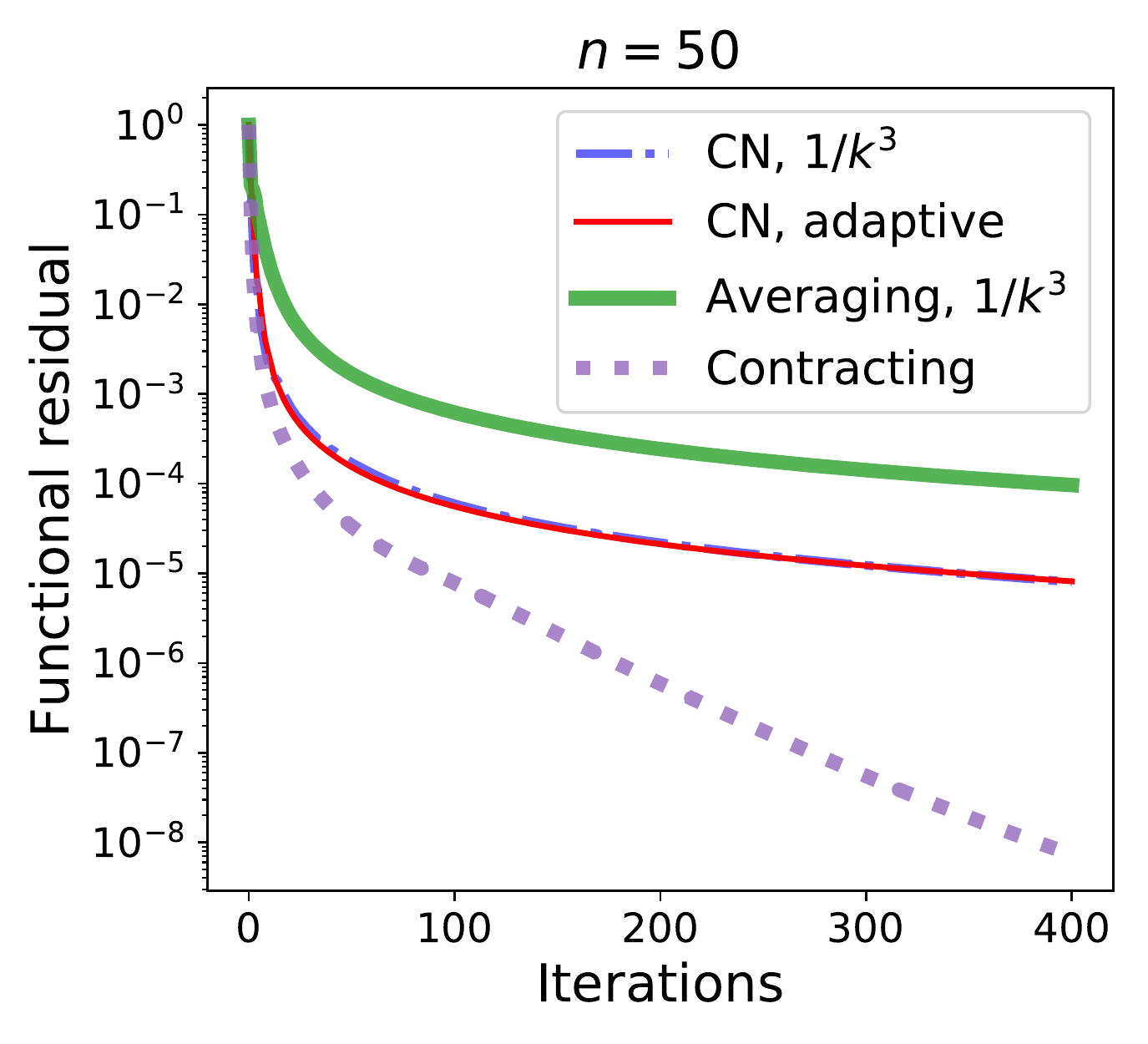}
	\end{minipage}
	\begin{minipage}{0.30\textwidth}
		\centering
		\includegraphics[width=\textwidth ]{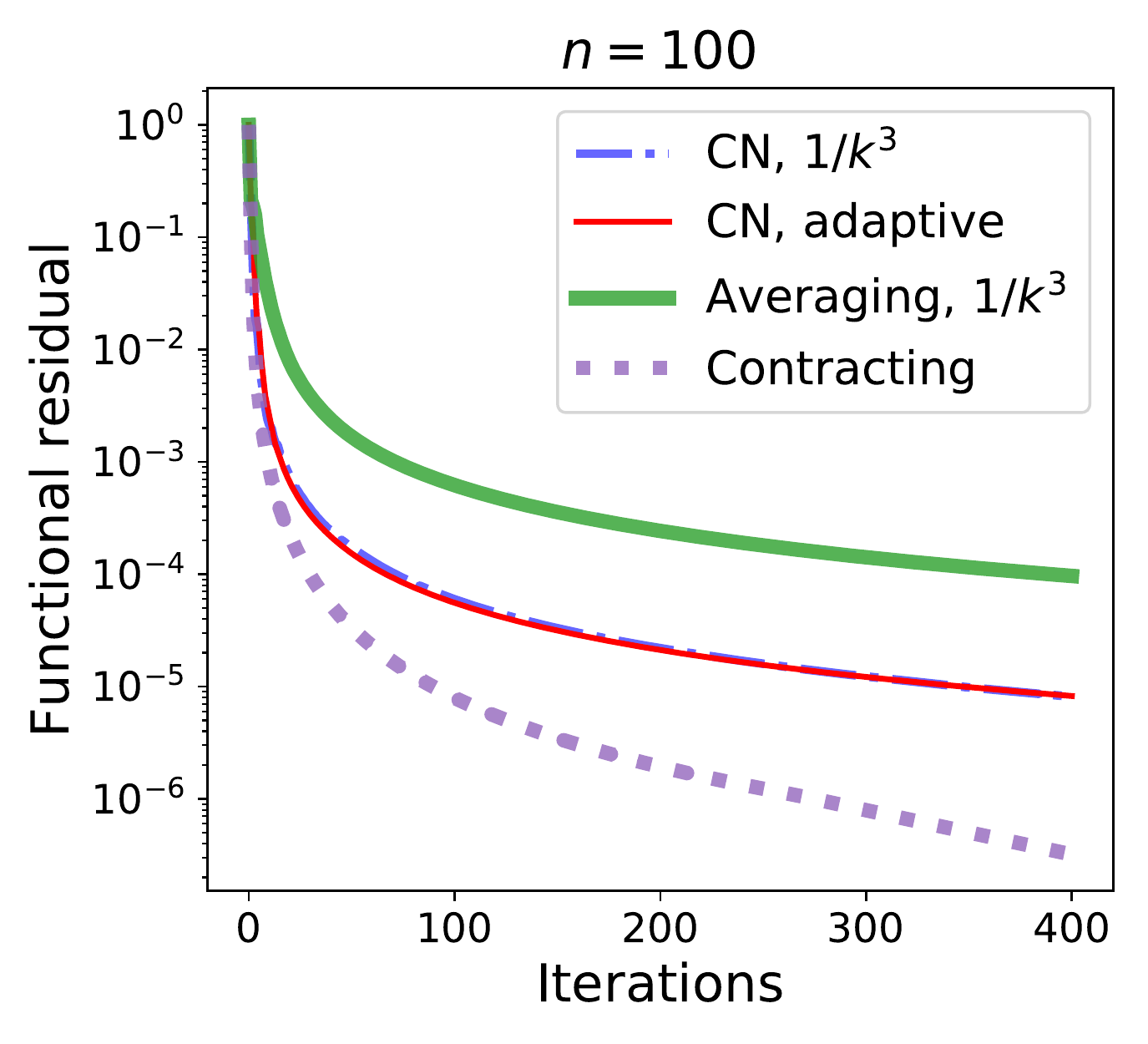}
	\end{minipage}
	\begin{minipage}{0.30\textwidth}
		\centering
		\includegraphics[width=\textwidth ]{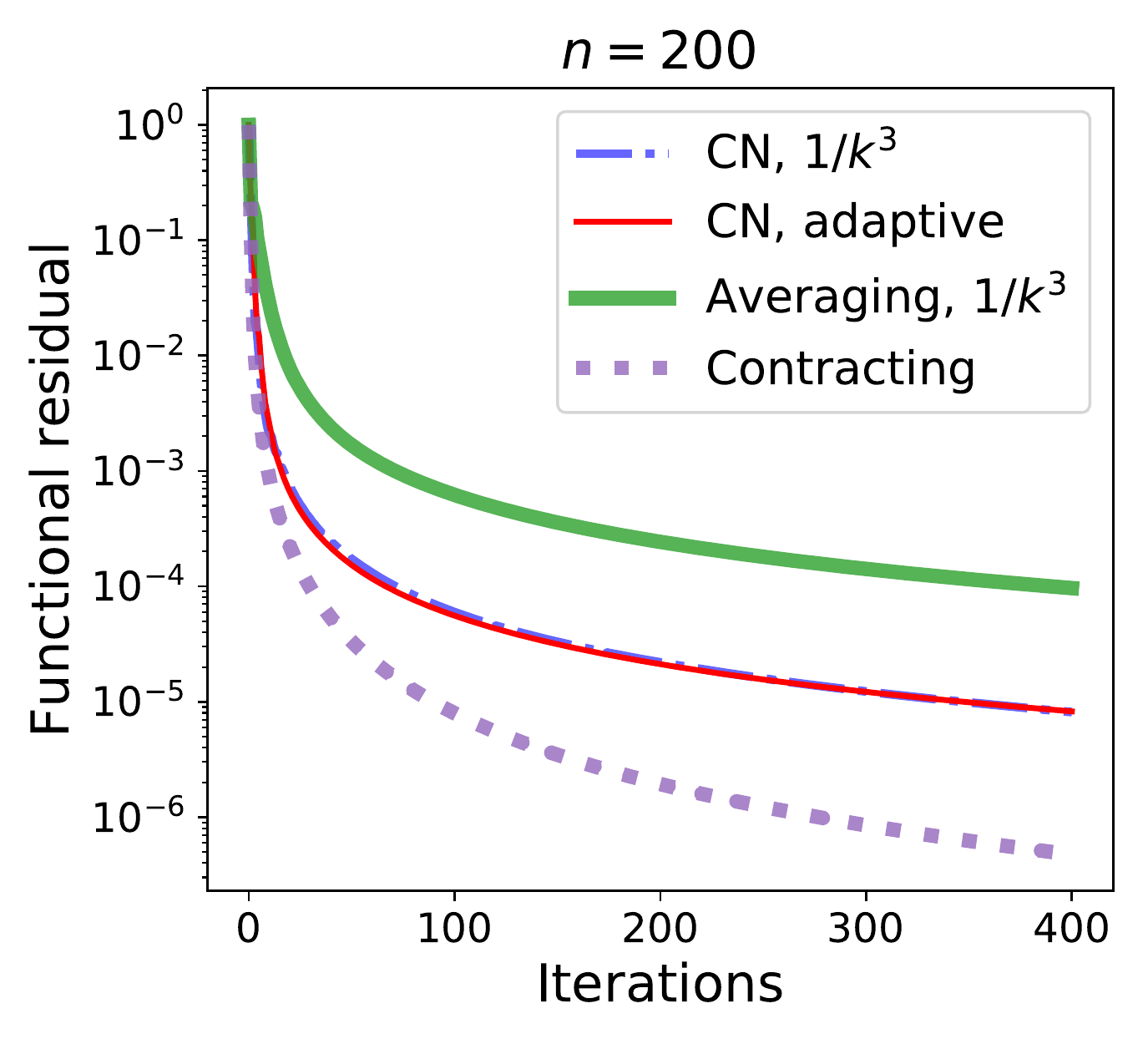}
	\end{minipage}

	\begin{minipage}{0.30\textwidth}
		\centering
		\includegraphics[width=\textwidth ]{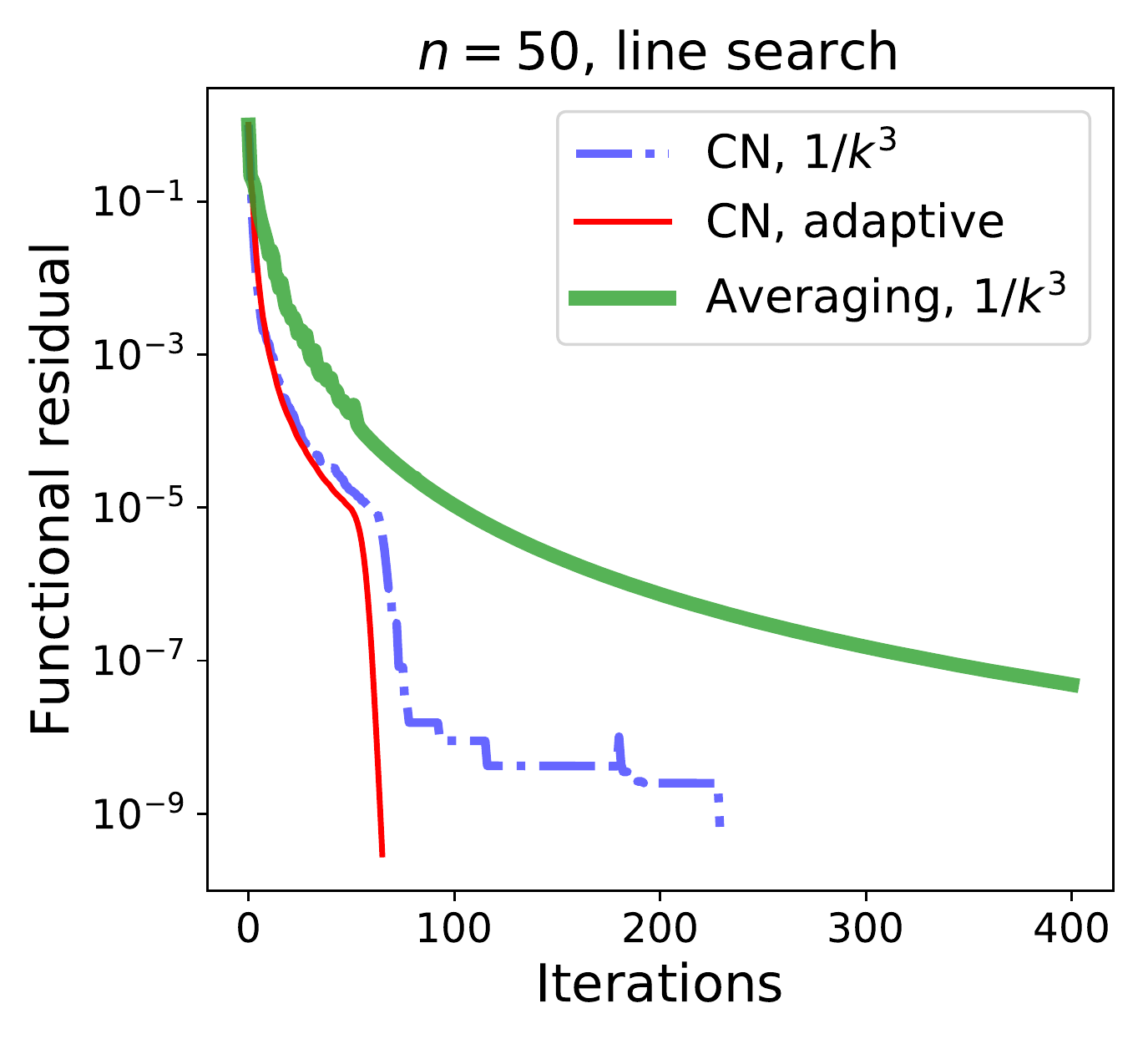}
	\end{minipage}
	\begin{minipage}{0.30\textwidth}
		\centering
		\includegraphics[width=\textwidth ]{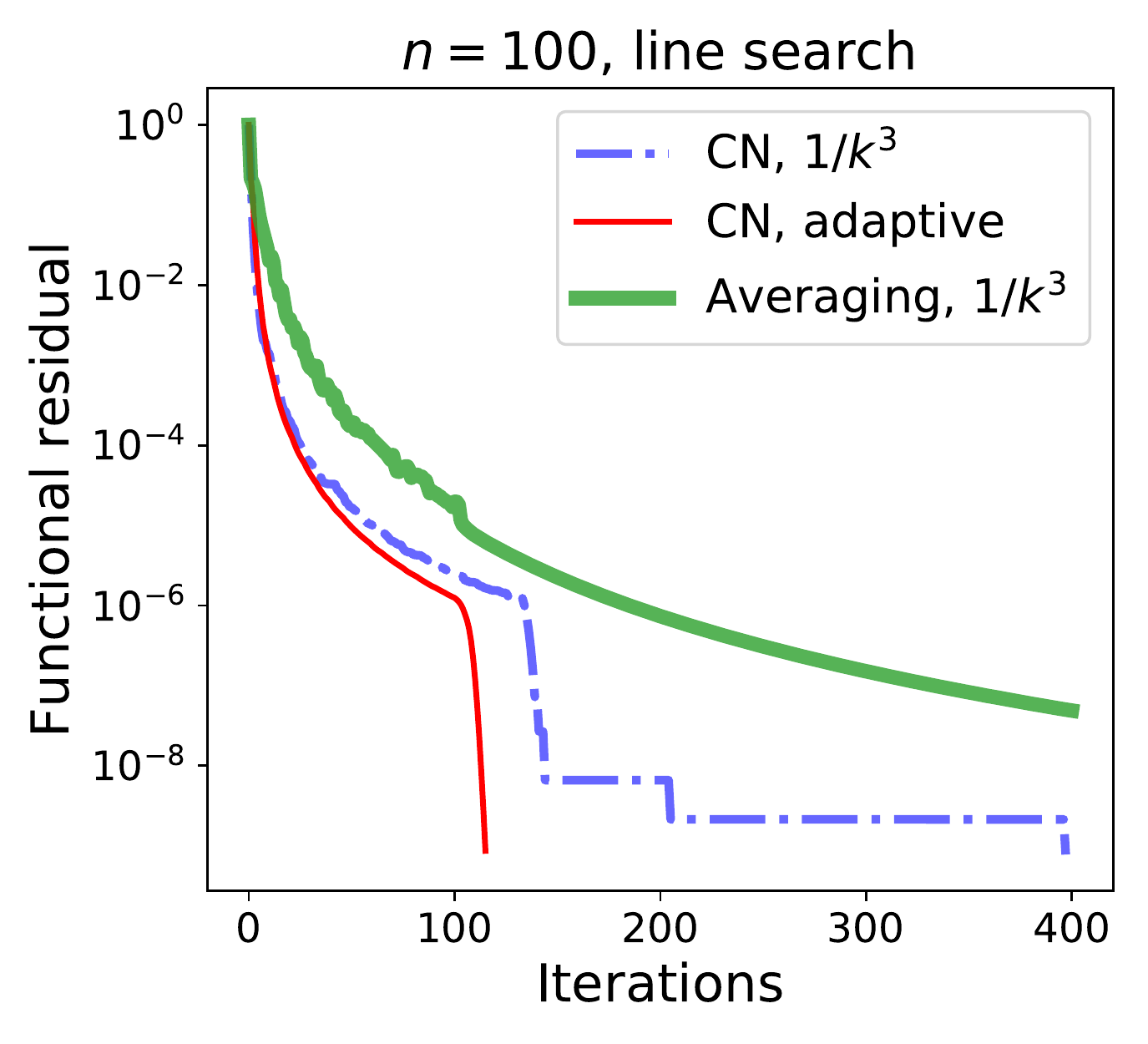}
	\end{minipage}
	\begin{minipage}{0.30\textwidth}
		\centering
		\includegraphics[width=\textwidth ]{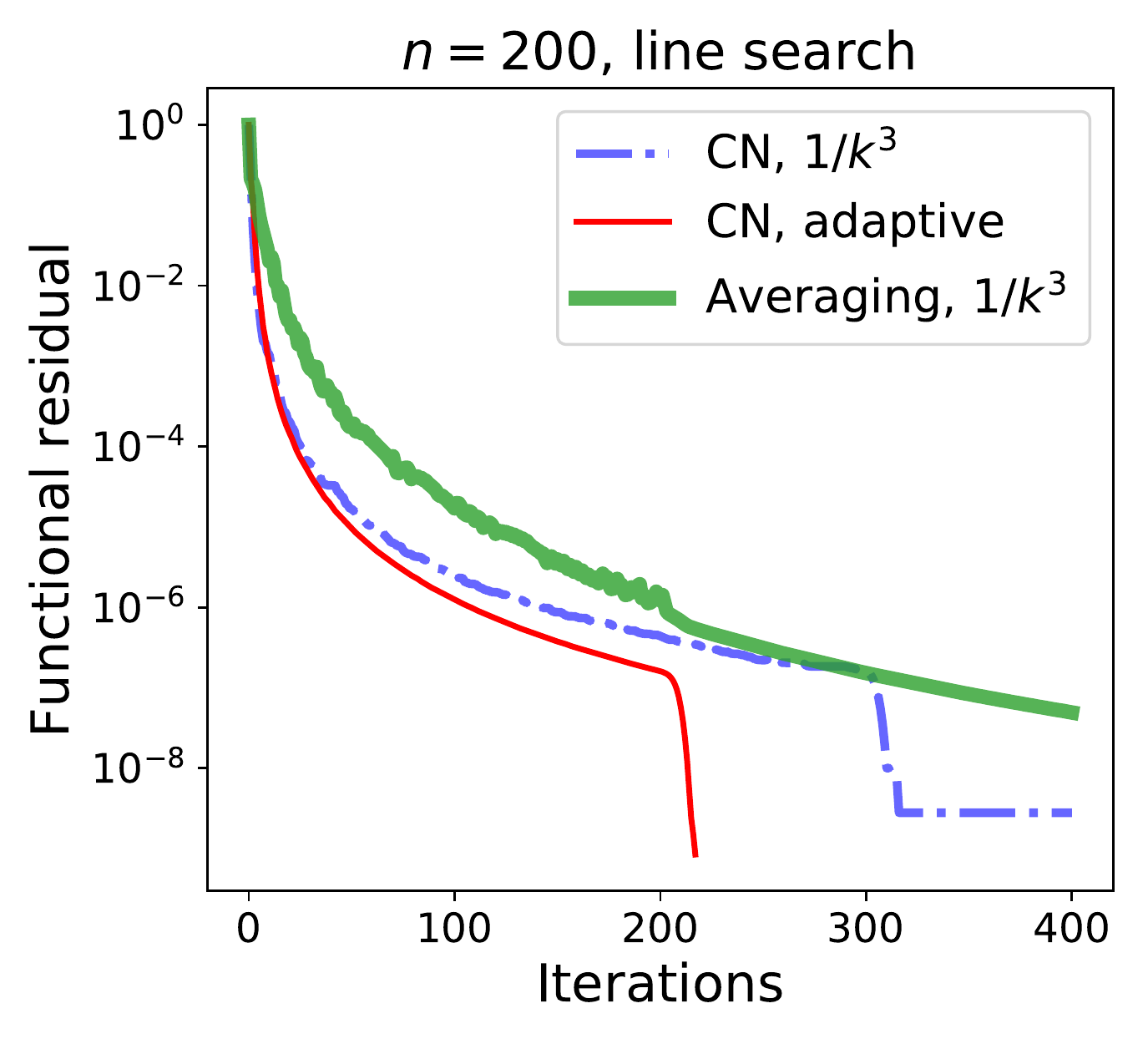}
	\end{minipage}

	\caption{ Methods with averaging and acceleration.} 
\label{fig:CN_averaging}
\end{figure}

For the first three algorithms, we also compare the constant choice for the regularization parameter: $H = 1$ (on the top graphs),
and a simple line search\footnote{Namely, we multiply $H$ by the factor of two, until condition $F(T_{H, \delta}(x^k)) \leq \Omega_H(x^k; T_{H, \delta}(x^k))$ is satisfied. At the next iteration, we start
	the line search from the previous estimate of $H$, divided by two.} 
for choosing $H$ at every iteration (bottom). The results are shown on Figure~\ref{fig:CN_averaging}.

We see, that all the methods have a \textit{sublinear} rate of convergence,
until the iteration counter is smaller than the dimension of the problem.
The use of the line search significantly helps for improving the rate.
Thus, it seems to be an important open question 
(which we keep for the further research) ---
to equip the contracting proximal scheme
(Algorithm~\ref{alg:Accelerated}) with a variant of line search as well.

\newpage

\section{Auxiliary Results}

\BL
For every $a, b \geq 0$ and integer $p \geq 1$ it holds
\beq \label{PoweredSumBound}
\ba{rcl}
(a + b)^{p - 1} & \leq & 2^{p - 2}a^{p - 1} + 2^{p - 2} b^{p - 1}.
\ea
\eeq 
\EL
\proof
For $p = 1$ it is trivial. For $p \geq 2$
we use convexity of $y(x) = x^{p - 1}$:
$$
\ba{rcl}
\bigl(\frac{a + b}{2}\bigr)^{p - 1}
& \leq & \frac{a^{p - 1}}{2} + \frac{b^{p - 1}}{2},
\ea
$$
which is~\eqref{PoweredSumBound}.
\qed

\BL
For every $s > 1$, it holds
\beq \label{SublinBound}
\ba{rcl}
\sum\limits_{i = 1}^k \frac{1}{i^s}
& \leq & 
\frac{s}{s - 1}.
\ea
\eeq
\EL
\proof
Indeed,
$$
\ba{rcl}
\sum\limits_{i = 1}^k \frac{1}{i^s} 
\;\; = \;\; 1 + \sum\limits_{i = 2}^k \frac{1}{i^s}
& \leq & 
1 + \int\limits_1^{+\infty} \frac{dx}{x^{s}}
\;\; = \;\; \frac{s}{s - 1}.
\ea
$$
\qed

\BL \label{LemmaLogSumExp}
For a given $a_i \in \E^{*}$, $1 \leq i \leq m$, consider
the log-sum-exp function
$$
\ba{rcl}
f(x) & = & \log\left( \sum\limits_{i = 1}^m e^{\la a_i, x \ra} \right), 
\quad x \in \E.
\ea
$$
Then, for Euclidean norm $\|x\| = \la Bx, x \ra^{1 / 2}, x \in \E$ with $B = \sum\limits_{i = 1}^m a_i a_i^{*}$
(assuming $B \succ 0$, otherwise we can reduce dimensionality of the problem),
we have the following estimates for the Lipschitz constants:
$$
\ba{rcccccl}
L_1 & = & 1, \qquad L_2 & = & 2, \qquad L_3 & = & 4. 
\ea
$$
\EL
\proof
Denote $\kappa(x) = \sum_{i = 1}^m e^{\la a_i, x \ra}$.
Then, for all $x \in \E$ and $h \in \E$, we have 
$$
\ba{rcl}
\la \nabla f(x), h \ra & = &
\frac{1}{\kappa(x)} \sum\limits_{i = 1}^m e^{\la a_i, x \ra} \la a_i, h \ra, \\
\\
\la \nabla^2 f(x) h, h \ra & = & 
\frac{1}{\kappa(x)} 
\sum\limits_{i = 1}^m e^{\la a_i, x \ra} ( \la a_i, h \ra - \la \nabla f(x), h \ra )^2
\;\; \leq \;\; 
\sum\limits_{i = 1}^m \la a_i, h \ra^2
\;\; = \;\; \|h\|^2,
\\
\\
D^3 f(x)[h]^3
& = & 
\frac{1}{\kappa(x)}
\sum\limits_{i = 1}^m e^{\la a_i, x \ra} ( \la a_i, h \ra - \la \nabla f(x), h \ra )^3
\;\; \leq \;\;
\la \nabla^2 f(x)h, h\ra \max\limits_{1 \leq i, j \leq m} \la a_i - a_j, h \ra 
\;\; \leq \;\; 
2 \|h\|^3, \\
\\
D^4 f(x)[h]^4
& = & 
\frac{1}{\kappa(x)}
\sum\limits_{i = 1}^m e^{\la a_i, x \ra} ( \la a_i, h \ra - \la \nabla f(x), h \ra )^4
- 3 \la \nabla^2 f(x)h, h \ra^2 \\
\\
& \leq &
D^3 f(x)[h]^3 \max\limits_{1 \leq i, j \leq m} \la a_i - a_j, h \ra 
\;\; \leq \;\;
4 \|h\|^4.
\ea
$$

\qed

\BL \label{LemmaStrongUniformConvex}
Let $\psi(x) = \frac{\mu}{2}\|x - x_0\|^2$, $\mu \geq 0$.
Consider the ball of radius $D$ around some fixed point $c$:
$$
\ba{rcl}
\mathcal{B} & = & \{ x : \|x - c\| \leq D \}.
\ea
$$
Then, for any $p \geq 1$, it holds, 
$$
\ba{rcl}
\psi(y) - \psi(x) - \la \nabla \psi(x), y - x \ra 
& \geq & \frac{\sigma_{p + 1}}{p + 1}\|y - x\|^{p + 1},
\qquad x, y \in \mathcal{B},
\ea
$$
with $\sigma_{p + 1} = \frac{(p + 1)\mu}{2^p D^{p - 1}}$.
Thus function $\psi$ is uniformly convex of degree $p + 1$ on a ball.
\EL
\proof
For all $x, y \in \mathcal{B}$, we have 
$$
\ba{rcl}
\|x - y\| & \leq & \|x - c\| + \|c - y\| 
\;\; \leq \;\; 2D.
\ea
$$
Therefore, 
$$
\ba{rcl}
\psi(y) - \psi(x) - \la \nabla \psi(x), y - x \ra
& = & \frac{\mu}{2}\|y - x\|^2 \\
\\
& = & \frac{\mu}{2\|y - x\|^{p - 1}}\|y - x\|^{p + 1} \\
\\
& \geq & 
\frac{\mu}{2^p D^{p - 1}} \|y - x\|^{p + 1}
\;\; = \;\;
\frac{\sigma_{p + 1}}{p + 1}\|y - x\|^{p + 1}.
\ea
$$

\qed

\section{Proof of Theorem~\ref{th:GlobalConvAlg1}}
\proof
Indeed, by Lemma~\ref{lemma:Global}, 
for every $y \in \dom \psi$ we have
\beq \label{eq-th1-gbound}
\ba{rcl}
F(x_{k + 1}) & \leq & F(T_{k + 1})
\;\; \refLE{OneStep} \;\;
F(y) + \frac{L_p\|y - x_k\|^{p + 1}}{p!} + \delta_{k + 1},
\quad k \geq 0.
\ea
\eeq
Let us introduce an \textit{arbitrary} sequence of positive increasing coefficients
$\{ A_k \}_{k \geq 0}$, $A_0 \Def 0$. Denote $a_{k + 1} \Def A_{k + 1} - A_k$.
Then, plugging $y = \frac{a_{k + 1}x^{*} + A_k x_k}{A_{k + 1}}$ into~\eqref{eq-th1-gbound},
we obtain
$$
\ba{rcl}
F(x_{k + 1}) & \leq & \frac{a_{k + 1}}{A_{k + 1}}F^{*} + \frac{A_k}{A_{k + 1}}F(x_k)
+ \frac{a_{k + 1}^{p + 1}}{A_{k + 1}^{p + 1}} \frac{L_p \|x_k - x^{*}\|^{p + 1}}{p!}
+ \delta_{k + 1},
\ea
$$
or, equivalently
$$
\ba{rcl}
A_{k + 1}(F(x_{k + 1}) - F^{*})
& \leq & 
A_k (F(x_k) - F^{*})
+ \frac{a_{k + 1}^{p + 1}}{A_{k + 1}^{p}} \frac{L_p \|x_k - x^{*}\|^{p + 1}}{p!}
+ A_{k + 1} \delta_{k + 1}.
\ea
$$
Summing up these inequalities, we get, for every $k \geq 1$
\beq \label{eq-th1-telescoped}
\ba{rcl}
A_k(F(x_k) - F^{*}) & \leq & \sum\limits_{i = 1}^k A_i \delta_i 
+ \frac{L_p}{p!} 
\sum\limits_{i = 1}^k \frac{a_i^{p + 1}}{A_i^p}\|x_i - x^{*}\|^{p + 1} \\
\\
& \leq &
\sum\limits_{i = 1}^k A_i \delta_i + 
\frac{L_p D^{p + 1}}{p!} \sum\limits_{i = 1}^k \frac{a_i^{p + 1}}{A_i^p},
\ea
\eeq
where the last inequality holds due to monotonicity of the method:
$$
\ba{rcl}
F(x_i) & \leq & F(x_0), \qquad i \geq 0.
\ea
$$
Finally, let us fix $A_{k} \equiv k^{p + 1}$. Then, by the mean value theorem,
for some $\xi \in [k - 1; k]$
$$
\ba{rcl}
a_k & = & A_k - A_{k - 1} \; \; = \; \; k^{p + 1} - (k - 1)^{p + 1} \\
\\
& = & (p + 1)\xi^{p} \; \; \leq \; \; (p + 1) k^p.
\ea
$$
Therefore, 
\beq \label{th1-ratio-bound}
\ba{rcl}
\sum\limits_{i = 1}^k \frac{a_i^{p + 1}}{A_i^p}
& \leq & 
\sum\limits_{i = 1}^k \frac{(p + 1)^{p + 1} i^{p(p + 1)}}{i^{(p + 1)p}}
\;\; = \;\; (p + 1)^{p + 1} k,
\ea
\eeq
and
\beq \label{th1-delta-sum}
\ba{rcl}
\sum\limits_{i = 1}^k A_i \delta_i
& = & 
\sum\limits_{i = 1}^k \frac{c i^{p + 1}}{i^{p + 1}}
\;\; = \;\; ck.
\ea
\eeq
Plugging these bounds into~\eqref{eq-th1-telescoped} completes the proof.
\qed

\section{Proof of Theorem~\ref{th:GlobalConvAlg2}}
\proof
First, by the same reasoning as in Theorem~\ref{th:GlobalConvAlg1}, 
we obtain the following bound, for every $k \geq 1$:
\beq \label{eq-th2-telescoped}
\ba{rcl}
A_k (F(x_k) - F^{*})
& \leq & \sum\limits_{i = 1}^k A_i \delta_i + \frac{L_p D^{p + 1}}{p!} 
\sum\limits_{i = 1}^k \frac{a_i^{p + 1}}{A_i^p},
\ea
\eeq
where $\{ A_k \}_{k \geq 0}$ is an \textit{arbitrary} sequence of increasing coefficients,
with $A_0 \Def 0$, and $a_{k} \Def A_{k} - A_{k - 1}$.

Substituting into~\eqref{eq-th2-telescoped} the values 
$\delta_{i} = c(F(x_{i - 2}) - F(x_{i - 1})), \; i \geq 2$, we have
\beq \label{eq-th2-main}
\ba{rcl}
A_k (F(x_k) - F^{*})
& \leq & A_1 \delta_1 + c \sum\limits_{i = 2}^k A_i (F(x_{i - 2}) - F(x_{i - 1})) + \frac{L_p D^{p + 1}}{p!} 
\sum\limits_{i = 1}^k \frac{a_i^{p + 1}}{A_i^p}, \quad k \geq 1,
\ea
\eeq
or, rearranging the terms, it holds for every $k \geq 2$:
\beq \label{eq-th2-main-2}
\ba{cl}
& (c + 1) A_k (F(x_k) - F^{*}) \\
\\
&\quad  \;\;\,\, \leq \;\;\,\, A_k ( F(x_k) - F^{*} ) + cA_k (F(x_{k - 1}) - F^{*}) \\
\\
&\quad \;\; \refLE{eq-th2-main} \;\;
\frac{L_p D^{p + 1}}{p!} \sum\limits_{i = 1}^k \frac{a_i^{p + 1}}{A_i^{p}}
+ c \sum\limits_{i = 1}^{k - 2} (A_{i + 2} - A_{i + 1})(F(x_i) - F^{*})
+ A_1 \delta_1 + cA_2(F(x_0) - F^{*}) \\
\\
&\quad \;\;\,\, = \;\;\,\,
\frac{L_p D^{p + 1}}{p!} \sum\limits_{i = 1}^k \frac{a_i^{p + 1}}{A_i^{p}}
+ c \sum\limits_{i = 1}^{k - 2} a_{i + 2} (F(x_i) - F^{*})
+ A_1 \delta_1 + cA_2(F(x_0) - F^{*}),
\ea
\eeq
and for $k = 1$ we have
\beq \label{eq-th2-base}
\ba{rcl}
A_1(F(x_1) - F^{*}) & \refLE{eq-th2-main} & 
A_1 \delta_1 + \frac{L_p D^{p + 1}}{p!} \frac{a_1^{p + 1}}{A_1^p}
\;\; = \;\;
A_1 \Bigl( \delta_1 + \frac{L_p D^{p + 1}}{p!} \Bigr).
\ea
\eeq
Now, let us pick $A_k \equiv k^{p + 2}$. Then,
$$
\ba{rcl}
a_k & \equiv & k^{p + 2} - (k - 1)^{p + 2} 
\;\; \leq \;\;
(p + 2) k^{p + 1},
\ea
$$
and
$$
\ba{rcl}
\sum\limits_{i = 1}^k \frac{a_i^{p + 1}}{A_i^p}
& \leq &
(p + 2)^{p + 1} \sum\limits_{i = 1}^k \frac{i^{(p + 1)^2}}{i^{(p + 2)p}}
\;\; = \;\;
(p + 2)^{p + 1} \sum\limits_{i = 1}^k i
\;\; \leq \;\;
(p + 2)^{p + 1} k^2.
\ea
$$
Therefore, \eqref{eq-th2-main-2} leads to
\beq \label{eq-th2-main-3}
\ba{cl}
& (c + 1) k^{p + 2} (F(x_k) - F^{*}) \\
\\
& \quad \;\;\,\, \leq \;\;\,\,
\frac{(p + 2)^{p + 1} k^2 L_p D^{p + 1}}{p!}
+ c(p + 2) \sum\limits_{i = 1}^{k - 2}
(i + 2)^{p + 1} ( F(x_i) - F^{*} )
+ \delta_1 + c 2^{p + 2}(F(x_0) - F^{*}),
\quad k \geq 2.
\ea
\eeq
And the statement to be proved is
\beq \label{eq-th2-prove}
\ba{rcl}
F(x_k) - F^{*}
& \leq & 
\frac{\beta}{k^{p + 2}}
+ 
\frac{\gamma L_p D^{p + 1}}{p! \, k^p}, \quad k \geq 1,
\ea
\eeq
where
\beq \label{eq-th2-beta-gamma}
\ba{rcl}
\beta & \Def & 
\frac{\delta_1 + c 2^{p + 2}(F(x_0) - F^{*})}{1 - c( (p + 2)^2 / (p + 1) - 1)}, 
\qquad
\gamma \;\; \Def \;\;
\frac{(p + 2)^{p + 1}}{1 - c( (p + 2)3^{p + 1} - 1 )}.
\ea
\eeq
Note, that from our assumptions, $c$ is small enough: $c \leq \frac{1}{(p + 2)3^{p + 1} - 1}$,
and~\eqref{eq-th2-beta-gamma} are correctly defined.

Let us prove~\eqref{eq-th2-prove} by induction. It holds for $k = 1$ by~\eqref{eq-th2-base}.
Assuming that it holds for all $1 \leq k \leq K - 2$, we have
$$
\ba{rcl}
F(x_K) - F^{*} 
& \overset{\eqref{eq-th2-main-3},\eqref{eq-th2-prove}}{\leq}  &
\frac{(p + 2)^{p + 1} L_p D^{p + 1}}{(c + 1) \, p! \, K^p }
+ \frac{c(p + 2)}{(c + 1)K^{p + 2}}
\sum\limits_{i = 1}^{K - 2}(i + 2)^{p + 1}
\Bigl( 
\frac{\gamma L_p D^{p + 1}}{p! \, i^p}
+ \frac{\beta}{i^{p + 2}} 
\Bigr)
+ \frac{\delta_1 + c2^{p + 2}(F(x_0) - F^{*})}{(c + 1) K^{p + 2}} \\
\\
& = &
\Bigl( 
\frac{(p + 2)^{p + 1}}{c + 1}
+ \frac{\gamma c(p + 2)}{(c + 1)K^2} \sum\limits_{i = 1}^{K - 2} 
\frac{(i + 2)^{p + 1}}{i^p}
\Bigr) \cdot \frac{L_p D^{p + 1}}{p! \, K^p} \\
\\
& \; &
\quad + \quad 
\Bigl(
\frac{\beta c(p + 2)}{(c + 1)} \sum\limits_{i = 1}^{K - 2} \frac{1}{i^{p + 2}}
+
\frac{\delta_1 + c2^{p + 2}(F(x_0) - F^{*})}{c + 1}
\Bigr) \cdot \frac{1}{K^{p + 2}} \\
\\
& \leq &
\frac{(p + 2)^{p + 1} + \gamma c (p + 2) 3^{p + 1}}{c + 1} 
 \cdot \frac{L_p D^{p + 1}}{p! \, K^p}
 +
\Bigl(
\frac{\beta c (p + 2)^2}{(c + 1)(p + 1)}
+ \frac{\delta_1 + c2^{p + 2} (F(x_0) - F^{*}) }{c + 1}
\Bigr)
\frac{1}{K^{p + 2}},
\ea
$$
where we have used in the last inequality the following simple bounds:
$$
\ba{rcl}
\sum\limits_{i = 1}^{K - 2} \frac{(i + 2)^{p + 1}}{i^p}
& \leq &
3^{p + 1} \sum\limits_{i = 1}^{K - 2} \frac{i^{p + 1}}{i^p} 
\;\; \leq \;\; 3^{p + 1} K^2, \\
\\
\sum\limits_{i = 1}^{K - 2} \frac{1}{i^{p + 2}}
& \refLE{SublinBound} & \frac{p + 2}{p + 1}.
\ea
$$
Therefore, to finish the proof, its enough to verify two equations:
$$
\ba{rcl}
\frac{\beta c(p + 2)^2}{(c + 1)(p + 1)} + \frac{\delta_1 + c2^{p + 2}(F(x_0) - F^{*})}{c + 1}
& = & \beta,
\quad
\text{and}
\quad
\frac{(p + 2)^{p + 1} + \gamma c(p + 2)3^{p + 1}}{c + 1} \;\; = \;\; \gamma.
\ea
$$
which are~\eqref{eq-th2-beta-gamma}.
\qed

\section{Proof of Theorem~\ref{th:GlobalStrongAlg2}}
\proof
Indeed, by Lemma~\ref{lemma:Global}, for every $y \in \dom \psi$ we have
\beq \label{th3:OneStep}
\ba{rcl}
F(x_{k + 1}) & \leq & F(y) + \frac{L_p \|y - x_k\|^{p + 1}}{p!} + \delta_{k + 1},
\quad k \geq 0.
\ea
\eeq
Let us substitute $y = \lambda x^{*} + (1 - \lambda ) x_k$
into~\eqref{th3:OneStep},
where 
$\lambda \equiv \omega_p^{-1/p} \in (0, 1]$.
This gives
$$
\ba{rcl}
F(x_{k + 1}) & \leq &
\lambda F^{*} + (1 - \lambda) F(x_k)
+ \frac{\lambda^{p + 1}  L_p\|x_k - x^{*}\|^{p + 1} }{p!} + \delta_{k + 1} \\
\\
& \leq &
\lambda F^{*} + (1 - \lambda) F(x_k) + 
\frac{\lambda^{p + 1} (p + 1)L_p}{\sigma_{p + 1} p!}(F(x_k) - F^{*}) + \delta_{k + 1},
\ea
$$
where we use uniform convexity in the last inequality. Therefore,
for every $k \geq 1$:
$$
\ba{rcl}
F(x_{k + 1}) - F^{*} & \leq & 
\Bigl(  1 - \omega_p^{-1/p} + \frac{\omega_p^{-1/p}}{p + 1}  \Bigr) (F(x_k) - F^{*})
+ \delta_{k + 1} \\
\\
& = &
\Bigl(1 - \frac{p}{p + 1} \omega_p^{-1/p}) (F(x_k) - F^{*})
+ c(F(x_{k - 1}) - F(x_k)) \\
\\
& \leq &
\Bigl(1 - \frac{p}{p + 1} \omega_p^{-1/p} + c \Bigr) (F(x_{k - 1}) - F^{*}),
\ea
$$
the last inequality uses monotonicity of the method: $F(x_k) \leq F(x_{k - 1})$
and the bound: $F^{*} \leq F(x_k)$.
\qed

\section{Proof of Theorem~\ref{th:LocalAlg2}}
\proof
Let us plug $y = x^{*}$ into~\eqref{OneStep}.
Thus, we obtain, for every $k \geq 1$:
$$
\ba{rcl}
F(x_{k + 1}) & \leq &
F^{*} + \frac{L_p\|x_k - x^{*}\|^{p + 1}}{p!} + \delta_{k + 1} \\
\\
& \leq & F^{*} + \frac{L_p}{p!}\bigl( \frac{2}{\sigma_2} \bigr)^{p + 1 \over 2}
(F(x_k) - F^{*})^{p + 1 \over 2} + \delta_{k + 1} \\
\\
& = & 
F^{*} + \frac{L_p}{p!}\bigl( \frac{2}{\sigma_2} \bigr)^{p + 1 \over 2}(F(x_k) - F^{*})^{p + 1 \over 2}
+ c(F(x_{k - 1}) - F(x_{k}))^{p + 1 \over 2} \\
\\
& \leq &
F^{*} + \Bigl( \frac{L_p}{p!}\bigl( \frac{2}{\mu} \bigr)^{p + 1 \over 2} + c \Bigr)
(F(x_{k - 1}) - F^{*})^{p + 1 \over 2},
\ea
$$
where monotonicity of the method: $F(x_k) \leq F(x_{k - 1})$
and the bound: $F^{*} \leq F(x_k)$ are used in the last inequality.
\qed

\section{Proof of Theorem~\ref{th:GlobalConvAlg3}}
\proof
The proof is similar to that one of Theorem~\ref{th:GlobalConvAlg1}.

By Lemma~\ref{lemma:Global}, for every $y \in \dom \psi$ we have
$$
\ba{rcl}
F(x_{k + 1}) & \refLE{OneStep} & F(y) + \frac{L_p\|y - y_{k}\|^{p + 1}}{p!} + \delta_{k + 1},  
\quad k \geq 0.
\ea
$$
Let us substitute $y = \lambda_k x_k + (1 - \lambda_k) x^{*}$, 
with $\lambda_k$ defined in the algorithm:
$$
\ba{rcl}
\lambda_k & \equiv & \bigl(\frac{k}{k + 1}\bigr)^{p + 1}.
\ea
$$
Thus we obtain
$$
\ba{rcl}
F(x_{k + 1}) & \leq &
(1 - \lambda_k) F^{*} + 
\lambda_k F(x_k) + (1 - \lambda_k)^{p + 1}\frac{L_p\|x_0 - x^{*}\|^{p + 1}}{p!} + \delta_{k + 1},
\ea
$$
or, equivalently
$$
\ba{rcl}
A_{k + 1} (F(x_{k + 1}) - F^{*}) & \leq & A_k ( F(x_k) - F^{*})
+ \frac{a_{k + 1}^{p + 1}}{A_{k + 1}^p} 
\frac{L_p \|x_0 - x^{*}\|^{p + 1}}{p!} + A_{k + 1} \delta_{k + 1},
\ea
$$
where $A_{k} \equiv k^{p + 1}$ and $a_k \equiv A_k - A_{k - 1}$ 
(so it holds: $\lambda_k \equiv {A_k} / {A_{k + 1}}$ and
$1 - \lambda_k \equiv {a_{k + 1}} / {A_{k + 1}}$).
Telescoping these inequalities and 
using the bounds~\eqref{th1-ratio-bound} and~\eqref{th1-delta-sum}
complete the proof.
\qed

\section{Proof of Theorem~\ref{th:GlobalAlg4}}
\proof
The proof is similar to that one of 
Theorem~1 from~\cite{doikov2019contracting},
where convergence rate of Contracting Proximal Method is established.
Additional technical difficulties, which are arising here, are caused by
using inexact solution of the subproblem, equipped with the stopping condition~\eqref{StopCondition}.

We denote the optimal point of $h_{k + 1}(\cdot)$ by $z_{k + 1} \Def \argmin_{y \in \E} h_{k+1}(y)$.
Since the next prox-center $v_{k + 1}$ 
is defined as an approximate minimizer, we have
\beq \label{eq:th6hmin}
\ba{rcl}
h_{k + 1}(v_{k + 1}) - h_{k + 1}(z_{k + 1}) & \leq & \zeta_{k + 1}.
\ea
\eeq
Function $h_{k + 1}(\cdot)$ is strongly convex with respect to $d(\cdot)$,
thus we have
\beq \label{eq:th6PointsBound}
\ba{rcl}
\zeta_{k + 1} & \refGE{eq:th6hmin} & h_{k + 1}(v_{k + 1}) - h_{k + 1}(z_{k + 1})
\;\; \geq \;\; \beta_d(z_{k + 1}; v_{k + 1}) \\
\\
& \geq & \frac{1}{2^{p - 1}(p + 1)} \|v_{k + 1} - z_{k + 1}\|^{p + 1}.
\ea
\eeq
Therefore,
\beq \label{eq:th6xiDef}
\ba{rcl}
\|v_{k + 1} - z_{k + 1} \| & \refLE{eq:th6PointsBound} & \xi_{k + 1} 
\;\; \Def \;\; 2^{p - 1 \over p + 1}(p + 1)^{1 \over p + 1} \zeta_{k + 1}^{1 \over p + 1}.
\ea
\eeq

Let us prove by induction the following inequality, for every $k \geq 0$:
\beq \label{eq:th6inductive}
\ba{rcl}
\beta_d(x_0; x) + A_k F(x) & \geq & \beta_d(v_k; x) + A_k F(x_k) + C_k(x),
\qquad x \in \dom F.
\ea
\eeq
where $C_k(x) \Def -\sum\limits_{i = 1}^k\bigl( \tau_i \|x - v_i\| + \zeta_i \bigr)$,
and $\tau_i \Def p2^{p - 2}\|z_{i} - x_0\|^{p - 1} \xi_{i} + 2^{p - 2} \xi_{i}^p$.

It obviously holds for $k = 0$. Assume that it holds
for the current iterate, and consider the next step $k + 1$:
\beq \label{eq:th6induction}
\ba{rcl}
& \; & 
\!\!\!\!\!\!\!\!\!\!\!\!\!\!\!\!\!\!  
\beta_d(x_0; x) + A_{k + 1} F(x) \\
\\
& = & 
\beta_d(x_0; x) + A_k F(x) + a_{k + 1} F(x) \\
\\
& \refGE{eq:th6inductive} &
\beta_d(v_k; x) + A_k F(x_k) + a_{k + 1} F(x) + C_k(x) \\
\\
& \geq &
\beta_d(v_k; x) + A_{k + 1} f\bigl( \frac{a_{k + 1} x + A_k x_k}{A_{k + 1}} \bigr)
+ a_{k + 1}\psi(x) + A_{k} \psi(x_k) + C_k(x) \\
\\
& = & h_{k + 1}(x) + A_k \psi(x_k) + C_k(x),
\ea
\eeq
where the last inequality holds by convexity of $f$.

Using strong convexity of $h_{k + 1}(\cdot)$ with respect to $d(\cdot)$, we obtain
\beq \label{eq:th6hStrongConv}
\ba{rcl}
h_{k + 1}(x) & \geq & h_{k + 1}(z_{k + 1}) + \beta_d(z_{k + 1}; x) \\
\\
& \refGE{eq:th6hmin} &
h_{k + 1}(v_{k + 1}) + \beta_d(z_{k + 1}; x) - \zeta_{k + 1} \\
\\
& = &
h_{k + 1}(v_{k + 1}) + \beta_d(v_{k + 1}; x)
+ \beta_d(z_{k + 1}; v_{k + 1})
+ \la \nabla d(v_{k + 1}) - \nabla d(z_{k + 1}), x - v_{k + 1} \ra - \zeta_{k + 1} \\
\\
& \geq & 
h_{k + 1}(v_{k + 1}) + \beta_d(v_{k + 1}; x)
- \| \nabla d(v_{k + 1}) - \nabla d(z_{k + 1})\|_{*} \cdot \|x - v_{k + 1}\| - \zeta_{k + 1}.
\ea
\eeq

Now, computing second derivative of $d(x) = \frac{1}{p + 1}\|x - x_0\|^{p + 1}$, 
we get
\beq \label{eq:th6dSecondDeriv}
\ba{rcl}
\nabla^2 d(x)
& = & (p - 1)\|x - x_0\|^{p - 3} B(x - x_0) (x - x_0)^{*}B + \|x - x_0\|^{p - 1}B \\
\\
& \preceq & p\|x - x_0\|^{p - 1}B.
\ea
\eeq
Therefore, 
\beq \label{eq:th6LipGrad}
\ba{rcl}
\| \nabla d(v_{k + 1}) - \nabla d(z_{k + 1}) \|_{*}
& = &
\| \int\limits_0^1 \nabla^2 d(z_{k + 1} + \tau(v_{k + 1} - z_{k + 1})) d\tau
(v_{k + 1} - z_{k + 1})  \|_{*} \\
\\
& \refLE{eq:th6xiDef} &
\xi_{k + 1} 
\int\limits_0^1 \| \nabla^2 d(z_{k + 1} + \tau(v_{k + 1} - z_{k + 1})) \| d\tau \\
\\
& \refLE{eq:th6dSecondDeriv} &
p  \xi_{k + 1} 
\int\limits_0^1 \|z_{k + 1} - x_0 + \tau(v_{k + 1} - z_{k + 1}) \|^{p - 1} d\tau \\
\\
& \overset{\eqref{PoweredSumBound},\eqref{eq:th6xiDef}}{\leq} & 
p  \xi_{k + 1} 
\int\limits_0^1 \bigl( 2^{p - 2}\|z_{k + 1} - x_0\|^{p - 1}
+ 2^{p - 2} \tau^{p - 1} \xi_{k + 1}^{p - 1} \bigr) d\tau \\
\\
& = & 
p 2^{p - 2} \|z_{k + 1} - x_0\|^{p - 1} \xi_{k + 1}
+ 2^{p - 2} \xi_{k + 1}^p \;\; \Def \;\; \tau_{k + 1}.
\ea
\eeq
Combining obtained bounds together, we conclude
$$
\ba{rcl}
\beta_d(x_0; x) + A_{k + 1}F(x) 
& \refGE{eq:th6induction} &
h_{k + 1}(x) + A_k \psi(x_k) + C_k(x) \\
\\
& \overset{\eqref{eq:th6hStrongConv},\eqref{eq:th6LipGrad}}{\geq} & 
h_{k + 1}(v_{k + 1}) + \beta_d(v_{k + 1}; x)
- \tau_{k + 1} \|x - v_{k + 1}\| - \zeta_{k + 1} + A_k \psi(x_k) + C_k(x) \\
\\
& = & 
h_{k + 1}(v_{k + 1}) + \beta_d(v_{k + 1}; x) + A_k \psi(x_k) + C_{k + 1}(x) \\
\\
& = &
A_{k + 1}f(x_{k + 1}) + a_{k + 1}\psi(v_{k + 1}) + \beta_d(v_k; v_{k + 1})
+ \beta_d(v_{k + 1}; x) + A_k \psi(x_k) + C_{k + 1}(x) \\
\\
& \geq &
A_{k + 1}F(x_{k + 1}) + \beta_d(v_{k + 1}; x) + C_{k + 1}(x).
\ea
$$
Thus,~\eqref{eq:th6inductive} is proven for all $k \geq 0$.

Let us plug $x := x^{*}$ into~\eqref{eq:th6inductive}. We obtain
\beq \label{eq:th6Main}
\ba{rcl}
\beta_d(v_k; x^{*}) + A_k(F(x_k) - F^{*}) 
& \refLE{eq:th6inductive} &
\beta_d(x_0; x^{*}) - C_k(x^{*}) \\
\\
& = &
\beta_d(x_0; x^{*}) + \sum\limits_{i = 1}^k \zeta_i
+ \sum\limits_{i = 1}^k \tau_i \|v_i - x^{*}\| \\
\\
& \refLE{SublinBound} &
\beta_d(x_0; x^{*}) + \frac{c(p + 2)}{p + 1} 
+ \sum\limits_{i = 1}^k \tau_i \|v_i - x^{*}\|
\;\; \Def \;\; \alpha_k,
\ea
\eeq
and to finish the proof, we need to estimate $\alpha_k$ from above.

By uniform convexity of $d(\cdot)$, we have
\beq \label{eq:th6AlphaLower}
\ba{rcl}
\frac{1}{2^{p - 1}(p + 1)}\|v_k - x^{*}\|^{p + 1}
& \leq & \beta_d(v_k; x^{*}) \;\; \refLE{eq:th6Main} \;\; \alpha_k.
\ea
\eeq
At the same time,
$$
\ba{rcl}
\alpha_k & = & \alpha_{k - 1} + \tau_{k}\|v_{k} - x^{*}\|
\;\; \refLE{eq:th6AlphaLower} \;\;
\alpha_{k - 1} 
+ 
2^{p - 1 \over p + 1} (p + 1)^{1 \over p + 1} \tau_k \alpha_k^{1 \over p + 1}.
\ea
$$
Dividing both sides of the last inequality by $\alpha_k^{1 \over p + 1}$,
and using monotonicity of $\{ \alpha_k \}_{k \geq 0}$, we get
$$
\ba{rcl}
\alpha_k^{p \over p + 1} & \leq & \alpha_{k - 1}^{p \over p + 1}
+ 2^{p - 1 \over p + 1} (p + 1)^{1 \over p + 1} \tau_k.
\ea
$$
Therefore,
\beq \label{eq:th6AlphaKBound}
\ba{rcl}
\alpha_k & \leq &
\left( 
\alpha_0^{p \over p + 1}
+ 2^{p - 1 \over p + 1}(p + 1)^{1 \over p + 1} \sum\limits_{i = 1}^k \tau_i
\right)^{p + 1 \over p}.
\ea
\eeq
To finish, it remains to bound the sum of $\tau_i$, which is
\beq \label{eq:th6TauIBound}
\ba{rcl}
\sum\limits_{i = 1}^k \tau_i
& = &
\sum\limits_{i = 1}^k
\left(
p2^{p - 2}\|z_i - x_0\|^{p - 1} \xi_i + 2^{p - 2} \xi_i^p 
\right) \\
\\
& \refLE{PoweredSumBound} &
\sum\limits_{i = 1}^k 
\left(
p4^{p - 2}\|z_i - x^{*}\|^{p - 1} \xi_i
+ p 4^{p - 2} \|x_0 - x^{*}\|^{p - 1} \xi_i + 2^{p - 2}\xi_i^p
\right) \\
\\
& \refEQ{eq:th6xiDef} &
p 4^{p - 2} \sum\limits_{i = 1}^k
\|z_i - x^{*}\|^{p - 1} \xi_i
\;\; + \;\;
p (p + 1)^{1 \over p + 1} 2^{p - 1 \over p + 1} 4^{p - 2} \|x_0 - x^{*}\|^{p - 1}
\sum\limits_{i = 1}^k \zeta_i^{1 \over p + 1} \\
\\
& \; & \quad
+
\quad
2^{p - 2}2^{(p - 1)p \over p + 1} (p + 1)^{p \over p + 1}
\sum\limits_{i = 1}^k \zeta_i^{p \over p + 1} \\
\\
& \refLE{SublinBound} &
p 4^{p - 2} \sum\limits_{i = 1}^k
\|z_i - x^{*}\|^{p - 1} \xi_i + \Delta_1,
\ea
\eeq
where
$$
\ba{rcl}
\Delta_1 & \Def &
p(p + 1)^{1 \over p + 1} 2^{p - 1 \over p + 1} 4^{p - 2} \|x_0 - x^{*}\|^{p -  1}
c^{1 \over p + 1} (p + 2)
\;\; + \;\;
2^{p - 2} 2^{(p - 1)p \over p + 1} (p + 1)^{p \over p + 1}
c^{p \over p + 1} \frac{(p + 2)p}{(p + 2)p - p - 1},
\ea
$$
and we need to bound $\|z_i - x^{*}\|$ from above.

Substituting $x := x^{*}$ into~\eqref{eq:th6induction}, we have
\beq \label{eq:th6Almost}
\ba{rcl}
\beta_d(x_0; x^{*}) + A_{k + 1}F^{*} 
& \refGE{eq:th6induction} &
h_{k + 1}(x^{*}) + A_k \psi(x_k) + C_k(x^{*}) \\
\\
& \geq &
h_{k + 1}(z_{k + 1}) + \beta_d(z_{k + 1}; x^{*}) + A_k \psi(x_k) + C_k(x^{*}) \\
\\
& \geq & 
A_{k + 1}F\bigl( \frac{a_{k + 1}z_{k + 1} + A_k x_k}{A_{k + 1}} \bigr)
+ \beta_d(z_{k + 1}; x^{*}) + C_k(x^{*}).
\ea
\eeq
So,
\beq \label{eq:th6GammaKBound}
\ba{rcl}
& \; & 
\!\!\!\!\!\!\!\!\!\!\!\!\!\!\!\!\!\! 
\frac{1}{2^{p - 1}(p + 1)}\|z_{k + 1} - x^{*}\|^{p + 1} \\
\\
& \leq & \beta_d(z_{k + 1}; x^{*}) 
\;\; \refLE{eq:th6Almost} \;\;
\beta_d(x_0; x^{*}) - C_k(x^{*}) 
\;\; \leq \;\; \alpha_k \\
\\
& \overset{\eqref{eq:th6AlphaKBound},\eqref{eq:th6TauIBound}}{\leq} &
\left(
\Delta_2 + \Delta_3 \sum\limits_{i = 1}^k\|z_i - x^{*}\|^{p - 1} \xi_i
\right)^{p + 1 \over p},
\ea
\eeq
with
$$
\ba{rcl}
\Delta_2 & \Def & 
\alpha_0^{p \over p + 1} 
+ 2^{p - 1 \over p + 1} (p + 1)^{1 \over p + 1} \Delta_1,
\quad \text{and} \quad
\Delta_3 \;\; \Def \;\; 2^{p - 1 \over p + 1} (p + 1)^{1 \over p + 1} p 4^{p - 2}.
\ea
$$
Therefore, for the monotone sequence 
$\gamma_k \; \Def \; \Delta_2 + \Delta_3 \sum\limits_{i = 1}^k \|z_i - x^{*}\|^{p - 1} \xi_i$,
it holds
$$
\ba{rcl}
\gamma_{k + 1} & = & \gamma_k 
+ \Delta_3 \|z_{k + 1} - x^{*}\|^{p - 1}\xi_{k + 1}
\;\; \refLE{eq:th6GammaKBound} \;\;
\gamma_k + \Delta_3 2^{(p - 1)^2 \over p + 1} (p + 1)^{p - 1 \over p + 1} 
\gamma_{k}^{p - 1 \over p} \xi_{k + 1} \\
\\
& \leq &
\gamma_k + \Delta_3 2^{(p - 1)^2 \over p + 1} (p + 1)^{p - 1 \over p + 1} 
\gamma_{k + 1}^{p - 1 \over p} \xi_{k + 1}.
\ea
$$
Dividing both sides by $\gamma_{k + 1}^{p - 1 \over p}$, 
and using monotonicity again, we obtain
\beq \label{eq:th6GammaOneDiff}
\ba{rcl}
\gamma_{k + 1}^{1 \over p}
& \leq &
\gamma_k^{1 \over p} + \Delta_3 2^{(p - 1)^2 \over p + 1} (p + 1)^{p - 1 \over p + 1} \xi_{k + 1},
\ea
\eeq
Telescoping which, gives
\beq \label{eq:th6GammaTelescoped}
\ba{rcl}
\alpha_k^{\frac{1}{p + 1}} & \refLE{eq:th6GammaKBound} &  \gamma_k^{1 \over p} 
\;\; \refLE{eq:th6GammaOneDiff} \;\; 
\gamma_0^{1 \over p}
+ \Delta_3 2^{(p - 1)^2 \over p + 1} (p + 1)^{p - 1 \over p + 1} \sum\limits_{i = 1}^k \xi_i \\
\\
& \leq &
\Delta_2^{1 \over p}
+ \Delta_3 2^{(p - 1)p \over p + 1} (p + 1)^{p \over p + 1} 
c^{1 \over p + 1} \frac{s}{s - p - 1} \\
\\
& = &
O\Bigl( \|x_0 - x^{*}\| + c^\frac{1}{p + 1} \Bigr).
\ea
\eeq
Finally,
$$
\ba{rcl}
F(x_k) - F^{*} & \refLE{eq:th6Main} & \frac{\alpha_k}{A_k} 
\;\; \refLE{eq:th6GammaTelescoped} \;\;
O\Bigl(  \frac{L_p( \|x_0 - x^{*}\|^{p + 1} + c)}{k^{p + 1}}  \Bigr).
\ea
$$

Lastly, let us prove bound~\eqref{eq-NStepsBound} for the number of 
tensor steps, needed to find $v_{k + 1}$. We minimize $h_{k + 1}(\cdot)$,
starting from the previous prox-point $v_k$.
We denote the first component of $h_{k + 1}(\cdot)$, by:
$$
\ba{rcl}
g_{k + 1}(x) & \Def & A_{k + 1} f\bigl( \frac{a_{k + 1}x + A_k x_k}{A_{k+1}} \bigr),
\ea
$$
which is \textit{contracted} version of the smooth part of our objective $F(x)$.
Direct computation gives the following relation between Lipschitz constants 
for the derivatives of $g_{k + 1}$ and $f$:
\beq \label{eq:th6LipBound}
\ba{rcl}
L_p(g_{k + 1}) & = & \frac{a_{k + 1}^{p + 1}}{A_{k + 1}^p} L_p(f)
\;\; = \;\; 
\frac{ ((k + 1)^{p + 1} - k^{p + 1})^{p + 1} }{(k + 1)^{p(p + 1)}}
\;\; \leq \;\;
\frac{( (p + 1)(k + 1)^p)^{p + 1}}{(k + 1)^{p(p + 1)}}
\;\; = \;\;
(p + 1)^{p + 1}.
\ea
\eeq
Therefore, condition number $\omega_p$~\eqref{eq-CondNumber} for $h_{k + 1}$ 
is bounded by an \textit{absolute constant}, and we need to estimate only the value
under the logarithm in~\eqref{LinearComplexity}.
Due to Lemma~\ref{lemma:Global}, one monotone inexact tensor step $M := M_{H,\delta}(v_k)$
for function $h_{k+1}(\cdot)$ with constant $H := pL_p(g_{k + 1})$ gives
\beq \label{eq:th6hGlobal}
\ba{rcl}
h_{k + 1}(M) 
& \overset{\eqref{OneStep},\eqref{eq:th6LipBound}}{\leq} & 
h_{k + 1}(y) + \frac{(p + 1)^{p + 1}\|y - v_k\|^{p + 1}}{p!}
+ \delta, \qquad y \in \dom \psi.
\ea
\eeq
We substitute $y := x^{*}$ (minimizer of $F$) into~\eqref{eq:th6hGlobal}, and thus we obtain
$$
\ba{rcl}
h_{k + 1}(M) - h_{k + 1}^{*}
& \refLE{eq:th6hGlobal} &
h_{k + 1}(x^{*}) - h_{k + 1}^{*} 
+ \frac{(p + 1)^{p + 1}\|v_k - x^{*}\|^{p + 1}}{p!} + \delta \\
\\
& \refLE{eq:th6induction} &
A_{k + 1}F^{*} - A_k \psi(x_k) + \beta_d(x_0; x^{*}) 
- C_k(x^{*}) - h_{k + 1}^{*} 
+ \frac{(p + 1)^{p + 1}\|v_k - x^{*}\|^{p + 1}}{p!} + \delta \\
\\
& \overset{\eqref{eq:th6Main},\eqref{eq:th6AlphaLower}}{\leq} &
A_{k + 1}F^{*} - A_k \psi(x_k) - h_{k + 1}^{*}
+ \Bigl(1 + \frac{(p + 1)^{p + 2}2^{p - 1}}{p!}\Bigr) \alpha_k + \delta \\
\\
& = &
A_{k + 1} F^{*} - \min\limits_{y} \Bigl\{ h_{k + 1}(y) + A_k \psi(x_k) \Bigr\}
+ \Bigl(1 + \frac{(p + 1)^{p + 2}2^{p - 1}}{p!}\Bigr) \alpha_k + \delta \\
\\
& \leq &
A_{k + 1} F^{*} - \min\limits_{y} \Bigl\{ 
A_{k + 1}F\bigl( \frac{a_{k + 1}y + A_k x_k}{A_{k + 1}}\bigr) + \beta_d(v_k; y)   
 \Bigr\}
+ \Bigl(1 + \frac{(p + 1)^{p + 2}2^{p - 1}}{p!}\Bigr) \alpha_k + \delta \\
\\
& \leq & 
A_{k + 1} F^{*} - \min\limits_{y} 
\Bigl\{ A_{k + 1} F\bigl( \frac{a_{k + 1}y + A_k x_k}{A_{k + 1}}\bigr) \Bigr\}
+ \Bigl(1 + \frac{(p + 1)^{p + 2}2^{p - 1}}{p!}\Bigr) \alpha_k + \delta
\\
\\
& = &
\Bigl(1 + \frac{(p + 1)^{p + 2}2^{p - 1}}{p!}\Bigr) \alpha_k + \delta 
\quad \refLE{eq:th6GammaTelescoped} \quad
O\Bigl( \|x_0 - x^{*}\|^{p + 1} + c + \delta \Bigr).
\ea
$$

So, if we set $\delta := c$ and 
perform just one step of the monotone inexact tensor method for $h_{k + 1}(\cdot)$,
the remaining amount of steps $N_k$ needed to find $v_{k+1}$, 
such that~\eqref{StopCondition} holds, is bounded as:
$$
\ba{rcl}
N_k & \refLE{LinearComplexity} &
O\Bigl( \log \frac{h_{k + 1}(M) - h_{k + 1}^{*}}{\zeta_{k + 1}}  \Bigr)
\;\; \leq \;\;
O\Bigl( \log \frac{k(\|x_0 - x^{*}\|^{p + 1} + c)}{c} \Bigr).
\ea
$$
\qed

\end{document}